\documentclass{article}
\usepackage[top=3cm,bottom=3cm]{geometry}
\usepackage[mathscr]{euscript}
\usepackage{amssymb,amsthm,amsmath,amsxtra,dsfont,appendix,cancel,tensor}
\usepackage{mathrsfs}
\usepackage[nottoc]{tocbibind}
\usepackage{graphicx,tikz,esint}
\usetikzlibrary{arrows}
\usepackage{hyperref} \hypersetup{colorlinks=true, citecolor=blue, linkcolor=black}

\numberwithin{equation}{section} 

\renewcommand{\limsup}[1]{\underset{{#1}}{\overline{\operatorname{lim}}}}

\newcommand{\supp}{\operatorname{supp}}
\newcommand{\tr}{\operatorname{Tr}}
\newcommand{\Var}{\operatorname{Var}}
\renewcommand{\P}[1]{\mathbb{P}\left[#1\right]}
\newcommand{\E}[1]{\mathbb{E}\left[#1\right]}

\newcommand{\1}{\mathds{1}}
\newcommand{\A}{\mathrm{A}}
\newcommand{\C}{\mathbb{C}}

\newcommand{\Cu}{\operatorname{C}}

\newcommand{\F}{\mathrm{F}_V}
\newcommand{\G}{\mathrm{G}}
\newcommand{\g}{\mathcal{G}}

\newcommand{\I}{\mathbf{I}}
\newcommand{\J}{\mathbf{J}}  
\renewcommand{\k}{\bf k}%{\underline{k}}    
\renewcommand{\L}{\mathbf{L}}
\newcommand{\m}{M}
\newcommand{\M}{\mathbf{M}}

\newcommand{\N}{\mathbb{N}}
\newcommand{\No}{\mathcal{N}}
\renewcommand{\o}{\mathrm{r}}
\newcommand{\p}{\mathbf{P}}
\newcommand{\R}{\mathbb{R}}

\newcommand{\var}{\operatorname{Var}}
\newcommand{\Sy}{\mathbb{S}}
\newcommand{\T}{\mathbf{T}} 
\newcommand{\U}{\mathcal{U}}  
\newcommand{\W}{\mathrm{w}} 
\newcommand{\X}{\mathrm{X}} 
\newcommand{\Z}{\mathbb{Z}}

\newtheorem{theorem}{Theorem}[section]
\newtheorem{definition}[theorem]{Definition}
\newtheorem{proposition}[theorem]{Proposition}
\newtheorem{corollary}[theorem]{Corollary}
\newtheorem{lemma}[theorem]{Lemma}

\theoremstyle{definition} \newtheorem{remark}{Remark}[section]

\title{ \vspace{-.5cm} 
Limit theorems for biorthogonal ensembles and related combinatorial identities}
\date{}
\author{Gaultier Lambert\thanks{KTH Royal Institute of Technology, Department of Mathematics, glambert@kth.se. Supported by the grant KAW 2010.0063 from the Knut and Alice Wallenberg Foundation. }}

\begin{document}

\maketitle

\begin{abstract} \normalsize

\noindent We study the fluctuations of certain biorthogonal ensembles for which the underlying family $\{P,Q\}$ satisfies a finite-term recurrence relation of the form $x P(x) = \mathbf{J}P(x)$. 
For polynomial linear statistics of such ensembles, we reformulate the cumulant method introduced in~\cite{Soshnikov_00a} in terms of counting certain lattice paths on  the adjacency graph of the recurrence matrix $\mathbf{J}$. 
In the spirit of~\cite{BD_13}, we show that the asymptotic fluctuations  are described by the right-limits of the  matrix $\mathbf{J}$. Moreover, whenever the right-limit is a Laurent matrix, we show that the CLT is equivalent to  Soshnikov's main combinatorial lemma. We discuss several applications to unitary invariant Hermitian random matrices. In particular, we provide a general central limit theorem (CLT) and a law of large numbers in the one-cut regime.  We also prove a CLT for the square singular values of the product of independent complex rectangular Ginibre matrices, as well as for  the Laguerre and Jacobi biorthogonal ensembles introduced in~\cite{Borodin_99}, and we explain how to recover the equilibrium measure from the asymptotics of the recurrence coefficients.
Finally, we discuss the connection with the Strong Szeg\H{o} limit theorem where this combinatorial method originates.\\
\end{abstract}

\vspace{.4cm}

{\bf Keywords:} biorthogonal ensembles; cumulant method; central limit theorems; law of large numbers; strong Szeg\H{o} theorem.  \\

{\bf Mathematics Subject Classification.} 60B20, 60G55, 60F05, 82B41

\vspace{.4cm}

%\tableofcontents

\section{Introduction}
%
%%% Matrix models
%
Random matrix ensembles were introduced by Wigner who proposed to describe  the distribution of scattering resonances of heavy nuclei by the eigenvalues of certain random matrices, \cite{Wigner_51}. Now-days, random matrix models  and their generalizations have found a wide range of applications in physics and mathematics, \cite{OHRMT}. Mathematically, unitary invariant ensembles are easier to analyze because of their determinantal structure. The connection between these ensembles and orthogonal polynomials  was developed  by Gaudin and Metha in the early stages of the theory, \cite{MG_60}. Namely, if the manifold of $N\times N$ Hermitian matrices, denoted $\mathcal{H}_N$, is equipped with the probability measure
\begin{equation} \label{HM}
d \mathbb{P}_{N,V}(M) =Z_{N,V}^{-1}   e^{-N \tr V(M)} dM ,
\end{equation}
where $dM$ is the Lebesgue measure on $\mathcal{H}_N$ and $V:\R\to\R$  is a  function which satisfies the condition
\begin{equation}\label{confinement} 
 \lim_{|x| \to\infty} \frac{V(x)}{\log(1+|x|^2)} = +\infty . 
 \end{equation}
The eigenvalues of a random matrix sampled from $\mathbb{P}_{N,V}$ have a joint density on $\R^N$ with respect to the product of $N$ copies of the measure $d\mu_N/dx =  e^{-N V(x)} $ which is given by
 \begin{equation} \label{jpdf}
\varrho_N(x_1,\dots,x_N) = \frac{1}{N!}  \det_{N\times N}\big[  K_N(x_i,y_j) \big] .
\end{equation}
The function $K_N$ is expressed in terms of the orthogonal polynomials  $(P_k^N)_{k=0}^\infty$  with respect to the measure $\mu_N$ as follows,
\begin{equation}\label{kernel}
 K_N(x,y)= \sum_{k=0}^{N-1} P^N_k(x)P^N_k(y) 
\end{equation}
and, by convention: 
\begin{equation} \label{OP}
 P^N_k(x)=\kappa_{k}^N x^k + \cdots  \hspace{1cm}\text{and}\hspace{1cm} 
\int P_k^N(x) P_j^N(x) d\mu_N(x)= \delta_{k,j}  ,
\end{equation}
where  $\kappa_k^N>0$. The most well known example is the Gaussian Unitary Ensemble (GUE) which corresponds to the potential $V(x)=2x^2$, in which case $(P_k^N)_{k=0}^\infty$ are rescaled Hermite polynomials. 
%
%%%  Biorthogonal determinantal processes
%
Generally, if $\mathfrak{X}$ is a complete separable metric space equipped with a Radon measure $\mu$, a determinantal process is a point process on  $(\mathfrak{X},\mu)$ whose correlation functions with respect to $\mu$ are of the form 
 \begin{equation} \label{correlation}
\rho_k(x_1,\dots,x_k) =  \det_{k\times k}\big[  K(x_i,y_j) \big]  ,
\end{equation}
where the function $K$ is called the correlation kernel. There is an extensive literature available on determinantal processes, e.g$.$ \cite{Soshnikov_00b,HKPV_06, Johansson_05}, and the determinantal form of the correlation functions implies  that many natural observables  can be expressed explicitly in terms of the kernel $K$,  for instance the moment generating function of a linear statistic, see formula (\ref{Laplace_3}) below. 
Using the orthogonality relations, it is simple to check by integrating $N-k$ variables that the correlation functions of  the j.p.d.f$.$ (\ref{jpdf}) satisfies (\ref{correlation}) for all $k\le N$. Hence, the eigenvalues of the Hermitian ensemble $\mathbb{P}_{N,V}$ form a determinantal process with correlation kernel (\ref{kernel}) with respect to $\mu_N$. In this article, we will consider a generalizations of these ensembles introduced by Borodin in \cite{Borodin_99} which are  called {\it Biorthogonal Ensembles} (BOE).
%
%% Biorthogonal ensembles
%
\begin{definition}\label{ensemble}
A biorthogonal family $\{ P_k, Q_k\}_{k=0}^\infty$ is a set of measurable functions defined on $(\mathfrak{X},\mu)$ such that $P_kQ_n \in L^1(\mu)$ for all $n,k \in \N_0$ and  
\begin{equation} \label{biorthogonal}
\int P_k(x) Q_n(x) d\mu(x)= \delta_{k,n}  .
\end{equation}
A biorthogonal ensemble $\Xi_N=(\xi_k)_{k=1}^N$ is a determinantal process on $(\mathfrak{X},\mu_N)$ with correlation kernel  
\begin{equation}\label{kernel_0}
 K_N(x,y)= \sum_{k=0}^{N-1} P^N_k(x)Q^N_k(y)  , 
\end{equation}
where $\{ P_k^N, Q_k^N\}_{k=0}^\infty$ is a biorthogonal family. 
\end{definition}

In the following, we will only consider the case $\mathfrak{X}=\R$, although the results could be formulated in a broader context. %An alternative definition of a biorthogonal ensemble is that it is a point process with the j.p.d.f$.$ (\ref{biorthogonal_1}). %see section~\ref{sect:cumulants}.
Interestingly, biorthogonal ensembles arise in many different contexts beyond Hermitian ensembles, e.g.  non-intersecting paths, domino tilings, etc. We do not intend to review the theory here and we refer to \cite{Konig_05,Borodin_11,BD_13, Kuijlaars_10} for further discussion.  
In many examples, there is more structure than just the biorthogonality relation (\ref{biorthogonal}). For instance, for the Hermitian models (\ref{HM}), the functions $Q_k^N=P_k^N$ are polynomials of degree $k$ and the kernel (\ref{kernel}) is symmetric.  In this case, the determinantal  process $\Xi_N$ is entirely characterized by the reference measure  $d\mu_N/dx =  e^{-N V(x)}$  and $K_N$ is called the {\it Christoffel-Darboux kernel}.  In the following, we will consider the general setting where $\mu_N$ is an arbitrary Borel measure on $\R$ which satisfies, for all $k\in\N$,
\begin{equation} \label{reference}
 \int_\R |x|^k d\mu_N(x) <\infty  .
 \end{equation}

This condition guarantees that the polynomials (\ref{OP}) exist and $(P_k^N)_{k=0}^\infty$ is an orthonormal basis of the space of polynomials $\R\langle x\rangle$ equipped with the inner product inherited from $L^2(\mu_N)$.
In the literature, such a point process is usually called the {\it Orthogonal Polynomial Ensemble} (OPE) with reference measure $\mu_N$ and it is of interest to know which properties of the measure $\mu_N$ lead  to  $\Xi_N$ having certain universal features as $N\to\infty$.  For instance, one may consider the well-known problem of finding for which potentials $V(x)$ the Hermitian model (\ref{HM}) falls in the GUE universality class. In particular, we will revisit the question of universality of fluctuations  for a linear statistic
\begin{equation} \label{statistics} 
\Xi_N(f)= \sum_{k=1}^N f(\xi_k) ,
\end{equation} 
where $\{ \xi_k\}_{k=1}^N$ refers to the point configuration of a biorthogonal ensemble. 
%
%% Breuer-Duits theory
%
Our approach is motivated by~\cite{BD_13} where the authors prove that it is possible to express the moment generating function of the random variable $\Xi_N(f)$ as  a Fredholm determinant involving the Jacobi matrix $\J$ of the measure $\mu_N$. 
 If $F$ is a polynomial and $|\lambda|$ is sufficiently small,
\begin{equation} \label{Laplace}
 \E{e^{ \lambda \Xi_N(F)}} =  \det\big[ \I +  \p_N \big( \exp\lambda F(\J) -1 \big) \p_N \big]  ,
 \end{equation}  
 where $\p_N$ is the projection from $l^2(\mathbb{N}_0)$ onto $\text{span}(e_0,\dots, e_{N-1})$.
Their insight was that the concepts of right-limit and Laurent matrices  (definitions~\ref{R_limit} and~\ref{Laurent} respectively)  establish the link with the theory of Toeplitz determinants.
%
%% Combinatorial formalism
%
Our  approach is different and consists in expressing the cumulants of the random variable $\Xi_N(F)$  in terms of lattice-paths on the graph of the Jacobi matrix $\J$; see lemma~\ref{thm:cumulants}.
Given a test function $f$, the cumulants of a  linear statistic $\Xi_N(f)$  are defined by the  power series
 \begin{equation} \label{Laplace_0}
\log \E{e^{ \lambda \Xi_N(f)}} =   \sum_{n=1}^\infty  \Cu^n_N\left[ f\right] \frac{\lambda^n}{n!}  .
\end{equation}    

For determinantal processes, the cumulant method was introduced by Soshnikov in~\cite{Soshnikov_00a} to prove CLTs for linear statistics of eigenvalues of random matrices from  the classical compact groups.  For Haar distributed unitary matrices  (CUE), using Fourier analysis, Soshnikov reduced the problem of computing the limits of cumulants   to certain combinatorial identities known as the {\it Main Combinatorial Lemma} (MCL - theorem~\ref{thm:G}). 
In this paper,  we apply a similar method to unitary invariant random matrices and, in general, to any  BOE whose functions 
$(P^N_k)_{k=0}^\infty$ satisfy a {\it recurrence relation} of the form:
 \begin{equation} \label{recurrence}
  x\begin{pmatrix}
 P^N_0(x) \\  P^N_1(x) \\  P^N_2(x) \\ \vdots
 \end{pmatrix}
 = \mathbf{J} \begin{pmatrix}
 P^N_0(x) \\  P^N_1(x) \\  P^N_2(x) \\ \vdots
 \end{pmatrix} 
 \hspace{1cm} \forall x\in\R .
 \end{equation}
If  there exists a sequence of (infinite) matrices $\J^{(N)}$ such that (\ref{recurrence}) holds for all $N\in\N$ and the matrices $\J^{(N)}$ have a fixed number of non-zero diagonals (i.e.~$\J^{(N)}$ are band matrices). Then, we prove that, if $\J^{(N)}$ has a right-limit  which is a Laurent matrix as $N\to\infty$,  polynomial linear statistics of the biorthogonal ensemble $\Xi_N$ satisfy a CLT. 
Remarkably, the fact that the cumulants $\Cu^n_N\left[ F\right] $ converge to zero as $N\to\infty$ for all $n > 2$ and for any polynomial $F$ boils down to Soshnikov's MCL.
In particular, this shows that  the combinatorial structure of  unitary invariant ensembles in the one-cut regime and the circular unitary ensemble (CUE) is the same.\\

Going beyond \cite{BD_13}, we also compute the limits of the cumulants of linear statistics when the right-limit $\M$ of the recurrence matrices $\J^{(N)}$ is arbitrary; see theorem~\ref{thm:weak_cvg}.  Generically, in this regime,  we cannot expect the fluctuations to be Gaussian and our formulae are significantly more complicated. In fact, it would be worth to investigate further the law of the limiting random variable $\X(f)$ given by formula (\ref{weak_limit})  below.  In particular, in view of the results \cite{Pastur_06, Shcherbina_13, BG_13b} for OPEs in the multi-cut regime, it would be interesting to understand the case where the right-limit $\M$  is periodic or quasi-periodic along its main diagonals, see section~\ref{sect:OPE}. This will be treated in another article.\\
% This generalization is a significant step since it provides the first example of a general limit theorem for linear statistics of biorthogonal ensembles which are not Hermitian.
 
%
%% Overview
%
The rest of this paper is organized as follows. The main results,  which consist of limit theorems for  the biorthogonal ensembles satisfying a finite-term recurrence relation, are formulated in section~\ref{sect:results} and the proofs are given in section~\ref{sect:proof}. 
In section~\ref{sect:discussion}, we provide a short introduction to the theory of BOEs and we review the connection between the main combinatorial lemma of Soshnikov, \cite{Soshnikov_00a}, the Strong  Szeg\H{o} limit theorem and our CLT.
Several applications are discussed in section~\ref{sect:examples}. First, we recall from~\cite{BD_13} that, as a by-product of theorem~\ref{thm:CLT}, we obtain universality for OPEs when the reference measure $\mu_N$ satisfies the so-called one-cut condition. In this case,  we show that the CLT holds for $C^1$ test function as well, see theorem~\ref{thm:OPE}. In section~\ref{sect:product_Ginibre}, based on the results of~\cite{AIK_13,KZ_14}, we prove a CLT for square singular values of the product of independent rectangular  complex Ginibre matrices as the dimensions of the matrices go to infinity in an arbitrary way. In section~\ref{sect:BME}, we obtain a similar CLT for the Laguerre and Jacobi {\it Muttalib-Borodin ensembles}. 
These BOEs were introduced in~\cite{Muttalib_95} as a model for disordered conductors in the metallic regime and have been studied from a mathematical perspective in~\cite{Borodin_99}.
In section~\ref{sect:LLN}, we show how to recover the equilibrium measure of a BOE from the asymptotics of its recurrence coefficients. Combined with an upper-bound on the variance, this implies a law of large numbers for smooth linear statistics.  For any product of square Ginibre matrices, we also discuss how to compute the moments of the equilibrium measure (called the Fuss-Catalan distribution) and how to estimate the tail of the largest square singular value. 
Finally, in the appendix~\ref{A:MCL}, we give a proof of  Soshnikov's main combinatorial lemma which emphasizes on the connection with the Dyson-Hunt-Kac (DHK) formulae and the Bohnenblust-Spitzer combinatorial lemma.

\section{Main results} \label{sect:results}
%
%% Cumulants formulation
%
Let us consider a BOE denoted $\Xi_N=(\xi_k)_{k=1}^N$  (cf$.$ definition~\ref{ensemble}) and suppose that the functions  $(P^N_k)_{k=0}^\infty$ satisfy a recurrence relation of the form (\ref{recurrence})  for a given sequence of  matrices $\J^{(N)}$. We also assume that there exists $\W>0$ such that for all $N\in\N$,
\begin{equation} \label{band}
\J^{(N)}_{ij} = 0  \hspace{1cm}\text{if}\  |i-j| \ge \W . 
\end{equation}

The main examples are OPEs which are discussed in details in section~\ref{sect:OPE}. 
In the following, the matrix $\J$ is called {\it the recurrence matrix}  and we omit the superscript $N$ when there is no ambiguity about the dimension. 
% discuss uniqueness in $\J$
 In~\cite{BD_13}, $\mathbf{J}$ is interpreted as an operator acting  formally on $l^2(\mathbb{Z})$ and the authors established conditions under which the limit of the Fredholm determinant (\ref{Laplace}) exists and can be evaluated; see section~\ref{sect:Szego} for further details. We have a different approach and interpret  $\mathbf{J}$  as the {\it weighted adjacency matrix} of a directed  graph
$\mathcal{G}(\J)$.
In general, to any (infinite) matrix $\mathbf{M}$  corresponds a weighted graph $\mathcal{G}(\mathbf{M})=(V,\vec{E})$ where $V\subseteq \Z$ and $\vec{E}= \big\{ (i,j)\in V\times V : \mathbf{M}_{ij} \neq 0 \big\}$.
The vertex set $V$ inherits the total ordering of $\Z$.
 Edges are oriented,  i.e.~ $(i,j)\neq (j,i)$, and they are weighted or labelled by the entries of $\mathbf{M}$. % see e.g.~figure~(\ref{G_J}). 
% In fact, if $\mathbf{M}$ is not symmetric, we will forget about the  edges orientation. % 
Using this formalism, we express the cumulants $\Cu^n_N[F]$ of a linear statistic $\Xi_N[F]$ for an arbitrary polynomial $F\in\R\langle x\rangle$  in terms of sums over lattice-paths on the adjacency  graph of the  matrix $F(\J)$.
In the following, we let 
$$\Lambda_n = \bigcup_{\ell=1}^{n-1} \big\{ {\bf n} \in \Z^\ell : 0<  n_1<\cdots<n_\ell < n \big\} .$$
 For any $n\ge 2$, the set $\Lambda_n$ is isomorphic to the set of compositions of the integer $n$ under the change of variables
 $\Psi :  \bigcup_{\ell=1}^{n-1} \big\{ {\k} \in \N^{\ell+1} :   k_1+\cdots+k_{\ell+1} = n \big\} \mapsto  \Lambda_n$ given by
 \begin{equation} \label{CV}
 n_1=k_1, \hspace{.3cm} n_2 = k_1+k_2, \hspace{.3cm} \dots\ , \hspace{.3cm}  n_{\ell} = k_1+\cdots+k_\ell . 
\end{equation}
Then, for any ${\bf n} = (n_1,\dots, n_\ell)\in\Lambda_n$, we define
\begin{equation} \label{mho}
 \mho({\bf n}) =  \frac{(-1)^{\ell+1}}{\ell+1} {n \choose \k}  =  \frac{(-1)^{\ell+1}}{\ell+1}  \frac{n!}{k_1! \cdots k_{\ell+1} !}  .    
 \end{equation}

Given a band matrix $\mathbf{M}$, let $\o \in \Z$ be the root of $\mathcal{G}(\M) = (V,\vec{E})$.
We define for any composition ${\bf n}\in \Lambda_n$, 
$ \Gamma^{\bf n}_{\o}(\mathbf{M})$ 
to be the set of all paths $\pi= \big( \pi(0) , \pi(1),\dots, \pi(n)\big)$ on the graph $\mathcal{G}(\M)$ 
such that \begin{align*}
&(i)\ \pi(0)=\pi(n) < \o \\
&(ii)\ \max\big\{\pi(n_1),\cdots,\pi(n_\ell)\big\} \ge \o \ .
\end{align*}

\begin{lemma}\label{thm:cumulants}
For any polynomial $F\in\R\langle x \rangle$. If $\M= F(\J)$, we have  for  $n\ge 2$,
\begin{equation}\label{cumulants}
\Cu^n_N[F] =  \sum_{{\bf n} \in \Lambda_n} \mho({\bf n}) 
\sum_{\pi \in \Gamma^{\bf n}_{N}(\mathbf{M})}\
 \prod_{i=1}^{n} \M_{\pi(i-1)\pi(i)} .
\end{equation}
\end{lemma}

\proof Section~\ref{sect:path}. \qed\\

%
%% Unbounded operators
%

A band matrix $\M$ is bounded (as an operator) if and only if all its entries are uniformly bounded. 
In this case, we can immediately deduce from lemma~\ref{thm:cumulants} some estimates for the cumulants $\Cu^n_N[F]$, as well as for the moment generating function of the linear statistic $\Xi_N[F]$.
 It is also important to consider the case where the recurrence matrix $\J^{(N)}$ is unbounded, for instance to treat OPEs with respect to a measure $\mu_N$ supported on the whole real line.
 Note that formula (\ref{cumulants}) depends only on finitely many entries of $\J$ and it is easy to verify that, if there exits constants $C>0$ such that when $N$ is sufficiently large,\begin{equation*}
 \big| \J^{(N)}_{N+i, N+j} \big| \le C  \hspace{1cm}\forall  |i|, |j| \le N   ,
 \end{equation*}
 then the sequence of random variables $\Xi_N(F)$ is tight. So,  given a subsequence $(N_k)_{k\in\N}$, we provide a sufficient condition for the convergence in distribution of $\Xi_{N_k}(F)$ as $k\to\infty$.

%In section~\ref{sect:path}, lemma~\ref{thm:bound}, using this condition,  we derive a uniform estimate for the cumulants of the random variable $\Xi_N(F)$. Namely, for any polynomial $F$, there exists a constant $\mathrm{C}_F >0$  such  that for any $n\ge 2$, 
%\begin{equation} \label{cumulant_bound}
 %\big| \Cu^n_N[ F] \big| \le n! \exp(n \mathrm{C}_F)   \ ,
 %\end{equation}

%Hence, if (\ref{H1}) holds, the power series (\ref{Laplace_0}) converges uniformly in the disk  $|\lambda| < e^{-\mathrm{C}_F}$ and

%\begin{remark} It is important to consider unbounded matrix $\J$  to treat OPEs with respect to a measure $\mu_N$ supported on the whole real line. The hypothesis (\ref{H1})
 %is not restrictive since, provided that the measures $\mu_N$ are properly rescaled, in general the recurrence coefficients grow at most polynomially.
   %\end{remark}

%
%% Right-limit + local convergence
%

\begin{definition}\label{R_limit}
Let $(\mathbf{M}^{(N)})_{N\in\N}$ be a sequence of (infinite) matrices. We say that $\mathbf{M}$  is a  right-limit of $\mathbf{M}^{(N)}$ along a subsequence $N_k$ if for all $i,j \in \Z $, 
$$\mathbf{M}_{ij}=\lim_{k\to\infty}  \mathbf{M}^{(N_k)}_{N_k+i , N_k+j}  . $$  
Then, we will denote 
$ \displaystyle \mathbf{M}^{(N_k)} \overset{\ell}{\rightarrow} \mathbf{M}$ as $k\to\infty$. Moreover, we say that $\mathbf{M}$  is the  right-limit of ~$\mathbf{M}^{(N)}$ if  
$ \displaystyle \mathbf{M}^{(N)} \overset{\ell}{\rightarrow} \mathbf{M}$ as $N\to\infty$.
\end{definition}

This concept comes from the spectral theory of Jacobi operators and its relevance to fluctuations  of unitary invariant random matrices has been pointed out in~\cite{BD_13}.  
In the context of this paper, it is interesting to reformulate definition~\ref{R_limit} in terms of convergence of graphs. A sequence of weighted graphs $\mathcal{G}_k$ rooted at  $\o_k$ is said to {\it converge locally} to a rooted graph $\mathcal{G}_\infty$ if, for every $n\in\N$, the $n$-neighborhood of the root $\o_k$ stabilizes and the weights in this neighborhood converge as $k\to\infty$; see e.g$.$ \cite{Virag_14} section~3 for another application of graph convergence in the context of random matrix theory.
Hence, a sequence  $(\mathbf{M}^{(N)})_{N\in\N}$ of band matrices has a  right-limit $\mathbf{M}$ along a subsequence $(N_k)_{k\in\N}$   if and only if the sequence of graphs $\mathcal{G}(\mathbf{M}^{(N_k)})$ rooted at $N_k$ converges locally to $\mathcal{G}(\mathbf{M})$ rooted at 0. \\

Motivated by lemma~\ref{thm:cumulants}, let us define for all $n\in\N$, 
\begin{equation}\label{path}
\varpi^{n}_{\o}(\mathbf{M}) =   \sum_{{\bf n} \in \Lambda_n} \mho({\bf n}) 
\sum_{\pi \in \Gamma^{\bf n}_{\o}(\mathbf{M})}\
 \prod_{i=1}^{n} \M_{\pi(i)\pi(i-1)} .
\end{equation}
The definitions of right-limit and local convergence implies that, 
if $\displaystyle \mathbf{M}^{(N_k)} \overset{\ell}{\rightarrow} \mathbf{M}$ as $k\to\infty$,   
then for any $n\in\N$,
 \begin{equation}\label{R_limit_1}
 \lim_{k\to\infty} \varpi^{n}_{N_k}(\mathbf{M}^{(N_k)}) = \varpi^{n}_{0}(\mathbf{M})  ,
\end{equation}
and this yields a {\it limit theorem}.

\vspace{.5cm}

\begin{theorem} \label{thm:Laplace}
Let $\Xi_N$ be a BOE with a recurrence matrix $\J^{(N)}$ whose entries are uniformly bounded in $N$. Suppose that, for a given polynomial $F\in\R\langle x \rangle$, the matrix $F\big(\J^{(N)}\big)$  has a right-limit $\mathbf{M}$ along a subsequence $(N_k)_{k\in\N}$. Then there exists a constant $\mathrm{C}_F>0$ such that
\begin{equation} \label{Laplace_1}
 \lim_{k\to\infty} \left( \log \E{e^{ \lambda \Xi_{N_k}(F)} } - \lambda \E{\Xi_{N_k}(F)}\right)  = \sum_{n=2}^\infty  \varpi^{n}_{0}(\mathbf{M})\frac{\lambda^n}{n!} ,
 \end{equation}
uniformly in compact sets of  the disk  $|\lambda| < e^{-\mathrm{C}_F}$.
\end{theorem}

\proof Let $\M^{(N_k)} = F\big(\J^{(N_k)}\big)$. According to lemma~\ref{thm:cumulants}, for any $N \in\N$, the cumulants of the linear statistic $\Xi_N(F)$ are given by $\Cu^n_{N}[F] = \varpi^{n}_{N} (\mathbf{M}^{(N)})$.
Moreover, by (\ref{R_limit_1}), we have  for any $n\ge 2$, 
\begin{equation} \label{cumulant_lim}
\lim_{k\to\infty}\Cu^n_{N_k}[F] =\varpi^{n}_{0}(\mathbf{M}) \ .
\end{equation}
If we subtract the first term $ \Cu^1_{N_k}[F]= \E{\Xi_{N_k}(F)}$ from the series (\ref{Laplace_0}), the estimate (\ref{cumulant_bound}) proved in section~\ref{sect:path}  guarantees that we can pass to the limits term by term. This yields formula  (\ref{Laplace_1}). \qed\\

If we are only interested in the convergence in distribution of the linear statistic  $\Xi_N[F]$, we do not have to require that the recurrence matrices $\J^{(N)}$ are bounded. Moreover, if  the recurrence matrix has a right-limit,  we obtain the following result.

%%% Limit theorems

\begin{theorem}\label{thm:weak_cvg}  Let $\Xi_N$ be a BOE with a recurrence matrix $\J^{(N)}$. If $\J^{(N)}$ has a right-limit $\mathbf{L}$ along a subsequence $(N_k)_{k\in\N}$ which is a bounded matrix, then
for any polynomial $F\in\R\langle x \rangle$, we have
\begin{equation*}
 \Xi_{N_k}(F)-\E{\Xi_{N_k}(F)} \Rightarrow \mathrm{X}(F)  \hspace{1cm}\text{as }k\to\infty \ .
\end{equation*}
In addition, if we let  $\M=F(\L)$, the Laplace transform of the random variable $\mathrm{X}(F)$ is given by
\begin{equation} \label{weak_limit}
 \E{e^{ \lambda \mathrm{X}(F)} } = \exp\bigg(\sum_{n=2}^\infty  \varpi^{n}_{0}(\mathbf{M})\frac{\lambda^n}{n!}\bigg) \ .
 \end{equation}
\end{theorem}

\proof It is straightforward to verify that if the matrix $\J^{(N)}$ has a right-limit $\L$ along a subsequence $(N_k)_{k\in\N}$, then for any polynomial $F$,
\begin{equation} \label{local_cvg_1}
\lim_{k\to\infty}  F\left(\J^{(N_k)}\right)_{N_k+i , N_k+j} = F(\L)_{ij}  \hspace{1cm}  \forall i,j \in \Z \ . 
\end{equation}
Moreover, since we suppose that the matrix $\L$ is bounded, we have $ \big| \varpi^{n}_{0}\big(F(\L)\big) \big| \le  n!  \exp( n \mathrm{C}_F)$ according to lemma~\ref{thm:bound}.  By a standard argument, e.g.~\cite[Lemma~4.8]{JL_15}, this estimate implies that there exists a unique centered random variable $\mathrm{X}(F)$ whose cumulants are given by $\varpi^{n}_{0}\big(F(\L)\big)$ for all $n\ge 2$.
Finally,  (\ref{cumulant_lim}) implies that, when centered, the linear statistic $\Xi_{N_k}(F)$ converges in distribution (and in the sense of moments) to  $\mathrm{X}(F)$ as $k\to\infty$. \qed

\clearpage

The assumption that the right-limit $\L$ is bounded is not necessary to define the cumulants $\varpi^{n}_{0}(\mathbf{M})$ of the random variable $\X(F)$, see formula (\ref{path}).  However, it guarantees that the Laplace transform (\ref{weak_limit}) is convergent in a small disk $|\lambda| < e^{-\mathrm{C}_F}$. Moreover, in general, it holds as long as the point process $\Xi_N$ is properly scaled so that its empirical measure converges to a compactly supported equilibrium measure. 
In particular, this is well-known to be the case for the ensemble (\ref{HM}) provided that the condition (\ref{confinement}) holds, see section~\ref{sect:OPE} for additional details. \\

%
%% Regularized determinant + quasi-periodic operators
%

%An analogous {\it limit theorem} was established in~\cite{BD_13} using asymptotics of regularized Fredholm determinants. However, for technical reasons, their proof works only for a non-varying measure $\mu_N= \mu$; see theorem~2.4 and remark~2.6 therein. In particular, using the conventions of~\cite[chap.~9]{Simon_05}, we have
%\begin{equation} \label{regularized_det}
%\sum_{n=2}^\infty  \varpi^{n}_{0}(\mathbf{M})\frac{\lambda^n}{n!} = \lim_{m\to\infty} \log\left(  \sideset{}{_2}\det\left[\I+  \mathbf{Q}_m\big(e^{\lambda \mathbf{M}}-1\big) \mathbf{Q}_m  \right] \right) \ , 
%\end{equation}
%where $\mathbf{Q}_m$  is the projection on $l^2(\Z)$ onto $\text{span}(e_{-1},\dots,e_{-m})$. For any $m \ge 1$, the  determinant on the RHS of (\ref{regularized_det}) is analytic for sufficiently small $|\lambda|$ and the coefficients of its Taylor series can be computed using formula (\ref{Laplace_5}), replacing $\J$ by $\M$ and the projection $\p_N$ by  $\mathbf{Q}_m$,    or using the Plemelj-Smithies formulae. Then, adapting the proof of lemma~\ref{thm:cumulants}, it is not difficult to check that  these coefficients  converge to $\varpi^{n}_{0}(\mathbf{M})$ for any $n\ge 2$ as $m\to\infty$.  Moreover, it follows from lemma~\ref{thm:bound} that the series 
%(\ref{regularized_det}) converges uniformly in a small disk around 0.\\
 
Formula (\ref{Laplace_1}) shows that, for the class of BOEs we consider, the fluctuations are described by the right-limits of  the recurrence matrix $\J$. For the matrix model (\ref{HM}), it is known that for a multi-cut potential $V$,  the entries of the matrix $\J$  oscillates and a generic  eigenvalues statistic does not converge in distribution  as $N\to\infty$. However, $\J$  has right-limits along certain subsequences  and theorem~\ref{thm:Laplace} implies that, along such a subsequence, for any polynomial $F$, $\Xi_{N_k}(F)$ converges weakly as $k\to\infty$. 
This mechanism was investigated heuristically by Pastur in~\cite{Pastur_06} using semiclassical formulae for orthonormal polynomials and he showed that, when the asymptotic eigenvalues distribution fills up several intervals, the fluctuations are generically not Gaussian.
In fact, such general limit theorems have recently been proved for $\beta$-ensembles in~\cite{BG_13a,BG_13b}, see section~\ref{sect:OPE} for more details.
According to formula (\ref{Laplace_1}), the limit law of the random variable $\Xi_{N_k}(F)$ is Gaussian only if   $\varpi^{n}_{0}(\mathbf{M})$ vanishes for all $n>2$.
Hence, the right-limit $\mathbf{M}$ must have a very special structure to obtain such subtle cancellations.  In fact, the natural condition is that  $\mathbf{M}$ is constant along its main diagonals.

%% Main combinatorial lemma + Gaussian fluctuations

\begin{definition} \label{Laurent}
 A Laurent polynomial is a function $ s(z)=\sum_{k\in\Z} \widehat{s}_k z^k $ such that only finitely many coefficients  $\widehat{s}_k \in \R$ are non-zero.  We let $\mathbf{L}(s)=\big(\widehat{s}_{j-i}\big)_{i,j \in\Z}$    
be the Laurent matrix with symbol $s(z)$.
\end{definition}

Given a Laurent matrix   $\mathbf{L}$, the graph $\mathcal{G}(\L)$ is translation-invariant and it allows us to simplify formula (\ref{path}). Namely, it is proven in section~\ref{sect:fluctuation} that
\begin{equation} \label{L_path}
\varpi^{n}_{0}\big(\mathbf{L}(s)\big) =    
\sum_{\begin{subarray}{c}\omega_1+\cdots+\omega_n = 0\\ \omega_i \in\Z \end{subarray}} 
\widehat{s}_{\omega_1}\cdots \widehat{s}_{\omega_n}\    \G_n(\omega_1,\dots,\omega_n)  \ ,
\end{equation}
where for any $x\in\R^n$,   
\begin{equation} \label{G}
 \G_n(x_1,\dots, x_n) =  \sum_{{\bf n} \in \Lambda_n} \mho({\bf n}) 
  \max\left\{0,\sum_{i=1}^{n_1} x_{i},\sum_{i=1}^{n_2} x_{i},\cdots, \sum_{i=1}^{n_{\ell}} x_{i}  \right\} \ ,
 \end{equation}
 and $\mho({\bf n})$ is given by (\ref{mho}).
 Formula (\ref{L_path}) is very similar to the limit that Soshnikov obtained for cumulants of linear eigenvalues statistics of Haar distributed random matrices; see formula (\ref{cumulant_CUE}) for comparison. For OPEs, this analogy is motivated by the fact that these models fall in the sine process universality class at local and mesoscopic scales, see~\cite{L_15a} and reference therein. However, it was rather unexpected that, at the global scale, we obtain the same combinatorial structure as for the CUE.
As a consequence of the Main Combinatorial Lemma of ~\cite{Soshnikov_00a}, lemma~\ref{thm:G} below, for any Laurent polynomial $s(z)$, one has
\begin{equation} \label{MCL}
\varpi^{n}_{0}\big(\mathbf{L}(s)\big) =\delta_{n,2} \frac{1}{2}    \sum_{\omega \in \Z}
\widehat{s}_{\omega} \widehat{s}_{-\omega} |\omega| \ .
\end{equation}

 This means that, if the matrix $F\big(\J^{(N)}\big)$ has a right-limit along a subsequence $(N_k)_{k\in\N}$ which is a Laurent matrix, then the 
linear statistic $\Xi_{N_k}(F)$ converges in distribution to a Gaussian random variable as $k\to\infty$.  
 In random matrix theory, the MCL and its variations have become an important devices to prove that fluctuations of eigenvalues statistics are asymptotically Gaussian. Besides   the classical compact groups, the techniques developed in~\cite{Soshnikov_00a} also applies to the sine process and other generalizations,~\cite{Soshnikov_01}.  The MCL was also used in \cite{RV_07a} to keep track of eigenvalues fluctuations for the Ginibre ensemble. In fact, the approach of Rider and  Vir\`{a}g is rather similar to the one used in this article and it can be generalized to any complex OPE with respect to a rotationally invariant measure.   
 A similar lemma appeared first in the work of Spohn on linear statistics of the sine process seen as the invariant measure for Dyson's Brownian motion; see Lemma~2 in \cite{Spohn_87}.  
However, the first correct proof of Spohn's lemma was given in~\cite{Soshnikov_00a} as well. 
 Another related lemma (lemma~\ref{thm:RS}) appeared in \cite{RS_96} to show that, conditionally on the Riemann Hypothesis, the correlation functions for the spacings between $N$ consecutive appropriately rescaled  zeroes of the Riemann Zeta function converge to that of the GUE as $N\to\infty$. In the appendix~\ref{A:MCL}, we show that all these lemmas are consequences of the DHK formulae (\ref{DHK}) which were used  in~\cite{Kac_54} to provide one of the first proofs of the Strong Szeg\H{o} theorem; cf.~formula~(\ref{Szego}) below. \\
 
 Even though it appears in several different contexts, it was pointed out in~\cite{JL_15} that the combinatorial structure behind the MCL is  very sensitive. Namely, if we modify the correlation kernel (\ref{kernel_0}), e.g. by removing the mode $N-m$ for a given $m\ge 2$, we generically obtain a non-Gaussian process in the limit. In this case, the limits (\ref{cumulant_lim}) still hold but we need to `remove' all the paths  such that $\pi(0)=m$ or $\pi(n_j)=-m$ for some $j\in\{1,\dots,\ell\}$ from the definition of 
$\varpi^{n}_{0}(\mathbf{M})$ and, even if $\M$ is a Laurent matrix,  there is no cancellation like (\ref{MCL}). \\

%
%% Central limit theorem
%
As a consequence of formula  (\ref{MCL}), we obtain a CLT for biorthogonal ensembles.
 
 \begin{theorem}\label{thm:CLT}
 Let $\Xi_N$ be a BOE with recurrence matrix $\J^{(N)}$. Suppose that there exists a subsequence $(N_k)_{k\in\N}$ and a Laurent polynomial $s(z)$  such that $\J^{(N_k)}  \overset{\ell}{\rightarrow}  \L(s)$. Then, for any polynomial $F\in \R\langle x\rangle$, 
$$\Xi_{N_k}(F) - \E{\Xi_{N_k}(F)} \Rightarrow \No\big(0, \|F\|^2_s \big) \hspace{1cm}\text{as }k\to\infty \ .$$
The variance is given by
\begin{equation}\label{variance_s}
 \|F\|^2_{s} = \sum_{k=1}^\infty k \widehat{F(s)}_k \widehat{F(s)}_{-k}
 \end{equation}
where
\begin{equation} \label{Fourier}
 \widehat{F(s)}_k = \frac{1}{2\pi i} \oint F\big(s(z)\big) z^{-k} \frac{dz}{z} \ . 
 \end{equation}
 \end{theorem}
\noindent In formula $(\ref{Fourier})$ and in the rest of this article, $\displaystyle \oint$ denotes an integral over the closed contour $\big\{|z|=1\big\}$. 
 
 \proof Section~\ref{sect:fluctuation}. \qed\\

Theorem~\ref{thm:CLT}  first appeared in~\cite{BD_13} as Corollary~2.2 and the authors  provided applications to OPEs, random lozenge tilings,  and to specific instances of two-matrix models. Their results for OPEs are reviewed in  section~\ref{sect:OPE}. We also provide an extension of theorem~\ref{thm:CLT} to rather general test functions for the varying weights (\ref{HM})  when the support of the equilibrium measure is connected; see theorem~\ref{thm:OPE}.
 Actually, this extension follows rather directly from the variance estimates derived in~\cite{L_15a}.   
%For ,  it was proved in~\cite{BD_13} that the eigenvalues of $M$ form a BOE whose recurrence matrix is not symmetric and has $\text{deg}(V)+1$ non-zero main diagonals. Under suitable assumptions on the potential $V$, this matrix has a right-limit which is Laurent matrix and there is a CLT. 
In section~\ref{sect:product_Ginibre}, we give a new application of theorem~\ref{thm:CLT}  to  a BOE which consists of  square singular values of a product of rectangular Ginibre matrices. We obtain a family of CLTs which depend on the ratios between the dimensions of the matrices, see theorem~\ref{thm:MOP}. In section~\ref{sect:BME}, 
we prove that a special cases of this theorem also describes the asymptotic fluctuations of linear statistics of the Laguerre and Jacobi Muttalib-Borodin ensembles, see theorem~\ref{thm:Laguerre} and~\ref{thm:Jacobi} respectively.

%%% Work

\begin{remark} \label{rk:symbol}
In the formulation of theorem~\ref{thm:CLT}, different Laurent polynomials can  give raise to the same variance. In particular, it is easily verified that  $\| F\|_s =\| F\|_{\bar s}$ where $\bar{s}(z)= s(1/z)$ and that, for any $r \neq 0$,  $\| F\|_s =\| F\|_{s_r}$ where $s_r(z)= s(r z)$. These facts can be deduced from the definitions by observing that we have $\L( \bar s) = \L(s)^{\mathrm T}$ and that  $ \L(s_r)= \mathbf{R} \L(s) \mathbf{R}^{-1}$ where the matrix $\mathbf{R}_{kj }= r^{k} \delta_{k j}$. These two transformation arise when considering the biorthogonal families $\{Q_k, P_k\}$ and $\{ r^k P_k, r^{-k} Q_k\}$  instead of $\{P_k, Q_k\}$ to define the process $\Xi_N$. Moreover, according to formula (\ref{kernel_0}), these transformations do not change the correlation kernel $K_N$.  Finally, it is also worth observing that up to scaling and translation of the process $\Xi_N$, the right-limit of the recurrence matrix $\J$ is given by 
$\big( \L(s) - b \mathbf{I}\big)/a = \L(\tilde s)$ where $\tilde{s}(z)= \frac{s(z)-b}{a}$ and $\mathbf{I}$ denotes the identity matrix. 
\end{remark}

\begin{remark} \label{rk:non-Hermitian}
 If the Laurent matrix $\L(s)$ is self-adjoint, then $\widehat{ F(s)}_{-k} =\overline{\widehat{F(s)}_k}$ and by Devinatz's formula, Proposition~6.1.10 in \cite{Simon_05},  
\begin{equation} \label{Devinatz}
 \| F\|_{s}^2 = \frac{1}{8\pi^2} \iint_{[-\pi,\pi]^2} \left| \frac{F\big(s(e^{i\theta})\big) - F\big(s(e^{i\phi})\big)}{e^{i\theta} -e^{i\phi}} \right|^2 d\theta d\phi  \ .
 \end{equation}
This provides an expression for the variance $ \| F\|_{s}^2$ directly in terms of the test function $F$. Moreover, since 
$s(e^{i\theta})$ is real-valued when $\theta\in [-\pi,\pi]$, this formula makes sense for any function $F\in C(\R)$.
For BOEs, the right-limit $\L(s)$ need not be self-adjoint, so that in general $\widehat{ F(s)}_{-k} \neq \overline{\widehat{F(s)}_k}$ and it is not even evident from formula $(\ref{variance_s})$ that $\| F\|_s \ge 0$.  Moreover, since $s(z)$ is not real-valued on the contour $\big\{|z|=1\big\}$,  it is not clear either how to extend formula (\ref{Fourier}) to test  functions which are not analytic.
\end{remark}

\section{Background and discussion} \label{sect:discussion}

In this section, we provide some background on biorthogonal ensembles. In particular, we give a quick proof of formula (\ref{Laplace}) which is the starting point of our cumulants analysis. Then, we review Soshnikov's CLT for linear statistics of Haar-distributed unitary random matrices and we discuss the connection between Soshnikov's results, our results and the Strong Szeg\H{o} theorem.

\subsection{The cumulants method}  \label{sect:cumulants}

Let  $\Xi_N$ be a BOE in the sense of definition~\ref{ensemble}.  By formula (\ref{correlation}), the joint density of the process is given by
\begin{align} \notag
\varrho_N(x_1,\dots,x_N) &=  \frac{1}{N!}  \det_{N\times N}\big[  K(x_i,x_j) \big] \\
& \label{biorthogonal_1}
= \frac{1}{N!}  \det_{N\times N}\big[  P^N_{j}(x_i) \big] \det_{N\times N}\big[  Q^N_{j}(x_i) \big]   ,
\end{align}
where the index $i=1,\dots, N$ and  $j=0,\dots, N-1$. For any $g\in L^\infty(\mu_N)$, we have
\begin{align*} \E{\prod_{k=1}^N g(\xi_k)}& = \int \varrho_N(x_1,\dots,x_N)   \prod_{k=1}^N  g(x_k) d\mu_N(x_k) 
%& =  \frac{1}{N!}  \int \det_{N\times N}\big[ g(x_i)  P^N_{j}(x_i) \big] \det_{N\times N}\big[  Q^N_{j}(x_i) \big]     \prod_{k=1}^N   d\mu_N(x_k) \ ,
\end{align*}
 and by (\ref{biorthogonal_1}) and Andr\'{e}ief's formula:
 \begin{equation*}
 \E{\prod_{k=1}^N g(\xi_k)} =  \det_{N\times N}\left[ \int g(x)  P^N_{i}(x) Q^N_{j}(x) d\mu_N(x) \right]  .
 \end{equation*}
Formally, if we take $g(x)= e^{ f(x)}$ for some test function $f:\R\to\C$ we obtain
 \begin{equation} \label{Andreief}
 \E{e^{\Xi_N(f)}} =  \det_{N\times N}\left[ \int e^{f(x)}  P^N_{i}(x) Q^N_{j}(x) d\mu_N(x) \right] .
 \end{equation}
On the other hand, the recurrence relation  (\ref{recurrence}) implies that for any $n\in \N$ and $j\ge 0$,
$$  x^n P^N_{i}(x) = \sum_{j \ge 0}  \J^n_{ij} P^N_j(x)  .  $$
Note that, since we assume that the matrix $\J$ is a band matrix, the previous sum has only finitely many terms. Besides, when  $\J$ is bounded,  for any real-analytic function $g$, we obtain
$$  g(x) P^N_{i}(x) = \sum_{i \ge 0}  g(\J)_{ij} P^N_j(x)  .  $$ 
In particular, because of the biorthogonal structure (\ref{biorthogonal}), 
 \begin{equation} \label{functional_calculus}
  g(\J)_{ij} = \int  g(x) P^N_{i}(x)Q^N_{j}(x)  d\mu_N(x)   .
  \end{equation}
  Combining formulae (\ref{Andreief}) and (\ref{functional_calculus}),   we obtain
\begin{equation} \label{Breuer_Duits}
 \E{e^{ \Xi_N(f)}} =  \det_{N\times N}\big[   \exp  f(\J) \big]  ,
 \end{equation}
 where we take the determinant of the principal $N\times N$ minor of the infinite matrix $ \exp  f(\J)$.
  In order to investigate the asymptotics of formula (\ref{Breuer_Duits}) as $N\to\infty$, it is convenient rewrite this determinant over $l^2(\N_0)$.  If  $\p_N$ denotes the projection onto $\text{span}(e_1,\dots, e_N)$, we have
\begin{equation} \label{Laplace_2}
 \E{e^{ \Xi_N(f)}} =  \det\big[ \I-\p_N +  \p_N \exp\big(f(\J) \big) \p_N \big] .
 \end{equation}  
This yields formula (\ref{Laplace}).
In fact,  if the operator  $\J$ is bounded and  the test function $f$ is an entire function, then the RHS of (\ref{Laplace_2}) is the Fredholm determinant of a  finite rank operator acting on  $l^2(\N)$. This determinant is unitary equivalent to
\begin{equation}  \label{Laplace_3}
 \E{e^{ \Xi_N(f)}} =  \det\big[1 + \mathcal{K}_N \big( e^f -1 \big)\big]_{L^2(\mu_N)} ,
 \end{equation}    
where $\mathcal{K}_N$ is the integral operator acting on $L^2(\mu_N)$ with kernel $K_N$, (\ref{kernel_0}), and $e^f-1$ is interpreted as a multiplier. Formula (\ref{Laplace_3}) is well-known, see e.g.~\cite{Johansson_05} formula (2.33) for another derivation. 
  Using~formula (5.12) in \cite{Simon_05}, we can express the cumulants of the random variable $\Xi_N(f)$ in terms of the Jacobi matrix $\J$. Namely, if $\lambda \in \C $,
 \begin{equation*} 
\log  \det\big[ 1 +  \p_N \big( \exp\lambda f(\J) -1 \big) \p_N \big]  = \sum_{\ell=1}^\infty \frac{(-1)^{\ell+1}}{\ell}  \tr \left[ \big(\p_N (\exp\lambda f(\J) -1)\p_n \big)^\ell \right]  ,
 \end{equation*}   
and, if we expand the exponential, we obtain
 \begin{equation} \label{Laplace_4}
\log \E{e^{ \lambda \Xi_N(f)}} =    \sum_{n=1}^\infty  \Cu^n_N\left[ f\right] \frac{\lambda^n}{n!}  ,
\end{equation}  
where
\begin{equation} \label{Laplace_5}
\Cu^n_N[f] = \sum_{\ell=1}^{n} \frac{(-1)^{\ell+1}}{\ell} 
\sum_{\begin{subarray}{c}k_1+\cdots+k_\ell = n \\ k_i \ge 1 \end{subarray}} \frac{n!}{k_1!\cdots k_\ell!}
  \tr\big[\p_Nf(\J)^{k_1}\p_N\cdots \p_N f(\J)^{k_\ell} \p_N\big]   .
\end{equation}

If the Jacobi matrix $\J$ is unbounded, when the test function $f$ is a polynomial, it is not a priori plain that the series (\ref{Laplace_4}) is convergent.  However, since $\J$ is a band matrix, the cumulants $\Cu^n_N[f]$ are finite for all $n,N\in\N$. The method of cumulants to investigate fluctuations of linear statistics of large random matrices goes back to Costin and Lebowitz, \cite{CL_95}, and  has  been subsequently developed by many authors.
Much in the spirit of this article is the work of Soshnikov on Haar distributed random matrices from the Unitary group $\mathcal{U}_N$ which is presented in the next section. The combinatorial approach developed in \cite{RV_07a} for polynomial linear statistics of the Ginibre ensemble has a similar flavor too.  Namely, for a rotationally  invariant measure on the plane, it is possible to exploit the obvious recursive structure of the OPs $(\kappa_k^N z^k)_{k=0}^\infty$ to compute the cumulants of any linear statistic of the form $\sum c_k z^{\alpha_k}\bar{z}^{\beta_k}$ where $\alpha_k, \beta_k \in \N_0$. In this case, Soshnikov's main combinatorial lemma also comes into play. 
%n section~\ref{sect:path},  we provide a combinatorial interpretation of the traces on the RHS of  (\ref{Laplace_5}) and, by keeping track of cancellations, we can pass to the limit as $N\to\infty$. 

% based on a  combinatorial devise to compute the limit of cumulants of  polynomial statistics.  Since the Ginibre weight is rotationally invariant, the OPs are given by $(c_k^N z^k)_{k=0}^\infty$ and it is also possible to exploit the obvious recursive structure to compute the cumulants of a linear statistics of the form $\sum c_k z^{\alpha_k}\bar{z}^{\beta_k}$ where $\alpha_k, \beta_k \in \N_0$. 
%Surprisingly, the main combinatorial lemma is also required to show that the limit is Gaussian. However, this approach cannot go beyond the class of rotationally invariant  weights because, in general, the complex OPs do not satisfy a recurrence relation.  

\subsection{The Strong Szeg\H{o} theorem} \label{sect:Szego}  % 1952

The circular unitary ensemble is a complex biorthogonal ensemble where  $P_j(z)=\overline{Q_j(z)}= z^j$ and $\mu_N= \nu$ is the uniform probability measure on the unit circle $\mathbb{T}$. Hence, if $\mathcal{F}$ denotes the Fourier transform on $\mathbb{T}$,  we obtain 
$$  \int e^{\lambda f(z)}  P_{i}(z) Q_{j}(z) d\nu(z) = \mathcal{F}\big(e^{\lambda f}\big)_{j-i}  . $$
By formula (\ref{Andreief}), this implies that the Laplace transform of a linear statistic $\Xi_{\U_N}(f) = \sum_{k=1}^N  f(e^{i\theta_k})$ is given by a Toeplitz determinant with symbol $e^{\lambda f}$:
\begin{equation} \label{Toeplitz_det}
 \E{e^{\lambda\Xi_{\U_N}(f)}} = D_N\big[ e^{\lambda f}\big] .
\end{equation}
If
$\displaystyle \| f \|_{H^{1/2}(\mathbb{T})}^2= \sum_{n\in\Z} |n| \big|\hat{f}_n\big|^2 <\infty$ where $\hat{f}= \mathcal{F}(f)$,  then it follows from the Strong Szeg\H{o} Theorem that
\vspace{-.3cm}
\begin{equation}  \label{Szego}
\lim_{N\to\infty} \log D_N\big[ e^{\lambda f}\big] - \lambda N \int f(z) d\nu(z) = \frac{\lambda^2}{2}  \| f \|_{H^{1/2}(\mathbb{T})}^2  .
\end{equation}
The limits holds for real-valued test function $f\in L^1\cap H^{1/2}(\mathbb{T})$ and for any $\lambda \in \C$. This generalization of the Strong Szeg\H{o} Theorem  to complex-valued symbols is due to Johansson \cite{Johansson88}. 
According to formula (\ref{Toeplitz_det}), this has the following probabilistic interpretation:
\begin{equation} \label{Szego'}
\Xi_{\U_N}(f) - \E{\Xi_{\U_N}(f)} \Rightarrow \No\big( 0 , \| f \|_{H^{1/2}(\mathbb{T})}^2 \big)
\hspace{1cm}\text{as }N\to\infty  . 
\end{equation}

We refer to \cite[Chap.~6]{Simon_04a} and  the review \cite{DIK_13} for some historical background and interesting motivations about the Strong Szeg\H{o} limit theorem. The book \cite{Simon_04a} contains five different proofs and there is yet at least two others. The first by Deift using a Riemann-Hilbert formulation to compute the resolvent of the CUE correlation kernel, \cite{Deift99}. The second by Soshnikov which is based on the cumulant expansion, \cite{Soshnikov_00a}. In particular, this proof does not rely on the fact that the Laplace transform of $\Xi_{\U_N}(f)$ is given explicitly by a Toeplitz determinant but uses instead formula (\ref{Laplace_4}), translation invariance of the CUE kernel and some basic Fourier analysis. In fact, Lemma~1 in~\cite{Soshnikov_00a} implies that  
\begin{equation} \label{cumulant_CUE}
 \lim_{N\to\infty}  \Cu^n_{\U_N}[f] = 2
\sum_{\begin{subarray}{c}\omega_1+\cdots+\omega_n = 0\\ \omega_i \in\Z \end{subarray}}  \hat{f}(\omega_1)\cdots \hat{f}(\omega_n)\    \G_n(\omega_1,\dots, \omega_n)  ,
\end{equation}
where the function $\G_n$ is given by (\ref{G}) and the limit holds for any function $f\in H^{1/2}(\T)$ such that 
\begin{equation} \label{Kac_condition}\ \sum_{\omega\in\Z}   \big|\hat{f}(\omega)\big| <\infty  .
\end{equation}

Moreover,  according to definition (\ref{path}),  if $f(z)$ is a Laurent polynomial,  the RHS of formula (\ref{cumulant_CUE}) turns out to be equal  
$2\varpi^n_0\big(\L(f)\big)$ where $\L(f)$ denotes the Laurent matrix with symbol $f(z)$. 
The difficult part of Soshnikov's proof of the CLT (\ref{Szego'}) is the following Main combinatorial Lemma. 
  \begin{lemma}[MCL, \cite{Soshnikov_00a}] \label{thm:G}
For any $x\in\R^n$ such that $x_1+\cdots +x_n=0$,  we have
 \begin{equation*} \sum_{\sigma\in\Sy_n} \G_n(x_{\sigma(1)},\dots, x_{\sigma(n)}) 
 = \begin{cases} |x_1| &\text{if}\ n=2 \\ 0 &\text{else} \end{cases}  .
 \end{equation*} 
 \end{lemma}
 
 %\proof Appendix~\ref{A:MCL}. \qed\\ 

For completeness, a proof of the MCL inspired from~\cite{Soshnikov_00a} is given in the appendix~\ref{A:MCL}. We also discuss  the connection with the seminal work of  Kac, \cite{Kac_54}, on the Strong Szeg\H{o} theorem. \\

 There is an explicit formula  in terms of a Fredholm determinant for the error term in the Strong Szeg\H{o} theorem,  (\ref{Szego}), which is usually called the Borodin-Okounkov formula. Namely, if $\hat{f}_0=0$, we have
\begin{equation} \label{BO}
 D_N\big[ e^{f}\big]  = \mathrm{Z}(f) \det\big[\mathbf{Q}_N(1-H(f))\mathbf{Q}_N\big]_{l^2(\N_0)}  ,   
\end{equation}
where $\mathrm{Z}(f)= \exp \big( \sum_{k=1}^\infty k \hat{f}_k  \hat{f}_{-k}\big)$, $H(f)$ is an explicit trace-class operator acting $l^2(\N_0)$ and $ \mathbf{Q}_N$ denotes the projection onto $\operatorname{span}\{e_N,e_{N+1},\dots\}$.
 The identity (\ref{BO}) appeared in \cite{BO_00} with a remarkable proof based on the theory of Schur processes developed by the authors. Although it was proved many years earlier in a less explicit form in \cite{GC_79}. We refer to~\cite{Simon_05} section~6.2 for different proofs  and historical references.
  A different approach to the asymptotics of Toeplitz determinant using operator theoretic techniques was introduced by Widom in~\cite{Widom73}.  
Immediately after the work of Borodin and Okounkov, it was realized that Widom's method also allows to prove formula (\ref{BO}) and certain generalizations when $f\notin H^{1/2}(\mathbb{T})$; see  \cite{BW_00} or \cite{Basor_05} for a rather elementary introduction and further references. 
In the context of this article, this method was used in \cite{BD_13} to analyze the  biorthogonal ensembles which satisfy a recurrence relation of the form (\ref{recurrence}) for some band matrix $\J$. The authors were the first to realize that the existence of a right-limit for $\J$ implies a limit theorem for polynomial statistics.  In particular, they already obtained theorem~\ref{thm:CLT} and a weaker version of theorem~{\ref{thm:weak_cvg}}. One of their main tools  is the following asymptotic formula for certain Fredholm determinants,~\cite[Lem.~3.1]{BD_13}. It states that for any trace-class matrix $\mathbf{A}$ acting on $l^2(\N_0)$ and for any Laurent polynomial $s$, if $\T(s)$ is the Toeplitz matrix with entries $\T(s)_{ij}=\widehat{s}_{j-i}$ for all $i,j \in\N_0$, then
\begin{equation} \label{BD}
\lim_{N\to\infty}   \det\big[ \I+  \p_N( e^{\mathbf{T}(s)+\mathbf{A}} -1) \p_N \big] e^{-\tr\big[\p_N (\mathbf{T}(s)+\mathbf{A})\big]} = \mathrm{Z}(s)  .
\end{equation}\

The consequence of these asymptotics is well-illustrated by looking at the Chebyshev process of the second kind. Namely, we consider the OPE,  denoted $\Xi_N$, with respect to the reference measure $d\mu/dx= \frac{2}{\pi} \sqrt{1- x^2} \1_{|x|<1}$. 
The Jacobi matrix $\J$ of the measure $d\mu$ is the Toeplitz matrix with symbol $s(z)=\frac{z+1/z}{2}$;
\begin{equation*}
\mathbf{J} = \begin{pmatrix}
 0 & 1/2 & 0 & 0   \\
 1/2 & 0 & 1/2 & 0   &\mathbf{0} \\
  0 &  1/2&  0 & 1/2   \\ 
 &\mathbf{0}   & \ddots &\ddots&\ddots
\end{pmatrix}  .
\end{equation*}

It is easy to check that, for any polynomial $F\in \R\langle x\rangle$, $F(\J)=\T\big(F(s)\big) + \mathbf{A}$  where $\mathbf{A}$ is a finite rank matrix. Hence, since $\E{ \Xi_N(F)}= \tr[  \p_N F(\J)]$, according to formula (\ref{Laplace_2}),
\begin{equation*}
 \E{e^{ \Xi_N(F) - \E{ \Xi_N(F)}}} =   \det\big[ \I+  \p_N( e^{\mathbf{T}(s)+\mathbf{A}} -1) \p_N \big] e^{-\tr\big[\p_N (\mathbf{T}(s)+\mathbf{A})\big]}  
\end{equation*}
and, by (\ref{BD}), we obtain 
\begin{equation*}
 \lim_{N\to\infty}\E{e^{ \Xi_N(F) - \E{ \Xi_N(F)}}} = \mathrm{Z}\big(F(s)\big) .  \end{equation*}
 Since, by definition, $ \mathrm{Z}\big(\lambda F(s)\big) = \lambda^2 \|F\|^2_{s}$ for any $\lambda\in\R$, we get a special case of theorem~\ref{thm:CLT}.  The general case follows from a universality principle.
Breuer and Duits  showed that, if two recurrence matrices have the same right-limit, then for a given polynomial $f$, the cumulants  (\ref{Laplace_5}) also  have the same limits. Hence, if  $\J^{(N_k)}  \overset{\ell}{\rightarrow}  \L(s)$ where $\L(s)$ is a Laurent matrix,  we may replace  $\J^{N}$  by a Toeplitz matrix $\T(s)$ to compute the limit of the moment generating function of a given polynomial statistic and the central limit  theorem~\ref{thm:CLT} follows directly from formula (\ref{BD}).

\section{Examples}\label{sect:examples}

\subsection{Orthogonal Polynomial Ensembles} \label{sect:OPE}
%
%%  Jacobi operators
%
Let $(\mu_N)_{N\in\N}$ be a sequence of Borel measures on $\R$ satisfying the condition (\ref{reference}) for all $N\in\N$ and consider the determinantal process $\Xi_N$ with correlation kernel of the form (\ref{kernel}) where $P_k^N$ are the OPs with respect to $\mu_N$.
One can either consider discrete measures, e.g.~coming from applications to random tilings of planar domains, or absolutely continuous measures coming from unitary invariant random matrices. We will focus on the second class and, in particular, on the exponential weights  $\mu_N/dx = e^{-NV(x)}$, see (\ref{HM}),  which are  natural generalization of the GUE ($V(x)=2x^2$).  One can also consider non-varying measures with a compact support, e.g.~the Jacobi ensembles, whose behavior is similar and also fits in the context of the discussion below. 
%Alternatively, OPs can be defined by a three term recurrence relation for all $N,n\in\N$,
In any case, the orthogonal polynomials  satisfy a  three term recurrence relation:
\begin{equation}\label{recurrence'}
xP^N_n(x) = a^N_n P^N_{n+1}(x) + b^N_n P^N_n(x) +  a^N_{n-1} P^N_{n-1}(x)   ,
\end{equation}
which is of the type (\ref{recurrence}) where
\begin{equation*}
\mathbf{J} = \begin{bmatrix}
 b_0^N & a_0^N & 0 & 0 & 0  \\
 a_0^N & b_1^N & a_1^N & 0 & 0  &\mathbf{0} \\
  0 &  a_1^N  &  b_2^N  &  a_2^N  & 0 \\ 
 0 & 0 &  a_2^N &  b_3^N &  a_3^N \\
 &\mathbf{0} &  & \ddots &\ddots&\ddots
\end{bmatrix}
\end{equation*}
is called the  Jacobi matrix. Note that  $\mathbf{J}$ may depend on the dimension $N \in \N$, and  we may denote $\mathbf{J}= \mathbf{J}(\mu_N)$ to emphasis the dependency.  
We can interpret $\mathbf{J}$ as a symmetric operator acting  formally on $l^2(\N_0)$, then 
the reference measure $\mu_N$ is the spectral measure of $\mathbf{J}$ at $e_0$. Hence, $\J$ is bounded if and only if the measure $\mu_N$ is compactly supported. In addition, if the support of the measure $\mu_N$ is unbounded, $\J$ is essentially self-adjoint if and only if the moment problem for $\mu_N$ is determinate; see \cite[chapter~2]{Deift_99}.  
In the framework of this paper, the definiteness of $\J$ is not relevant since we view of $\mathbf{J}$ as the {\it adjacency matrix} of a  weighted graph:
%\begin{figure}[h]
%\begin{equation*} \label{G_J}
\begin{center}
\begin{tikzpicture}
 \node at (-1,0) {$\mathcal{G}(\J) =$} ;
 \foreach \x in {0,1,...,5}{\node(\x) at (1.6*\x,0) [draw,circle,inner sep=3pt ,fill]  {} ; }
 \node(6) at (10,0) {} ;
 \foreach \x in {0,1,...,5}{ \node at (1.6*\x,-.5) {\x} ; }
\path
 \foreach \x in {0,1,...,5}{ (\x)  edge [loop above]  node {$b^N_\x$} ()  }
   (0)      edge                 node [above]  {$a^N_0$}     (1)
  (1)      edge                 node [above]  {$a^N_1$}     (2)
  (2)      edge                 node [above]  {$a^N_2$}     (3)
  (3)      edge                 node [above]  {$a^N_3$}     (4)
  (4)      edge                 node [above]  {$a^N_4$}     (5)
  (5)      edge[dashed]                 node   {}     (6) ;   
  \end{tikzpicture}  . 
 % \end{equation} 
%\end{figure}
\end{center}

%
%% Equilibrium point configuration
%
The first fact is that a typical point  configuration of the process $\Xi_N$ reaches an equilibrium as $N\to\infty$. Namely, there exists a probability measure $\mu_*$ on $\R$ such that for any  $f\in C\cap L^\infty(\R)$, 
\begin{equation} \label{equilibrium}
 N^{-1} \Xi_N(f) \underset{N\to\infty}{\longrightarrow}  \int f(x) d\mu_*(x) dx   . 
\end{equation}
The convergence holds  almost surely and $\mu_*$ is called the {\it equilibrium measure}.  It is defined as the minimizer of a weighted logarithmic energy, it is compactly supported and absolutely continuous. In particular, for the matrix ensembles (\ref{HM}), the equilibrium density, denoted  $\rho_V$, exists as long as the potential $V(x)$ is lower-semicontinuous and  satisfies the condition (\ref{confinement}),  see e.g.~\cite{Johansson_98} or \cite[chap~6]{Deift_99}.  The law of large numbers (\ref{equilibrium})  usually follows from concentration estimates, see e.g.~\cite{BD_14,Hardy_15} in the context of this paper or \cite{AGZ} for more general results. In particular, the lattice path method developed in section~\ref{sect:path} is inspired from~\cite{Hardy_15} where it  was used to  estimate the variance of any polynomial linear statistic and show that the zeroes of the polynomial $P_N^N$ are also distributed according to $\mu_*$ as $N\to\infty$. Actually, Hardy's result implies that if there exists two continuous functions $b:\R_+ \to\R $ and $a:\R_+ \to [0,\infty) $  such that
\begin{equation} \label{equilibrium_1}
 a^{N}_{n} \to a(t)    \hspace{1cm}\text{and}\hspace{1cm}    b^{N}_{n} \to b(t) 
\end{equation} 
as $n/N \to t$, then the equilibrium density is given by 
\begin{equation} \label{equilibrium_0}
 \rho_V(x) = \int_0^1 \1_{\beta_-(t)<x<\beta_+(t)} \frac{ dt}{\sqrt{(\beta_+(t) - x)(x-\beta_-(t))}}  \hspace{.8cm}\text{where }\
 \beta_\pm(t)= b(t)\pm 2a(t)  . 
 \end{equation}

The assumption (\ref{equilibrium_1}) is not necessary for the existence of the equilibrium density, however it holds for many of the classical examples from the literature both in the continuous and discrete setting, see \cite{KV_99}.
Moreover, it implies that the recurrence matrix $\J(\mu_N)$ has a right-limit in the sense of definition~\ref{R_limit} which is a Laurent matrix with symbol $s(z)=az + b + az^{-1}$ where $a=a(1)$ and $b=b(1)$. In this case, theorem \ref{thm:CLT} describes the fluctuations around the equilibrium configuration given by $\rho_V$. 
%
%% Multicut regime +  Quasi periodic right-limit
%
Under the assumption \eqref{equilibrium_1}, we give a direct proof of formula  \eqref{equilibrium_0} in section~\ref{sect:LLN} using the lattice-path method. \\

For general OPEs, the fluctuations of linear statistics are expected to remain bounded for large~$N$. However, there is a CLT
 only if the support of the equilibrium measure is  connected. Otherwise, generically, the variance is quasi-periodic and there are limits only along  some particular subsequences. Moreover, the asymptotic distribution need not be Gaussian,  see~\cite{Pastur_06} or \cite[Chapter~14]{Pastur_Shcherbina} for a nice heuristic argument. This oscillatory behavior is explained by the fact that the fluctuations of the number of eigenvalues  in the different components of the support do not stabilize as $N\to\infty$. In \cite{Shcherbina_13, BG_13b}, such results have been obtained using the $1/N$  expansion of the integral
\begin{equation} \label{partition_function} 
Z_{N,V}(\beta) = \int e^{-\beta \mathscr{H}_V(x_1,\dots, x_N)/2} dx_1\cdots dx_N  , \end{equation}
where, if $V$ is a real-analytic and confining potential, the {\it Hamiltonian} is given by 
\begin{equation} \label{Hamiltonian} 
\mathscr{H}_V(x_1,\dots, x_N)= -\sum_{i\neq j} \log|x_i - x_j| + N \sum_{i} V(x_i)  .
\end{equation}

These results  go beyond the context of this paper, since there are valid for arbitrary $\beta>0$ when there is no determinantal structure. Moreover, the work of Borot and Guionnet, \cite{BG_13a,BG_13b}, goes beyond the CLT and established the existence of the all order asymptotic expansion for the partition function (\ref{partition_function}).  Formally, we briefly summarize their results regarding the asymptotic distribution of linear statistics. 
%Observe  that, if $f$ is a  real-analytic test function,
%$$  \E{e^{ \Xi_N(f)}}  = \frac{Z_{N}^{\beta}\big[V+  2f /N\big]}{ Z_{N}^{\beta}[V]} \ . $$
%
Let $n \ge 1$ and  suppose that $ \supp(\rho_V) =  \bigcup_{j=0}^n I_j$ where $I_j$ are closed non-empty intervals. If we let 
$\displaystyle  \vec{\epsilon} =\big(  \rho_V(I_1),\cdots,  \rho_V(I_n) \big) \in (0,1)^n$ be the so-called filling fractions and $\vec{N}= N\vec{\epsilon}$, then
\begin{equation}\label{BG}
  \E{e^{ \Xi_N(f)}} e^{- N  \langle f, \rho_V \rangle} \underset{N\to\infty}{\sim} e^{- (\beta/2-1)  v_V[f]   + \frac{\beta}{4} \mathscr{Q}[f,f]} \frac{\vartheta_{\vec{N}}(c_{\beta,V} +  \frac{\beta}{2} u_V [f]\hspace{.05cm} | \tau ) }{\vartheta_{\vec{N}}(c_{\beta,V} | \tau)}  ,
  \end{equation}
where $c_{\beta,V} \in  \R^n$,  $u_V, v_V$ are linear functionals, $\mathscr{Q}$ is a quadratic functional, $\tau$ is a positive definite $n \times n$  matrix, and 
\begin{equation}\label{theta}
 \vartheta_{\vec{N}}(v | \tau) =  \sum_{\  k \in \vec{N}+ \Z^n} e^{-\frac{1}{2} \langle \tau  k  , k \rangle +  \langle v , k\rangle } \ . 
 \end{equation}
For now on, we suppose that $\beta=2$ since it corresponds to the Hermitian matrix models (\ref{HM}). This values is special, since on the r.h.s.~of (\ref{BG}), the linear term vanishes and $c_{2,V} \equiv 0$.  
Then the asymptotics \eqref{BG} imply that the random variable $ \Xi_N(f)- N  \langle f, \rho_V \rangle$ has a limit and this limit is Gaussian if and only if $u_V [f]=0$, i.e.~the test function $f$ lies in the kernel of an $n$-dimensional linear system. Otherwise, there are limits only along subsequences such that  
\begin{equation}\label{charges}
(\{N_k \epsilon_1\}, \dots, \{N_k \epsilon_n\}) \to  \vec{q} \hspace{.7cm}\text{as}\hspace{.3cm} k \to\infty  ,
\end{equation}
where $\{x\}$ denotes the integer part of $x>0$. Note that, by periodicity, we may replace   
$\vec{N}$ by $(\{N \epsilon_1\}, \dots, \{N \epsilon_n\})$ in the definition of the Theta function (\ref{theta}). Then, the condition (\ref{charges}) implies that $\Xi_N(f) - N  \langle f, \rho_V \rangle$ converges in distribution to the sum of two independent random variables. A real-valued Gaussian with mean-zero and variance  $\mathscr{Q}[f,f]$ and a random variable $\Gamma$ which is  the projection on the vector $u_V[f]$ of a {\it discrete Gaussian random vector} with mean $\vec{q}$  and covariance matrix $\tau$. Namely, if $Y \in \Z^n$ is a random variable with distribution
$$ \P{Y=k} = e^{-\frac{1}{2} \langle k  , \tau k \rangle } , $$ 
then $\Gamma= \langle \vec{q} + Y , u_V[f] \rangle$ has the moment generating function
$$ \E{e^{\lambda \Gamma}} =  \frac{\vartheta_{\vec{q}}( \lambda u_V [f]\hspace{.05cm} | \tau ) }{\vartheta_{\vec{q}}(0 | \tau)} 
=  \frac{\vartheta_{0}(\tau \vec{q}+ \lambda u_V [f]\hspace{.05cm} | \tau ) }{\vartheta_{0}( \tau \vec{q}\hspace{.05cm} | \tau)}    . $$

In~\cite{Pastur_06},   Pastur already noticed that  (\ref{charges}) is the natural condition for the Laplace transform of $\Xi_N(f)$ to converge and that the limiting laws depend on the parameter $\vec{q} \in [0,1]^n$. Like theorem~\ref{thm:weak_cvg}, his results  are based on the hypothesis that the Jacobi matrix $\J$ has right-limits.   
Namely, from the asymptotics of Deift et al. \cite{Deift_al_99_b},  Pastur suppose that the recurrence coefficients satisfy for any $k\in\Z$,
\begin{equation}\label{Pastur}
a_{N+k}^N  \sim \mathcal{R}( N\vec{\epsilon} + k \vec{\alpha}) \ , \hspace{1cm}
b_{N+k}^N \sim \mathcal{S}( N\vec{\epsilon} + k \vec{\alpha})  ,
\end{equation} 
where $ \mathcal{R}: \mathbb{T}^n \to \R$  and $ \mathcal{S}: \mathbb{T}^n \to \R$  are continuous functions, $\vec{\alpha}\in \R^n$, and       
 $\mathbb{T}= \R/\Z$. This hypothesis implies that the Jacobi  matrix $\J(\mu_N)$ has a right-limit $\M_{\vec{q}}$  along the subsequence (\ref{charges}). Hence, by theorem~\ref{thm:weak_cvg}, once centered, the sequence  $\Xi_{N_k}(F)$  converges in distribution as $k\to\infty$.   Generically, the matrix $\M_{\vec{q}}$ is quasi-periodic and there is no reason why the cumulants $ \varpi^{n}_{0}(\M_{\vec{q}})$ would vanish for all $n \ge 3$. In fact, it would be very interesting to investigate the moment generating function (\ref{weak_limit}) in such cases to see if one can recover the description (\ref{BG}).
 In particular, in the case of the quartic potential $V (z) =z^4/4-tz^2/2$ where the asymptotics of the recurrence coefficients are known explicitly and the right-limits of the Jacobi matrix are two-periodic along  its diagonals, see \cite[Theorem~1.1]{BI_99}. \\
  
 %
 %% One-cut regime + Johansson's theorem
 %
 In  the one-cut situation, i.e.~when  $ \supp(\rho_V)= [ b -2a, b+2a]$ with $b\in\R$ and $a>0$, the situation is simpler, the recurrence coefficients have a limit:
 \begin{equation} \label{equilibrium_2}
 a^{N}_{n} \to a    \hspace{.8cm}\text{and}\hspace{.8cm}  b^{N}_{n} \to b \hspace{.8cm} \text{as} \hspace{.5cm}  \frac{n}{N}\to 1  .
\end{equation}
  For the exponential weights $\mu_N/dx = e^{-NV(x)}$, these limits are consequences of the work of Johansson \cite{Johansson_98,KS_10}. For a non-varying measure of the form $d\mu_N/dx = w \1_{I}$ where $I$ is a finite interval and $w>0$ on $I$, they follow from the celebrated  Denisov-Rakhmanov theorem, \cite[Theorem~1.4.2]{Simon11}.
The asymptotics~(\ref{equilibrium_2}) imply that $\J(\mu_N)  \overset{\ell}{\rightarrow} \L(s)$ as $N\to\infty$ where $L(s)$ is a Laurent matrix with symbol $s(z)=az+b +a z^{-1}$. Hence, by theorem~\ref{thm:CLT}, we obtain a CLT as $N\to\infty$, 
\begin{equation} \label{CLT}
\Xi_{N}(F) - \E{\Xi_{N}(F)} \Rightarrow \No\big(0,  \|F\|_s^2 \big)  .
\end{equation}
 %where the variance $\|F\|_s^2$ is given by formulae (\ref{variance_s} - \ref{Fourier}).  
 
 For general $\beta>0$ and  $C^2$ test functions, the counterpart of this CLT also appeared in \cite{Johansson_98}  if  the potential $V(x)$ is an even degree polynomial. Then, it was generalized to  any analytic potential in~\cite{KS_10}. For certain tridiagonal $\beta$-ensembles, similar CLTs have also been derived using the method of moments and lattice-path counting, see \cite{DE_06, DP_12}.
This shows that fluctuations of Hermitian random matrices are universal. By an affine transformation of the potential, we can always suppose that  $a=1/2$ and $b=0$. Then, according to formulae  (\ref{variance_s} - \ref{Fourier}),  the limiting variance of the linear statistic $\Xi_{N}(F)$ is  given by
\begin{equation}\label{variance_T}
 \Sigma(F)^2 = \frac{1}{4}  \sum_{k=1}^\infty k \big|\mathrm{c}_k(F) \big|^2
 \end{equation}
 where, if $T_k$ denote the Chebyshev polynomials of the first kind ($T_k(\cos \theta) = \cos k\theta$), then 
\begin{equation}\label{Fourier_T}
\mathrm{c}_k(F)=\frac{2}{\pi} \int_{-1}^1 F(x) T_k(x) \frac{dx}{\sqrt{1-x^2}} 
\end{equation}
 are the Fourier-Chebyshev coefficients of the function $F$.  Another explicit formula for the variance (\ref{variance_T}) can be deduced from Devinatz's formula (\ref{Devinatz}). One has
\begin{equation}
 \Sigma(F)^2 = \frac{1}{8\pi^2} \iint_{[-\pi,\pi]^2} \left| \frac{F\big(\cos\theta\big) - F\big(\cos\phi\big)}{e^{i\theta} -e^{i\phi}} \right|^2 d\theta d\phi  \ ,
 \end{equation}
 and by symmetry,
$$ \Sigma(F)^2= \frac{1}{4\pi^2} \iint_{[0,\pi]^2} \left| F(\cos\theta) - F(\cos\phi) \right|^2 
\bigg\{ \frac{1}{|e^{i\theta} -e^{i\phi}|^2}+  \frac{1}{|e^{i\theta} -e^{-i\phi}|^2} \bigg\}  d\theta d\phi \ . $$
Making the change of variables $x=\cos\theta$ and $y=\cos\phi$, we obtain 
\begin{align*}  \frac{1}{|e^{i\theta} -e^{i\phi}|^2}+  \frac{1}{|e^{i\theta} -e^{-i\phi}|^2}
&=   \frac{1}{(x-y)^2 + (\sqrt{1-x^2}- \sqrt{1-y^2})^2 }+  \frac{1}{(x-y)^2 + (\sqrt{1-x^2}+ \sqrt{1-y^2})^2 } \\
&=  \frac{4-4xy}{4(x-y)^2}
\end{align*}
after some elementary but non-obvious simplifications. This implies that
\begin{equation}\label{variance} 
 \Sigma(F)^2 = \frac{1}{4\pi^2} \iint_{[-1,1]^2} \left| \frac{ F(x) - F(y)}{x-y} \right|^2 \frac{1-xy}{\sqrt{1-x^2}\sqrt{1-y^2}} dxdy 
\end{equation}
which is a weighted $H^{1/2}$ Sobolev semi-norm. Optimally, the CLT (\ref{CLT}) should hold for any function $f:\R\to\R$ such that  $ \Sigma(f)<\infty$ and has some reasonable growth at $\infty$. However, because one needs to control the fluctuations the edges of the spectrum, it is not known if this holds even for the GUE.
%  In~\cite{BD_13}, when $\mu_N$ have a fixed  compact support, the authors discussed  the extension of (\ref{CLT}) to test function $f \in C^1(\R)$ with polynomial growth. 
In~\cite{L_15a}, some general variance estimates for OPEs have been derived in the one-cut regime. Namely, if the measures $\mu_N$ have densities $w_N$ and the OPs satisfy the  Plancherel-Rotach asymptotics:
\begin{equation} \label{semiclassical}
 P_{N-\delta}^N(x)= \sqrt{\frac{2}{\pi}} \frac{\cos\big[ N\pi \F(x) + \psi_\delta(x)\big]}{\sqrt{w_N(x)}(1-x^2)^{1/4}} + \underset{N\to\infty}{o(1)} 
\end{equation}
uniformly for all $x$ in compact subset of $[-1,1]$ for $\delta=0,1$, the functions $ \psi_0,\psi_1 \in C^1$ and $\F$ denotes the integrated density of states, i.e.~$\F' = \rho_V$ on $[-1,1]$ and $\F(1)=0$.  Then, for any function
$f\in C^1(\R)$ such that there exists  $Q, k>0$ and  $|f'(x)|\le Q |x|^{k}$ for all $|x| \ge 1$,  we obtain 
\begin{equation}\label{variance_estimate}  \limsup{N\to\infty} \Var\big[\Xi_N (f)\big] \le  \frac{16}{\pi^2}  \iint\limits_{[-1,1]^2} \left| \frac{f(x)- f(y)}{x-y} \right|^2 \frac{dxdy}{\sqrt{1-x^2}\sqrt{1-y^2}}  \ . 
\end{equation}

In fact, the estimate (\ref{variance_estimate}) is valid for more general test function $f$ but we prefer to keep a simple formulation.  For the weight $w_N(x)= e^{-NV(x)}$, if $\supp \rho_V  =[-1,1]$, the asymptotics (\ref{semiclassical}) have been derived in \cite{Deift_al_99_b} using the Riemann-Hilbert steepest descent method.
  For non-varying measures,  $d\mu_N/dx = w \1_{[-1,1]}$ where $w>0$ and satisfies mild regularity condition on $[-1,1]$, the Plancherel-Rotach asymptotics  are classical  and $\F(x)=\arccos(x)$, see \cite[Theorem~12.1.4]{Szego}. % In fact, in the later case, a similar estimate can be obtained without using the precise  asymptotics (\ref{semiclassical}), see~\cite{BD_13} formula (5.7). 
 %\begin{equation}  \limsup{N\to\infty} \Var\big[\Xi_N (f)\big] \le C\sup_{x,y \in [-1,1]} \left| \frac{f(x)- f(y)}{x-y} \right|^2  \ , 
%\end{equation}
 %
  As discussed in~\cite[Appendix~A]{L_15a}, the estimate (\ref{variance_estimate}) allows us to show that (\ref{CLT}) is also valid for $C^1$~test~functions.

\begin{theorem}\label{thm:OPE}   Let $V:\R\to\R$ be a real analytic function which satisfies the condition $(\ref{confinement})$ and such that the equilibrium measure satisfies $\supp(\rho_V)=[-1,1]$.  
 If $(\lambda_1, \dots, \lambda_N)$ denote the eigenvalues of a random matrix distributed according to  $\mathbb{P}_{N,V}$, then for any $f\in C^1(\R)$ such that there exists $Q,k>0$  and  $|f'(x)|\le Q |x|^{k}$ for all $|x| \ge 1$, we have
\begin{equation} \label{clt_2'}
\sum_{k=1}^N f(\lambda_k)  - \E{\sum_{k=1}^N f(\lambda_k)}  \underset{N\to\infty}{\Longrightarrow}   \No\big( 0, \Sigma(f)^2 \big)    \ .
\end{equation}
\end{theorem}

%
%% The Chebyshev noise
%
Theorem~\ref{thm:OPE} has the following interpretation: viewing $\Xi_N = \sum_{k=1}^N \delta_{\lambda_k}$ as 
a random measure on $\R$, once centered, it converges in distribution to a Gaussian field $\mathfrak{G}$ 
%(random tempered distribution) 
supported on  the interval $[-1,1]$  with covariance 
\begin{equation}\label{noise}
\E{ \mathfrak{G}(f) \mathfrak{G}(g)} = \frac{1}{4} \sum_{k=1}^\infty k \mathrm{c}_k(f)\mathrm{c}_k(g) \ . 
\end{equation}
The  random process   $\mathfrak{G}$ is realized as the tempered distribution $  \mathfrak{G}(x) = \frac{1}{2} \sum_{k=1}^\infty  \sqrt{k} \xi_k T_k(x)$ where $\xi_1, \xi_2, \dots $ are i.i.d.~$\No(0,1)$. Although the previous sum diverges pointwise, for any function $f \in L^2([-1,1], \nu)$ where $d\nu= \frac{2}{\pi}(1-x^2)^{-1/2} \1_{|x|<1} dx$ such that $\sum_{k=1}^\infty  k \mathrm{c}_k(f)^2 <\infty $, the random variable
$$\mathfrak{G}(f):=  \int_{-1}^1 f(x) \mathfrak{G}(x)  d\nu(x)  = \frac{1}{2}   \sum_{k=1}^\infty  \sqrt{k} \xi_k \mathrm{c}_k(f) $$
 converges almost surely.\\

The proof of theorem~\ref{thm:OPE} and the estimate (\ref{variance_estimate}) are given in the appendix~A of \cite{L_15a}, we just sketch the main arguments.  First, (\ref{variance_estimate})  implies that the sequence of random variables $X_N= \Xi_N(f)- \E{\Xi_N(f)}$ is tight. Then, we  approximate the test function $f$ by a sequence of polynomial $(F_k)_{k\in\N}$ so that  
$$ \sup\big\{ | f'(x)-F'_k(x)| : |x| \le 1 \big\} \le  1/k \ . $$  
By theorem~\ref{thm:CLT}, $\Xi_N(F_k) - \E{\Xi_N(F_k)} \Rightarrow \mathfrak{G}(F_k)$ as $N\to\infty$  for any $k\in\N$ and   we conclude that the limiting  distributions of $X_N$ as $N\to\infty$  and $\mathfrak{G}(F_k)$ as $k\to\infty$ coincide and equal to $\mathfrak{G}(f)$. 

\begin{remark} The Chebyshev polynomials $(T_j)_{j=1}^\infty$ are orthonormal with respect to the measure $d\nu= \frac{2}{\pi}(1-x^2)^{-1/2} \1_{|x|<1} dx$. Thus, according to (\ref{Fourier_T}), $\mathrm{c}_k(T_j)= \delta_{k,j}$  and by formula (\ref{noise}), this implies that the random variables $\big( \frac{2}{\sqrt{j}} \mathfrak{G}(T_j)\big)_{j=1}^\infty$ are independent identically distributed $\No(0,1)$. This theorem first appeared in~\cite{Johansson_98} and we provide an alternative proof in the appendix~\ref{A:Chebyshev}  which is based on  the fact that,  when the condition  (\ref{equilibrium_2}) holds with the normalization $a=1/2$ and $b=0$, $\nu/2$ is the spectral measure of the right limit of the matrix  $\J(\mu_N)$.         
%e fact that the arcine distribution $d\nu= \frac{2}{\pi}(1-x^2)^{-1/2} \1_{|x|<1} dx$  arises in the covariance structure of $\mathfrak{G}$ may be expected since the condition (\ref{equilibrium_2}) with $b=0$ is equivalent to the fact that the graph $\g\big(\J(\mu_N)\big)$ rooted at $N$ converges locally to $\Z$ with weight $a$ on each edges and the spectral measure of $\Z$ is $\nu$.
\end{remark}

\subsection{Product of independent complex Ginibre matrices} \label{sect:product_Ginibre}

Recently, it was established that the square singular values of a product of rectangular complex Ginibre matrices forms a biorthogonal ensemble, \cite{AIK_13}.  In \cite{KZ_14}, it was further proved  that this point process can be interpreted as a multiple orthogonal polynomial ensembles of type II and the authors found explicit formulae for the corresponding polynomials and their recurrence coefficients. Then, they used 
a double contour integral representation for the correlation kernel to obtain its scaling limit at the hard edge. It describes a new universality class which generalizes the classical Bessel kernels.
In this section, we obtain the asymptotics of the recurrence coefficients and deduce from theorem~\ref{thm:CLT} that the fluctuations of the square singular values process are described by a CLT.\\
 
Let $X_j$ be an $N_j \times N_{j-1} $ random matrix whose entries are independent complex Gaussian with mean 0 and variance 1. We consider the square singular values of the product
\begin{equation} \label{W_matrix}
 W_N= X_m\cdots X_1 .
 \end{equation}
 
We denote $N_j = N+ \eta_j$ and we suppose that $\eta_0=0$ and $\eta_1,\dots, \eta_m \ge 0$, so that  $ W_N^{*} W_N$ is an $N\times N$ random matrix and, almost surely, all its eigenvalues are  positive. In~\cite{AIK_13}, it was proved that the square singular values $\mu_1,\dots, \mu_N$ of the matrix  $(\ref{W_matrix})$
 have a j.p.d.f.~on $(\R_+)^N$ which is given by
\begin{equation} \label{AIK}
\varrho_N(x_1,\dots,x_N) = \frac{1}{Z_N}  \det_{N\times N}\big[ x_i^{j-1} \big] \det_{N\times N}\big[ w_{j-1}(x_i) \big] ,
\end{equation} 
where, for any $n\in\N_0$,
$$w_n(x)= \frac{1}{2\pi i} \int_L \Gamma(s+\eta_1+ n) \prod_{l=2}^m \Gamma(s+\eta_l) x^{-s} dx  .  $$
The contour of integration is $L= \{ c + i y, y\in\R \}$ with $c>0$ and $w_0,w_1,\dots$ are Meijer G-functions, although we will not need this fact in the following. Akemann et~al.~used the determinantal structure to investigate the one-point correlation function and its moments.  %They proved that, in a certain scaling limit, it converges as $N\to\infty$ to a compactly supported measure called the {\it Fuss-Catalan distribution}. The largest square singular value  behaves asymptotically like $\mu_N\sim  \theta \prod_{j=1}^m N_j$ where $\theta >0$ and we consider the rescaled point process   
It becomes clear from their results that one should consider the rescaled square singular value process
\begin{equation} \label{M_Ginibre}
\Xi_N = \sum_{k=1}^N \delta_{\mu_k/M(N)} \hspace{.7cm} \text{where} \hspace{.7cm} M(N)= \prod_{j=1}^m N_j .
\end{equation}

In section~\ref{sect:LLN}, we will prove that the process $\Xi_N$ satisfies a law of large number analogous of (\ref{equilibrium}) and show that, under the assumption of~\cite{AIK_13}, the equilibrium measure is the Fuss-Catalan distribution; see theorem~\ref{thm:Fuss_Catalan}.  The next theorem describes the fluctuations  around this equilibrium configuration.

%The results of~\cite{AIK_13} suggest that should 

 %%% Classical theorem for MOP

\begin{theorem} \label{thm:MOP}
If there exists  $\gamma_1,\dots, \gamma_m \in [0,1]$ such that $N_j/ N \to 1/\gamma_j$ as $N\to\infty$, then
for any polynomial $F\in\R\langle x\rangle$, we have
\begin{equation} \label{CLT_MOP}
\Xi_{N}(F) - \E{\Xi_{N}(F)} \Rightarrow \No\left(0, \sum_{k=1}^\infty k \mathrm{C}^{(\gamma)}_{k}(F) \mathrm{C}^{(\gamma)}_{-k}(F) \right) 
\end{equation}
where, by convention $\gamma_0=1$, and  the coefficients of the variance  are given by
\begin{equation} \label{Fourier'}
\mathrm{C}^{(\gamma)}_k(F) = \frac{1}{2\pi i} \oint F\bigg(z^{-m}   \prod_{l=0}^m\big( z+ \gamma_l)\bigg) z^{-k} \frac{dz}{z} . 
 \end{equation}
\end{theorem}

\proof This CLT is a direct consequence of proposition~\ref{thm:KZ_2} below and theorem~\ref{thm:CLT}. \qed\\

This provides a new class of CLTs for smooth statistics of random matrices which generalize theorem~\ref{thm:OPE}.  Indeed, if $m=1$, $W_N^{*} W_N$ is a Wishart matrix and, if $\gamma_1=1$, the eigenvalues process $\Xi_N$ is described by the Laguerre ensemble which fits in the framework of section~\ref{sect:OPE}. 
 Although, as it was pointed out in remark~\ref{rk:non-Hermitian}, it is not obvious how one would generalize the CLT (\ref{Fourier'}) to test functions which are not analytic and it is a problem that we do not address here. 
%It would be interesting to see whether theorem~\ref{thm:MOP} is generic for multiple orthogonal polynomials which satisfy an $m$-term recurrence relation.  
The rest of this section gives an account on the biorthogonal structure of the ensemble (\ref{AIK}), the corresponding    
recurrence coefficients and their asymptotics. The result are summarized in proposition~\ref{thm:KZ_2} at the end of this section.
We follow the ideas of \cite{KZ_14}, although our  definitions are slightly different and  we obtain the asymptotics of the  recurrence coefficients in a more general case. \\

%Theorem~\ref{thm:MOP} is a direct application of theorem~\ref{thm:CLT}

%given polynomial $q_0, q_1,\dots$ such that $\text{deg } q_n = n$,  define
%$$ Q_n (x) = \int  q_n(s)  \prod_{l=1}^m \Gamma(s+ \nu_l) x^{-s} ds \ . $$

 %For simplicity, suppose that $X_1, \dots, X_m$ are independent $N\times N$ Ginibre matrices and let $(\xi_1,\dots, \xi_N)$ be the eigenvalues of the matrix
%$ X_m^* \cdots X_1^*X_1\cdots X_m$

Given a sequence of polynomials $q_0, q_1,\dots$ such that $\text{deg } q_n = n$,  define
\begin{equation} \label{Q_0}
 Q_n (x) =  \frac{1}{2\pi i} \int_L  q_n(s)  \prod_{l=1}^m \Gamma(s+ \eta_l) x^{-s} ds  . 
 \end{equation}
where $L= \{ c + i y, y\in\R \}$ with $c>0$. By performing linear combinations in the columns of the determinants of formula (\ref{AIK}), we can rewrite
\begin{equation} \label{AIK'}
\varrho_N(x_1,\dots,x_N) = \frac{1}{Z_N'}  \det_{N\times N}\big[ P_{j-1}(x_i) \big] \det_{N\times N}\big[ Q_{j-1}(x_i) \big]  ,
\end{equation} 
where $P_0,P_1,\dots$ are polynomials such that  $\text{deg } P_n = n$ and $\{P_n, Q_n\}_{n\ge 0}$ is a biorthogonal family (cf.~definition~\ref{ensemble}).
Then, it turns out that  $Z_N'=N!$ and the point process  $(\mu_1,\dots, \mu_N)$ is determinantal with correlation kernel $K_N(x,y)=\sum_{n=0}^{N-1} P_n(x)Q_n(y)$ on $L^2(\R_+)$, cf.~formula (\ref{biorthogonal_1}) below.  A natural choice is given by $q_0=1$ and for any $n\ge1$,
\begin{equation} \label{Q_1}
 q_n(s) =  (s-1)\cdots(s-n)= \frac{\Gamma(s)}{\Gamma(s-n)} ,
 \end{equation}
 where $\Gamma(s)$ denotes the Euler gamma function.
Then, the functions $Q_0,Q_1,\dots$ satisfy the Rodriguez's type formula: 
$$ Q_n (x) = \left(-\frac{d}{dx} \right)^n\big\{x^nQ_0(x) \big\}    \ . $$
Moreover, this also provides an integral formula for the Polynomial $P_n$. For any $x>0$,
\begin{equation} \label{P_1}
 P_n(x) =  \frac{1}{2\pi i}  \oint_{\Sigma_n} \frac{\Gamma(z-n)}{ \prod_{l=0}^m \Gamma(z+ \eta_l+1) } x^z dz  . 
 \end{equation}
The contour $\Sigma_n$ is closed around the poles of the function $\Gamma(z-n)$ at $0,1,\cdots, n$ in the region $\Re z> -1$,  so that  by computing the residues, we see that $P_n$ is indeed a polynomial of degree $n$. Note that, since $\eta_0=0$, the integrand has no other pole.  The fact that  $\{P_n, Q_n\}_{n\ge 0}$ is a biorthogonal family can easily be verified by adapting the proof of proposition~\ref{thm:KZ_1}. 
% checked below (proposition~\ref{thm:PQ}).
  Since $P_0,\dots, P_{n+1}$ is a basis of $\R_{n+1}\langle x\rangle$,  we have
\begin{equation} \label{KZ_0}
 xP_n(x) = \alpha_{n,1} P_{n+1}(x) + \sum_{k=0}^n \alpha_{n,-k} P_{n-k}(x)  . 
\end{equation}
where 
 $$ \alpha_{n,-k} = \int_0^\infty xP_n(x) Q_{n-k}(x) dx . $$
% The following explicit formula was first derived~\cite{KZ_14} using a different representation. 

\begin{proposition}\label{thm:KZ_1}
 Let $n \ge m$.  For any $ k>m$, we have  $ \alpha_{n, -k} = 0$ and  if $ -1 \le k \le m$,   
\begin{equation} \label{KZ_1}
 \alpha_{n,-k}  =  \frac{1}{(k+1)!}  \sum_{i=0}^{k+1} (-1)^{i} {k+1 \choose i}  \prod_{l=0}^m (n+\eta_l-i+1)  .
\end{equation}
\end{proposition}

\proof Using the Mellin-Barnes inversion formula, by formula (\ref{Q_0}), when $\Re z>-1$,
\begin{equation}\label{Mellin_inversion}
 \int_0^\infty  x^{z} Q_n (x) dx =  q_n(z+1)   \prod_{l=1}^m \Gamma(z+ \eta_l+1) . 
 \end{equation}
Using formula (\ref{Q_1}) and the functional equation of the gamma function, since $\eta_0=0$, it implies that if $\Re z>-1$, 
$$ \int_0^\infty  x^{z+1} Q_n (x) dx =  \frac{1}{\Gamma(z+2-n)} \prod_{l=0}^m (z+\eta_l+1)  \Gamma(z+ \eta_l+1)  , $$
According to formula (\ref{P_1}), we obtain for any $k\ge -1$,
\begin{align} \int_0^\infty  xP_{n+k}(x) Q_{n} (x) dx 
&\notag=   \frac{1}{2\pi i} \oint_{\Sigma_n}
  \prod_{l=0}^m (z+\eta_l+1)  \frac{ \Gamma(z-n-k)  }{\Gamma(z+2-n)} dz  \\
\alpha_{n+k,-k}&\label{Q_3} =  \frac{1}{2\pi i}  \oint_{\Sigma_n}  \frac{  \prod_{l=0}^m (z+\eta_l+1)}{(z-n+1)\cdots(z-n-k)} dz  .
\end{align}
If $ k >m$, the integrand is $O_{z\to\infty}(z^{-2})$ and all its poles lie inside $\Sigma_n$. Hence, we can move the contour to $\infty$ and the integral vanishes.  Otherwise, by the residue theorem,  we conclude that  for any $ -1 \le k \le m$,  
\begin{align}
\alpha_{n+k,-k}
& \notag = \sum_{j=0}^{k+1}   \prod_{l=0}^m (n+j+\eta_l)   \prod_{\begin{subarray}{c} i=0 \\ i \neq j \end{subarray}}^{k+1} \frac{1}{j-i}\\
&\label{Q_4} =  \sum_{j=0}^{k+1} \frac{(-1)^{k+1-j}}{j! (k+1-j)!}     \prod_{l=0}^m (n+j+\eta_l) .   
\end{align}
 We obtain formula (\ref{KZ_1}) by making the change of variables $n = n'-k$  (so that $\alpha_{n+k,-k} = \alpha_{n',-k}$) and $i= k+1-j$.    \qed\\

We have seen that $(\mu_1,\dots, \mu_N)$ is a determinantal process with kernel $K_N(x,y)=\sum_{n=0}^{N-1} P_n(x)Q_n(y)$. It follows that the rescaled process $\Xi_N$, normalized as (\ref{M_Ginibre}), has correlation kernel
\begin{equation} \label{M_Ginibre'}
\tilde{K}_N(x,y)=\sum_{n=0}^{N-1} \tilde{P}^N_n(x)\tilde{Q}^N_n(y)  ,
\end{equation}
where  
$$\tilde{P}_n^N = \kappa_n(N)  \sqrt{M} P_n\big( x  M\big) 
\hspace{.5cm}\text{and}\hspace{.5cm}
 \tilde{Q}_n^N = \kappa_n(N)^{-1}  \sqrt{M} Q_n\big( x M\big) . $$
  Then,  $\{\tilde{P}_n^N, \tilde{Q}_n^N\}_{n\ge 0}$ is still a biorthogonal family on $L^2(\R_+)$ and, by formula (\ref{KZ_0}), the polynomials  $\tilde{P}_0^N, \tilde{P}_1^N,\dots$  satisfy a recurrence relation of the form (\ref{recurrence}) where the matrix $\J$ is given by  
\begin{equation} \label{KZ_2}
\J_{n,n-k} = \frac{\kappa_n(N)}{\kappa_{n-k}(N)} \frac{ \alpha_{n,-k} }{M} \1_{-1 \le k\le m } . 
\end{equation}

By proposition~\ref{thm:KZ_1}, $\J$ has at most $m+1$ non-zero diagonals and, to obtain a limit theorem, it remains to check whether it  has a right-limit as $N\to\infty$. 
We will need the following elementary combinatorial lemma.

\begin{lemma} \label{Binomial_identity}
 For any $n\in\N_0$ and for any polynomial $R(x)$ of degree $\le n$,  
$$ \sum_{k=0}^n (-1)^{n-k}  {n\choose k } R(k) = R^{(n)}(0)  . $$
\end{lemma}

\proof Let $0 \le r \le n$, if we differentiate the binomial identity $r$ times, we obtain
$$  {n\choose r} (1+x)^{n-r} =   \sum_{k=0}^n {n\choose k } {k\choose r} x^{k-r}  ,$$
and, if we evaluate at $x=-1$, this implies that 
\begin{equation} \label{combinatorics_1}
 \sum_{k=0}^n (-1)^{k-r} {n\choose k } {k\choose r} = \delta_{r,n}  .
 \end{equation}
Since  $\left\{ \displaystyle{x\choose r} \right\}_{r=0,\dots, n}$ is a basis of $\R_n[\mathrm{X}]$ and  
 $\displaystyle R(x)=  R^{(n)}(0) {x\choose n} + \cdots $ , the lemma follows from equation (\ref{combinatorics_1}). \qed\\

\begin{proposition} \label{thm:KZ_3}
Let $\kappa_n(N)= N^{n}$ in formula (\ref{KZ_2}). Suppose that 
 $N_l/N = 1+ \eta_l/N \to 1/\gamma_l$ where $\gamma_1,\dots, \gamma_m \in [0,1]$
and there exists $t \ge 0$ such that $n/N \to t$ as $N\to\infty$.    
 Then, for all  $ -1 \le k \le m$,
  \begin{equation} \label{KZ_3}
\J_{n,n-k} = \sum_{\begin{subarray}{c} S \subseteq \{0,\dots, m\} \\ |S| = k+1 \end{subarray}}   \prod_{l \in S} \gamma_l   \prod_{l \notin S} \big(1+ \gamma_l(t-1)\big)
 + \underset{N\to\infty}{O}\big(N^{-1}\big)  . 
\end{equation}
\end{proposition}

\proof By formula \eqref{KZ_1},  we have for any  $ -1 \le k \le m$,
\begin{equation*}
 \alpha_{n,-k}  =  \frac{1}{(k+1)!}  \sum_{j=0}^{k+1} (-1)^{j} {k+1 \choose j}   \sum_{\begin{subarray}{c} S \subseteq \{0,\dots, m\} \end{subarray}}  (1-j)^{|S|}
  \prod_{l \notin S} (n + \eta_l) .
\end{equation*}
Applying lemma~\ref{Binomial_identity} to the polynomials $ (1-x)^{|S|}$, we obtain 
\begin{equation*}
 \alpha_{n+j,-k}  = 
  \sum_{\begin{subarray}{c} S \subseteq \{0,\dots, m\} \\ |S| = k+1 \end{subarray}}     \prod_{l \notin S} (n + \eta_l)  + \underset{N\to\infty}{O} \bigg(  \sum_{\begin{subarray}{c} S \subseteq \{0,\dots, m\} \\ |S| >k+1 \end{subarray}}  
  \prod_{l \notin S} (n + \eta_l)  \bigg) . 
\end{equation*}
Since $M= \prod_{l=1}^m N_l$, this implies that
\begin{equation*}
N^k \frac{ \alpha_{n+j,-k} }{M}
=  \sum_{\begin{subarray}{c} S \subseteq \{0,\dots, m\} \\ |S| = k+1 \end{subarray}} \prod_{l \in S} \frac{N}{N_l}       \prod_{l \notin S} \frac{n + \eta_l}{N_l}  + \underset{N\to\infty}{O} \big( N^{-1}  \big) . 
\end{equation*}
In particular, the error term follows from the fact that in the regime we consider:
$ \frac{n + \eta_l}{N_l} \to \gamma_l t + (1-\gamma_l)$ as $N\to\infty$. 
Then, the asymptotics \eqref{KZ_3} follow directly from formula  \eqref{KZ_2}. \qed\\ 

Since $\J_{n,n-k}=0$  if $k>m$ or $k<-1$ for all $n>0$,  proposition~\ref{thm:KZ_3} with the parameter $t=1$ shows that the matrix $\J$ has a right-limit  which is a Laurent matrix whose symbol depends explicitly on the parameters $(\gamma_1, \dots, \gamma_m)$. The next proposition summarizes our results.

\begin{proposition}\label{thm:KZ_2}  Let  $\mu_1,\dots, \mu_N$ be the square singular values of the matrix  $(\ref{W_matrix})$. The rescaled point process $(\ref{M_Ginibre})$ is a biorthogonal ensemble with correlation kernel
$(\ref{M_Ginibre'})$ where  $\tilde{P}^N_0, \tilde{P}^N_1,\dots$ are polynomials which satisfies a recurrence relation 
of the form $(\ref{recurrence})$. Moreover, the recurrence matrix $\J$  has a right-limit $\L(s)$ as $N\to\infty$ and $N_j/N \to 1/\gamma_j$ which is a Laurent matrix with symbol
 $ s(z) =   z^{-m}  \prod_{l=0}^m\big( z+ \gamma_l)$.
\end{proposition}

\subsection{The Jacobi and Laguerre Muttalib-Borodin ensembles} \label{sect:BME}

The Muttalib-Borodin ensembles have been introduced in~\cite{Muttalib_95}  to give a better description of the physics of disordered conductors in the metallic regime than  the classical random matrix theory. 
 Recently they have received  renewed attention in the mathematics literature, e.g.~\cite{Cheliotis_14, CR_14, FW_15, Zhang_15}. 
In particular, it was proved in~\cite{Cheliotis_14, FW_15} that these ensembles can be realized as the eigenvalue  density functions of certain random matrices. 
Muttalib argued that the distribution of the energy levels of such conductors should rather be modeled by a Boltzmann-Gibbs ensemble with density $\propto e^{-\beta \mathscr{H}^\theta_V(x_1,\dots, x_N)/2}$,  where the Hamiltonian is
$$  \mathscr{H}^\theta_V(x_1,\dots, x_N)= - \frac{1}{2}\sum_{i\neq j} \log|x_i - x_j| - \frac{1}{2}\sum_{i\neq j} \log|x_i^\theta - x_j^\theta| + N \sum_{i} V(x_i)  , $$
 for some small values of the parameter $\theta>0$.
 It was established in~\cite{Muttalib_95} that there exists polynomials $p_k$ and $q_k$  of degree $k$ (which depend on the parameters $N$ and $\theta$) such that 
for all $n, k \ge 0$,
\begin{equation} \label{OR}
\int  p_n(x) q_k(x^\theta)   e^{-NV(x)} dx = \delta_{n,k}  .
\end{equation}
So, when $\beta=2$, this point process is a biorthogonal ensemble with respect to the measure  $d\mu_N = e^{-NV(x)}dx$ on $\R$ with correlation kernel:
 \begin{equation*}
K_N(x, y) = \sum_{k=0}^{N-1} p_k(x) q_k(y^\theta) .
\end{equation*}

Moreover, if $\theta=1$, we have $q_k=p_k$ and  we recover the classical orthogonal polynomial ensemble $\mathbb{P}_{N,V}$, (\ref{Hamiltonian}).
In the following, we let $\theta = 1/m$ where $m\in\N$. Then, it is not difficult to check that  the orthogonality relation (\ref{OR}) implies that $x p_k(x) \in \text{span}\{ p_{k-m}, \dots, p_{k+1} \}$    
so that the sequence $\{ p_k \}_{k \ge m}$ satisfies an $m$-terms recurrence relation.
Generally, if $\theta \in \mathbb{Q}$, a more complicated recurrence relation still exists, see \cite[Theorem~1.1]{CR_14}, but it is not as useful in the context of this paper. \\

 In the following, we will consider the Laguerre and Jacobi weights investigated in~\cite{Borodin_99}. Borodin obtained explicit formulae for the correlation kernels  and the local limits at the hard edge and in the bulk. In addition, he derived expressions for the corresponding biorthogonal families that we will use to compute the recurrence coefficients and find their asymptotics  using the method introduced in section~\ref{sect:product_Ginibre}. For these  Muttalib-Borodin ensembles, the equilibrium measures has an explicit description~\cite{FW_15}. We shall prove that, when $\theta=1/m$, the asymptotic fluctuations around the equilibrium are also described by theorem~\ref{thm:MOP}. This raise the question whether, like OPEs,  the Muttalib-Borodin ensembles have a universal behavior, which depends only on $\theta$, as $N\to\infty$?   \\

 We begin with the Laguerre ensemble. Let $a>-1$ and define on $(0,\infty)^N$ the Hamiltonian:
 \begin{equation} \label{Laguerre}  \mathscr{H}^\theta_{Lag}(x_1,\dots, x_N)= - \frac{1}{2}\sum_{i\neq j} \log|x_i - x_j| - \frac{1}{2}\sum_{i\neq j} \log|x_i^\theta - x_j^\theta| - \sum_{i}\big( x_i - a \log (x_i) \big) .
 \end{equation}
We let $(\mu_1,\dots, \mu_N)$ be the corresponding point process. Its correlation kernel (with respect to the Lebesgue measure on $\R_+$)  is of the form
\eqref{kernel_0} and the biorthogonal family  is given by for all $x>0$, 
\begin{equation} \label{Laguerre_1} P_k(x) = \frac{1}{k!}  Z^a_k(x^\theta)
%\sum_{j=0}^n {n \choose j } \frac{(-1)^j x^{\theta j}}{ \Gamma(\theta j + a +1) } 
%&=  \frac{1}{2\pi i} \oint_{\Sigma_m} \frac{1 }{  \prod_{k=0}^{m}(z-k)\Gamma(a+1+z\theta)} x^{\theta z} dz \\
=\frac{1}{2\pi i} \oint_{\Sigma_k} \frac{\Gamma(z-k) }{  \Gamma(z+1)\Gamma(a+1+z\theta)} x^{\theta z} dz 
\end{equation}
and 
\begin{equation} \label{Laguerre_2}
 Q_k(x) = k! Y^a_k(x) x^a e^{-x} =   \frac{1}{2\pi i} \int_{L} \frac{ \Gamma( (w-1)\theta^{-1} + 1)}{ \Gamma( (w-1)\theta^{-1} + 1-k)} \Gamma(w+a) x^{-w} dw .
 \end{equation}
The contour $\Sigma_k$ is a closed loop around the points $0,1,\cdots, k$ and  $L= \{ c + i y, y\in\R \}$ with $c>0$. The functions $Z^a_k$ and $Y^a_k$ are polynomials of degree $k$  whose expression is given in~\cite[p724]{Borodin_99} and it is not difficult to derive formulae (\ref{Laguerre_1} - \ref{Laguerre_2}), see for instance  proposition~2.2  in~\cite{Zhang_15}. We notice that this biorthogonal family has the same structure as the family which arises  from describing the square-singular values process of product of Ginibre matrices. 
In fact, when $\theta \in \N$, an  explicit  connection was given in~\cite[sect.~3.2]{Zhang_15}.
%In particular, the method developed in section~\ref{sect:product_Ginibre} can be used to analyze the Laguerre Muttalib-Borodin ensemble as well.
According to formula (\ref{Laguerre_2}), the Mellin-Barnes inversion formula implies that for $\Re v >-1$,
$$\int_0^\infty x^v Q_n(x)dx =\frac{ \Gamma(v\theta^{-1} + 1)}{ \Gamma( v\theta^{-1} + 1-n)} \Gamma(v+1+a) .
$$
So, taking $v= \theta z+1$ with $\theta=1/m$, we obtain
\begin{equation}\label{Laguerre_3}
\int_0^\infty x Q_n(x) P_k(x)dx =  \frac{1}{2\pi i} \oint_{\Sigma_n} \frac{ \Gamma(z-k) \Gamma(z+m + 1)}{\Gamma( z + 1) \Gamma( z + m +1-n)} (z/m+1+a)  dz . 
  \end{equation}
Proceeding exactly like in the proof of proposition~\ref{thm:KZ_1}, we deduce from formula  (\ref{Laguerre_3}) that

\begin{equation} \label{Laguerre_4}
 xQ_n(x) =\sum_{k=-1}^n \alpha_{n,-k} Q_{n-k}(x)  . 
\end{equation}
where, if we let $\eta_0=0, \eta_1=1, \dots, \eta_{m-1}= m-1$ and $\eta_m =- m a$, we have 
\begin{equation} \label{Laguerre_5}
 \alpha_{n,-k}  =  \frac{1}{m} \sum_{i=0}^{m-k} \frac{(-1)^{m-k -i}}{{(m-k)!} } {m-k \choose i}  \prod_{l=0}^m (n+i- \eta_l)  .
\end{equation}

By lemma~\ref{Binomial_identity}, we can check that formula (\ref{Laguerre_5}) implies that for any $-1 \le k\le m$, the recurrence coefficients have the following asymptotics:
\begin{equation} \label{Laguerre_6}
 \alpha_{n,-k}  =  \frac{1}{m} {m+1 \choose m-k} n^{k+1} + \underset{n\to\infty}{O}(n^k)  . 
\end{equation}  

In the global regime, the form of the Hamiltonian (\ref{Laguerre}) suggests that we should consider the rescaled process $(\mu_1/N,\dots, \mu_N/N)$ in order to make the potential competitive, its correlation kernel is given by
\begin{equation*}
K^{Lag}_N(x,y)=\sum_{n=0}^{N-1} \tilde{P}_n^N(x)\tilde{Q}_n^N(y) ,
\end{equation*}
where  $\tilde{P}^N_n = \kappa_n(N)\sqrt{N} P_n\big( x N\big)  $ and $\tilde{Q}^N_n = \kappa_n(N)^{-1}  \sqrt{N} Q_n\big( xN\big) $. Given the asymptotics (\ref{Laguerre_6}), we choose the normalization  $\kappa_n(N)= N^{n}$ so that, according to formula (\ref{Laguerre_4}), the recurrence matrix $\J$ for the functions  $\tilde{Q}_0^N, \tilde{Q}_1^N,\dots$  satisfies
\begin{equation*}
\J_{n,n-k} =  N^{-k-1}\alpha_{n,-k} \1_{-1 \le k\le m } . 
\end{equation*}

This shows that $\J$ has a right-limit as $N\to\infty$ which is a Laurent matrix $\mathbf{L}(s)$ with symbol
\begin{equation*}  \label{Laguerre_7}   s(z)=  \frac{1}{m}  \sum_{k=-1}^m {m+1 \choose m-k}  z^{-k} = 
\frac{ z^{-m}  (z+1 )^{m+1}}{m} .  
\end{equation*}

 This should be compared with proposition~\ref{thm:KZ_2}.
According to remark~\ref{rk:symbol}, we see that,  up to a dilation, the fluctuations of linear statistics of the Laguerre Muttalib-Borodin ensemble with parameter $\theta=1/m$ and of the square-singular values  of the product of $m$ rectangular Ginibre matrices in the regime where the dimensions satisfy $N_j/N \to 1$ as $N \to \infty$ for $j=1,\dots, m$ are described by the same CLT. 
In particular, we obtain the following corollary of theorem~\ref{thm:CLT}.  

 \begin{theorem} \label{thm:Laguerre}
 Let $\Xi_N= \sum_{j=1}^N  \delta_{ m \mu_j/  N}$ where $(\mu_1, \dots, \mu_N)$ is a point process with density $\propto e^{- \mathscr{H}^{1/m}_{Lag}(x_1,\dots, x_N)}$ on $\R_+$, $(\ref{Laguerre})$. 
For any polynomial $F\in\R\langle x\rangle$, we have as $N\to\infty$, 
\begin{equation*} 
\Xi_{N}(F) - \E{\Xi_{N}(F)} \Rightarrow \No\left(0, \sum_{k=1}^\infty k \mathrm{C}^{(\bf 1)}_{k}(F) \mathrm{C}^{(\bf 1)}_{-k}(F) \right) \ ,
\end{equation*}
where the coefficients $\mathrm{C}^{(\bf 1)}_k(F) $ are given by formula $(\ref{Fourier'})$ with parameters $\gamma_0=\cdots=\gamma_m=1$. 
\end{theorem}

In the remainder of this section, we will consider another BOE defined  on $(0,1)^N$ by the Hamiltonian:
 \begin{equation*} \label{Jacobi}  \mathscr{H}^\theta_{Jac}(x_1,\dots, x_N)= - \frac{1}{2}\sum_{i\neq j} \log|x_i - x_j| - \frac{1}{2}\sum_{i\neq j} \log|x_i^\theta - x_j^\theta|  -(a-1) \sum_{i}\log (x_i)   .
 \end{equation*}
 The parameter  $a>0$ and we let $\theta=1/m$ for $m\in \N$. We denote by $(\nu_1,\dots, \nu_N)$ this point process and call it the Jacobi Muttalib-Borodin ensemble since it was introduced in~\cite[sect.~3]{Borodin_99}. For the general Jacobi weight, explicit formulae for the biorthogonal family  are also available, see~\cite{MT_82}. However, computing the asymptotics of the recurrence coefficients turn out to be more sophisticated and we will not formalize the argument. 
  Although the Jacobi potential is very different than the Laguerre potential, up to a scaling, the fluctuations of linear statistics of both BOEs are described by the same Gaussian process.

  \begin{theorem} \label{thm:Jacobi}
 Let $\beta_m = m(1+ 1/m)^{m+1}$ and $\Xi_N= \sum_{j=1}^N  \delta_{\nu_j/ \beta_m }$ where $(\nu_1, \dots, \nu_N)$ is the point process with density $\propto e^{- \mathscr{H}^{1/m}_{Jac}(x_1,\dots, x_N)}$.
For any polynomial $F\in\R\langle x\rangle$, we have as $N\to\infty$, 
\begin{equation*} 
\Xi_{N}(F) - \E{\Xi_{N}(F)} \Rightarrow \No\left(0, \sum_{k=1}^\infty k \mathrm{C}^{(\bf 1)}_{k}(F) \mathrm{C}^{(\bf 1)}_{-k}(F) \right) \ ,
\end{equation*}
where the coefficients $\mathrm{C}^{(\bf 1)}_k(F) $ are given by formula $(\ref{Fourier'})$ with $\gamma_0=\cdots=\gamma_m=1$. 
\end{theorem}

The remainder of this section is devoted to the proof of theorem~\ref{thm:Jacobi}. According to proposition~3.3 in~\cite{Borodin_99}, the  Jacobi Muttalib-Borodin ensemble has correlation kernel
 \begin{equation} \label{Jacobi_0}
K_N^{Jac}(x, y) = \sum_{n=0}^{N-1} \zeta_n(x) \psi_n(y)  ,
\end{equation}
where the biorthogonal family is given by 
 \begin{align} \notag \zeta_n(x) 
 &= \sum_{j=0}^n \frac{\prod_{k=0}^{n-1} (j + \theta k +a ) }{\prod_{k \neq j}(j-k)}  x^{j+a-1} \\ 
&  \label{Jacobi_1} =  \frac{ 1}{2\pi i} \oint_{\Sigma_n}   \frac{\prod_{j=0}^{n-1} (z + \theta j +a ) }{\prod_{j=0}^{n}(z-j)}  x^{z+a-1} dz  ,
\end{align}
and
 \begin{align} \notag 
 \psi_n(x) &=  \kappa_n \sum_{j=0}^n \frac{\prod_{k=0}^{n-1} (\theta  j +  k +a ) }{\theta^n \prod_{k \neq j}(j-k)}  x^{\theta j} \\
& \label{Jacobi_2} 
= \frac{ \kappa_n}{2\pi i} \oint_{\Sigma_n'}   \frac{\prod_{j=0}^{n-1} (\theta w +  j +a ) }{ \theta^n \prod_{j=0}^{n}(w-j)}  x^{\theta w} dw .
\end{align} 
The contour $\Sigma_n$ and $\Sigma_n'$ are closed loops around the points $0,1,\cdots, n$ in the region $\Re z > - a $, respectively $\Re w > - 1/\theta $, 
and $\kappa_n =  n(\theta+1)+a$. 
By the general theory, the functions $\psi_n$ satisfy a recurrence relation:
 \begin{equation}  \label{Jacobi_3'}
 x\psi_n(x) =\sum_{k=-1}^n \alpha_{n,k} \psi_{n+k}(x)  , 
\end{equation}
where
\begin{equation} \label{Jacobi_3}
\alpha_{n, k}= \int_0^1 x  \psi_n(x)  \zeta_{n+k}(x) dx  .
\end{equation}

\begin{lemma} \label{thm:Jacobi_recurrence}
If we denote $\eta_0= m a, \eta_1=1, \dots, \eta_{m}= m$, then
\begin{equation} \label{Jacobi_6}
\alpha_{n, k}=  \frac{ m \kappa_n }{ (m-k )!} \sum_{r=0}^{m-k}  (-1)^{r} {m -k \choose r } \frac{\prod_{j=0}^{m} (n-r + \eta_j ) }{\prod_{j=0}^{k+1}\big(n(m+1) -r+ m(j+a)\big)}  . 
\end{equation}
\end{lemma}

\proof
According to formulae (\ref{Jacobi_1} - \ref{Jacobi_2}), by Fubini's theorem, we have
\begin{equation}  \label{Jacobi_4}
 \int_0^1 x  \psi_n(x)  \zeta_{n+k}(x) dx  = \frac{  \kappa_n \theta^{-n}}{2\pi i } 
 \oint_{\Sigma_n'}   \frac{\prod_{j=0}^{n-1} (\theta w +  j +a ) }{ \prod_{j=0}^{n}(w-j)} \Upsilon_{n+k}(w)dw
 \end{equation}
where 
\begin{align*}  
\Upsilon_{n}(w) &=  \frac{1}{2\pi i} \oint_{\Sigma_n}   \frac{\prod_{j=0}^{n-1} (z + \theta j +a ) }{\prod_{k=0}^{n}(z-j)} \int_0^1  x^{z+a+\theta w} dx dz \\
& 
=  \frac{-1}{2\pi i} \oint_{\Sigma_n}   \frac{\prod_{j=0}^{n-1} (z + \theta j +a ) }{\prod_{j=0}^{n}(z-j)} \frac{dz}{z+a+\theta w +1}  . 
\end{align*}

In the previous formula, the integrand is of order $O_{z\to\infty}(z^{-2})$ and it is analytic outside of $\Sigma_n$ except for a simple pole located at  $z=-a-1 - \theta w$. Hence,  since $\theta= 1/m$, this implies that
\begin{equation} \label{Jacobi_5}
\Upsilon_{n}(w)
= \theta^n   \frac{\prod_{j=0}^{n-1} (w + m-j ) }{\prod_{j=1}^{n+1}( \theta w+ j+ a  )} .
\end{equation}

Then, by formulae (\ref{Jacobi_3}- \ref{Jacobi_5}), we obtain for any $n + k > m$,  
\begin{equation*} 
\alpha_{n, k}  = \frac{ \kappa_n \theta^k}{2\pi i}
  \oint_{\Sigma_n'}  \frac{  \prod_{j=0}^{n-1} (\theta w +  j +a )}{\prod_{j=1}^{n+k+1}(\theta w+j+ a)}   \frac{  \prod_{j=0}^{n+k-m -1}(w-j) }{ \prod_{j=0}^{n}(w-j)}  \prod_{j=1}^{m} (w + j ) dw  .
\end{equation*}

First note that, if $ k > m$, the integrand has no pole inside the contour $\Sigma_n'$ and $\alpha_{n, k} = 0$. On the other hand, if $k < -1$, since the integrand has no pole outside of the contour $\Sigma_n'$ and is of order $O_{w\to\infty}(w^{-2})$, then $\alpha_{n, k} = 0$ as well.  Finally, when $ -1 \le k \le m$, a residue computation yields
\begin{align*}
\alpha_{n, k}
 & = \frac{  \kappa_n \theta^{k}}{2\pi i}
  \oint_{\Sigma_n'}  \frac{ \theta w+ a }{\prod_{j=n}^{n+k+1}(\theta w+j+ a)}   \frac{   \prod_{j=1}^{m} (w + j )  }{ \prod_{j=n-m +k}^{n}(w-j)}dw \\
  &= \frac{ \kappa_n \theta^k }{ (m-k )!} \sum_{r=0}^{m-k}  (-1)^{r} {m -k \choose r } \frac{ \big( \theta( n-r) + a \big)\prod_{j=1}^{m} (n-r + j ) }{\prod_{j=n}^{n+k+1}\big(\theta (n-r)+ j+a\big)}  .
\end{align*}
Since $\theta=1/m$, this completes the proof. \qed\\

Based on lemma~\ref{thm:Jacobi_recurrence}, we can compute the asymptotics of recurrence coefficients $\alpha_{N,k}$ as $N\to\infty$. By  formula (\ref{Jacobi_3'}), the functions $\psi_n$ satisfy a recurrence relation of the form (\ref{recurrence}) where  
$$\J_{n,n+k} =\alpha_{n,k} \1_{-1 \le k\le m} $$ 
and, according to proposition~\ref{thm:Jacobi_asymptotics} below, 
this matrix $\J$ has a right-limit $\L(s)$ which is a Laurent matrix with symbol 
$$ s(z) = \sum_{k=-1}^m {m+1 \choose k+1}   \frac{m^{m-k}}{ (m+1)^{m+1}} z^{k} = z^{-1} \left( \frac{z+m}{m+1}\right)^{m+1}  .  $$ 

Thus, if we let $\beta_m = m(1+ 1/m)^{m+1}$, we see that
$$ s\left( \frac{m}{z}\right)  = \frac{z^{-m}(z + 1)^{m+1}}{\beta_m}  . $$
By remark~\ref{rk:symbol}, this implies that  the fluctuations of the rescaled Jacobi Muttalib-Borodin
$(\frac{\nu_1}{ \beta_m}, \dots, \frac{\nu_N}{\beta_m})$ are described by the Laurent matrix $\L(\tilde s)$ where 
$\tilde s(z)= z^{-m} ( z + 1)^{m+1}$. By theorem~\ref{thm:CLT}, this completes the proof of theorem~\ref{thm:Jacobi}. 
It remains to compute the asymptotics of the recurrence coefficients~(\ref{Jacobi_6}).

\begin{proposition} \label{thm:Jacobi_asymptotics}
For any $ -1 \le k \le m$, we have 
$$  \lim_{n\to\infty} \alpha_{n, k} = {m+1 \choose k+1}   \frac{m^{m-k}}{ (m+1)^{m+1}}  . $$
\end{proposition}

\proof  When $n$ is sufficiently large,  for all $0 \le r\le m+1$, 
$$ \frac{1}{\prod_{l=0}^{k+1}\big(n(m+1) -r+ m(l+a)\big)}  
= \sum_{j=0}^\infty \frac{n^{-k-2- j}}{(m+1)^{k+2+j}}
  \underbrace{\sum_{\begin{subarray}{c} j_0 +\cdots + j_{k+1}= j \\ j_0, \dots, j_{k+1} \ge 0 \end{subarray}} \prod_{l=0}^{k+1} \big( r- m(l+a) \big)^{j_l}}_{:= R_j(r)}  .
$$
In particular, $R_j(r)$ is a polynomial of degree $j$ whose leading coefficient is the number of nonnegative integer solutions of the equation $j_0 +\cdots + j_{k+1}= j $ which equals to $ {j + k +1 \choose j}$, see \cite[p~25]{Stanley_12}.
According to lemma~\ref{thm:Jacobi_recurrence}, we obtain 
$$ \alpha_{n, k}=  \frac{ m \kappa_n}{ (m-k )!} \sum_{r=0}^{m-k}  (-1)^{r} {m -k \choose r }
\sum_{j=0}^\infty \sum_{S \subseteq \{0,\dots, m\}} \frac{ n^{m-k-1-j - |S|} }{(m+1)^{k+2+j}} R_j(r) \prod_{i \in S}(\eta_i -r) .  $$

Since by definition, $ \frac{m \kappa_n}{ n} = (m +1) + \underset{n\to\infty}{O}(n^{-1}) $, we can rearrange this formula:
$$ \alpha_{n, k}=  \frac{1}{ (m-k )!} \sum_{j=0}^{m-k}
\sum_{\begin{subarray}{c} S \subseteq \{0,\dots, m\} \\ |S| \le m-k - j \end{subarray}}
 \frac{ n^{m-k-j - |S|} }{(m+1)^{k+j+1}}   
  \sum_{r=0}^{m-k}  (-1)^{r} {m -k \choose r }R_j(r) \prod_{i \in S}(\eta_i -r)  + \underset{n\to\infty}{O}\big( n^{-1}\big)  .$$

By lemma~\ref{Binomial_identity}, since $R_j(r)  \prod_{i \in S}(\eta_i -r)  $ is a polynomial of degree $|S|+n$ with leading coefficients $(-1)^{|S|}  {j + k +1 \choose j}$,  this implies that
$$ \alpha_{n, k}= \sum_{j=0}^{m-k}
\sum_{\begin{subarray}{c} S \subseteq \{0,\dots, m\} \\ |S| = m-k - j \end{subarray}}
 \frac{  (-1)^{j} }{(m+1)^{k+j+1}}   
  {j + k +1 \choose j} + \underset{n\to\infty}{O}\big( n^{-1}\big) , $$
so that
\begin{align*} \lim_{n\to\infty} \alpha_{n, k}
&= \sum_{j=0}^{m-k} \frac{ (-1)^j}{(m+1)^{k+1+j}}  {m+1 \choose m-k - j }  {j + k +1 \choose j} \\
&= {m+1 \choose k+1} \sum_{j=0}^{m-k} \frac{ (-1)^j}{(m+1)^{k+1+j}}  {m-k \choose  j }  .
\end{align*}
Using the binomial formula, this completes the proof. \qed

\section{Proofs of the Main results} \label{sect:proof}
\subsection{Paths formulation and the proof of  lemma~\ref{thm:cumulants}.} \label{sect:path}

 Let $\Xi_N$ be a BOE with a recurrence matrix $\J=\J^{(N)}$. In section~\ref{sect:cumulants}, we have proved that the cumulants of a linear statistic $\Xi_N(F)$, where $F$ is a polynomial, can be expressed in terms of the recurrence matrix $\J$.    
Namely, if we let $\M=F(\J)$, by formula (\ref{Laplace_5}),
\begin{equation} \label{cumulant_0}
\Cu^n_N[F] = - \sum_{\ell=1}^{n} \frac{(-1)^\ell}{\ell} 
\sum_{\begin{subarray}{c}k_1+\cdots+k_\ell = n \\ k_i \ge 1 \end{subarray}} \frac{n!}{k_1!\cdots k_\ell!}
  \tr\big[\p_N \M^{k_1}\p_N \cdots \p_N \M^{k_\ell} \p_N\big]  .
\end{equation}

Since $\J$ is a  band matrix and $F$ is a polynomial, the matrix elements of  $\M$ in the canonical basis of $l^2(\N_0)$ are finite and, by (\ref{band}), $\M$ has at most $\W^{\deg F}$ non-zero diagonals, so that the cumulants $\Cu^n_N[F]$ are well-defined for all $n,N\in\N$.  Our strategy is to express the
traces on the RHS of formula (\ref{cumulant_0}) in term of sum over paths on the  adjacency graph $ \mathcal{G}\big(\M \big)$.  First, we need to introduce further notation.
 For any $n>0$ and $j_1,j_2 \ge 0$, we denote
\begin{equation*}
\Gamma^{n}_{j_1\to j_2} = \{ \text{path }\pi \text{ of length }n \text{ on }\mathcal{G}(\M) \text{ such that } \pi(0)=j_1 \text{ and }\pi(n)= j_2\}   .
\end{equation*}
In general, if  $n= n_1+\cdots + n_\ell$, for any $j_1,\dots, j_{\ell+1} \ge 0 $, we define
\begin{align} \label{path_3}
 \Gamma_{j_1 \to j_2 \to \cdots \to j_{\ell+1}}^{n_1+n_2 + \cdots +n_\ell}
 =\big\{ &\text{path }\pi \text{ of length }n \text{ on }\mathcal{G}(\M) \text{ such that}\\
 & \notag \pi(0)=j_1,\ \pi(n_1)=j_2,\ \pi(n_1+n_2)=j_3,\  \cdots,\ \pi(n_1+\cdots + n_\ell)=j_{\ell+1} \big\} .
\end{align}

Recall that the edges of  $\mathcal{G}(\M)$ are oriented and the weight of the edge $e= (i,j)$ is given by $\M_e=\mathbf{M}_{ij}$.
Then, if $\pi$ is a path on the graph $\mathcal{G}(\M)$, we denote
$\displaystyle\M\{\pi\} = \prod_{e\in\pi} \M_e$ and for any collection $\Gamma$ of paths, we let
\begin{equation} \label{path_4}
 \M\big\{ \Gamma  \big\} =  \sum_{\pi \in \Gamma } \M\{\pi\} . 
 \end{equation}
In particular, if we denote by $ \M^{n}_{ij}$ the entries of the matrix $\M^n$, we get 
\begin{equation} \label{path_0}
 \M^{n}_{ij} = \sum_{\pi \in \Gamma_{i\to j}^n } \M\{\pi\}= \M\big\{\Gamma^n_{i\to j}\big\}  .
 \end{equation}
  For instance, this implies that for any $n\in\N$,
 \begin{equation} \label{trace_1}
  \tr\big[\p_N \M^{n}\p_N\big] = \sum_{j=0}^{N-1} \M^{n}_{jj}  = \sum_{j=0}^{N-1} \M\big\{\Gamma^n_{j\to j}\big\} . 
 \end{equation}
If $n=1$, this gives the expression of the first cumulant $\Cu^1_N[F] =\E{\Xi_N(F)}$. We can express in a  similar fashion
\begin{equation}  \label{trace_2}
  \tr\big[\p_N\M^{k_1}\p_N\cdots\p_N \M^{k_\ell} \p_N\big]  
  =\sum_{j_1=0}^{N-1}\cdots\sum_{j_\ell=0}^{N-1} \prod_{i=1}^\ell \M^{k_i}_{j_{i}j_{i+1}} ,
\end{equation} 
 where  $j_{\ell +1} = j_1$ by convention. 
If $\pi \in \Gamma_{i\to s}^{k_1}$ and $\kappa \in \Gamma_{s\to j}^{k_2}$ are two paths, then we define a new path $\pi\oplus\kappa\in \Gamma_{i\to s\to j}^{k_1+k_2}$ which is the concatenation of $\pi$ and $\kappa$. Moreover, by definition of the weight $\M\{\cdot\}$, we see that $\M\{\pi\oplus\kappa\}=\M\{\pi\}\M\{\kappa\}$ Conversely, any path $\nu\in \Gamma_{i\to s\to j}^{k_1+k_2}$ has a unique decomposition $\nu=\pi\oplus\kappa$ such that  $\pi \in \Gamma_{i\to s}^{k_1}$ and $\kappa \in \Gamma_{s\to j}^{k_2}$. Hence, we can write $ \Gamma_{i\to s\to j}^{k_1+k_2} = \Gamma_{i\to s}^{k_1}\oplus \Gamma_{s\to j}^{k_2}$ and if we iterate, we obtain
 $$ \Gamma_{j_1\to j_2}^{k_1}\oplus \Gamma_{j_2\to j_3}^{k_2}\oplus \cdots \oplus  \Gamma_{j_\ell\to j_{\ell+1}}^{k_\ell}
 = \Gamma_{j_1 \to j_2 \to \cdots \to j_{\ell+1}}^{k_1+k_2 + \cdots +k_\ell}  ,
 $$ 
 and
 $$ \M\big\{\Gamma_{j_1\to j_2}^{k_1}\big\}\M\big\{\Gamma_{j_2\to j_3}^{k_2}\big\}\ \cdots\  \M\big\{\Gamma_{j_\ell\to j_{\ell+1}}^{k_\ell}\big\}
 =\M\big\{ \Gamma_{j_1 \to j_2 \to \cdots \to j_{\ell+1}}^{k_1+k_2 + \cdots +k_\ell}\big\}  .
 $$  
By formula (\ref{path_0}), this implies that for any $k\in \N_0^m$, 
 $$\prod_{i=1}^\ell \M^{k_i}_{j_{i} j_{i+1}}
 =\M\big\{ \Gamma_{j_1 \to j_2 \to \cdots \to j_{\ell+1}}^{k_1+k_2 + \cdots +k_\ell}\big\} $$ 
 and, by formula (\ref{trace_2}),  
 \begin{equation}  \label{trace_3}
  \tr\big[\p_N\M^{k_1}\p_N\cdots\p_N \M^{k_\ell} \p_N\big]  
  = \sum_{j_1=0}^{N-1}\cdots\sum_{j_\ell=0}^{N-1}  \M\big\{ \Gamma_{j_1 \to j_2 \to \cdots \to j_{\ell+1}}^{k_1+k_2 + \cdots +k_\ell}\big\}    .
 \end{equation}
Observe that, if if $k_1 + \cdots +k_\ell=n$, then  for any $j_1\in\N_0$, 
$$ \bigcup_{j_2=0}^{N-1}\cdots\bigcup_{j_\ell=0}^{N-1} \Gamma_{j_1 \to j_2 \to \cdots \to j_{1}}^{k_1+k_2 + \cdots +k_\ell}
 = \big\{ \pi \in \Gamma_{j_1\to j_1}^n :\  \pi(k_1)<N,\ \pi(k_1+k_2)<N,\  \cdots,\ \pi(k_1+\cdots+ k_{\ell-1})<N \big\}   $$
 and 
 $$\sum_{j_2=0}^{N-1}\cdots\sum_{j_\ell=0}^{N-1}  \M\big\{ \Gamma_{j_1 \to j_2 \to \cdots \to j_{\ell+1}}^{k_1+k_2 + \cdots +k_\ell}\big\}  = \sum_{\pi\in \Gamma_{j_1 \to j_1}^n}  \M\{\pi\} \1_{\max\{\pi(k_1),\cdots,\pi(k_1+\cdots+k_{\ell-1})\} < N} .
 $$
Hence, we can rewrite the RHS of formula (\ref{trace_3}) as
 \begin{equation}  \label{trace_4}
  \tr\big[\p_N\M^{k_1}\p_N\cdots\p_N \M^{k_\ell} \p_N\big]  
  = \sum_{j_1=0}^{N-1} \sum_{\pi\in \Gamma_{j_1 \to j_1}^n}  \M\{\pi\} \1_{\max\{\pi(k_1),\cdots,\pi(k_1+\cdots+k_{\ell-1})\} < N}    .
 \end{equation}
According to formula  (\ref{cumulant_0}), we obtain for any $n\in\N$,
\begin{equation} \label{cumulant_1}
\Cu^n_N[F] =  -\sum_{\ell=1}^{n} \frac{(-1)^{\ell}}{\ell} 
\sum_{\begin{subarray}{c}k_1+\cdots+k_\ell = n \\ k_i \ge 1 \end{subarray}} \frac{n!}{k_1!\cdots k_\ell!}
\sum_{j=0}^{N-1} \sum_{\pi\in \Gamma_{j \to j}^{n}}  \M\{\pi\} \1_{\max\{\pi(k_1),\cdots,\pi((k_1+\cdots+k_{\ell-1}))\} < N} .
\end{equation}
Because of the identity
\begin{equation} \label{Mobius_1} \sum_{\ell=1}^{n} \frac{(-1)^\ell}{\ell} 
 \sum_{\begin{subarray}{c}k_1+\cdots+k_\ell = n \\ k_i \ge 1 \end{subarray}} \frac{n!}{k_1!\cdots k_\ell!} = \1_{n=1}  ,
\end{equation}
the contributions of most of the paths in  formula (\ref{cumulant_1}) cancel. 
Formula (\ref{Mobius_1}) follows directly from taking $F=1$ in \eqref{cumulant_0}; see also formula (\ref{Mobius_0}) in the appendix. Combining formulae (\ref{trace_1}) and (\ref{trace_4}), we have shown that for any composition $n=k_1 + \cdots +k_\ell$, 
 \begin{equation*} % \label{trace_4}
  \tr\big[\p_N\M^{k_1}\p_N\cdots\p_N \M^{k_\ell} \p_N\big]  
  =   \tr\big[\p_N\M^{n}\p_N\big]
  - \sum_{j=0}^{N-1} \sum_{\pi\in  \Gamma_{j \to j}^{n}} \M\{\pi\} \1_{\max\{\pi(k_1),\cdots,\pi(k_1+\cdots+k_{\ell-1})\} \ge N}    .
 \end{equation*}
Hence,  by (\ref{Mobius_1}), we obtain for any $n\ge 2$,
\begin{equation*} \label{cumulant_2}
\Cu^n_N[F] =  \sum_{\ell=2}^{n} \frac{(-1)^\ell}{\ell} 
\sum_{\begin{subarray}{c}k_1+\cdots+k_\ell = n \\ k_i \ge 1 \end{subarray}} \frac{n!}{k_1!\cdots k_\ell!}
 \sum_{j=0}^{N-1} \sum_{\pi\in \Gamma_{j \to j}^{n}}  \M\{\pi\} 
\1_{\max\{\pi(k_1),\cdots,\pi(k_1+\cdots+k_{\ell-1})\} \ge N}  .
\end{equation*}

Making the change of variables ${\bf n}=\Psi(\k)$ given by \eqref{CV} -- \eqref{mho} in the previous formula, we conclude that for any $n\ge 2$,
\begin{equation} \label{cumulant_3}
\Cu^n_N[F] =
\sum_{{\bf n}\in \Lambda_n} \mho({\bf n})
 \sum_{j=0}^{N-1} \sum_{\pi\in \Gamma_{j \to j}^{n}}  \M\{\pi\}
\1_{\max\{\pi(n_1),\cdots,\pi(n_{\ell})\} \ge N} .
\end{equation}

This  completes the proof of lemma~\ref{thm:cumulants}. 
Using the notation (\ref{path}), we see that $\Cu^n_N[F] = \varpi^{n}_{N}\big(F(\J)\big)$ and if we let
$$ \L_{ij} =  \J_{N+i ,N+j} \1_{ i,j  \ge - N } $$
then $\varpi^{n}_{N}\big(F(\J)\big)= \varpi^{n}_{0}\big(F(\L)\big)$. 
Moreover, if the entries of the recurrence matrix $\J$ are bounded (independently of $N$)  and satisfies (\ref{band}), we deduce from lemma~\ref{thm:bound} below that there exists a constant $ \mathrm{C}_F>0$ such that  for any $n,N\in \N$,
\begin{equation}\label{cumulant_bound}
 \big| \Cu^n_N[F]  \big| \le  n!  \exp( n \mathrm{C}_F) \ .
 \end{equation}

\begin{lemma} \label{thm:bound}  
Let $\L$ be a doubly-infinite bounded matrix  which satisfies the condition $(\ref{band})$. Then, for any polynomial $F$, there exists a constants $\mathrm{C}_F>0$ such that for any $n\ge 2$,
$$ \big| \varpi^{n}_{0}\big(F(\L)\big) \big| \le  n!  \exp( n \mathrm{C}_F) \ . $$
Moreover, if $ F(x)= \sum_{k=0}^{d} c_k x^k$ and $\mathrm{A}= \max\big\{| c_k| : k\le d\big\}$, then 
$\mathrm{C}_F \le  \log(2d \mathrm{A}) + d\log(\W C)$  where $C= \sup\{ |\L_{ij}| : i, j\in \Z \}$. 
\end{lemma}

\proof According to formula (\ref{path}), letting $\M=F(\L)$, we may rewrite
\begin{equation} \label{path'}
\varpi^{n}_{0}\big(\M\big) = 
 \sum_{{\bf n} \in \Lambda_n} \mho({\bf n})   \sum_{j<0}  \sum_{\pi\in \Gamma_{j \to j}^{n}}  \M\{\pi\}
\1_{\max\{\pi(n_1),\cdots,\pi(n_{\ell})\} \ge 0} \ .
\end{equation} 
First observe that, since the graph $\g(\M)$ has  degree at most $ \W^d$, only the paths which start  at a vertex $ j \ge - \W^{nd}/2$ contribute  to the sum (\ref{cumulant_2}), otherwise their maximum is non-negative. This implies that
\begin{equation*}
\big| \varpi^{n}_{0}\big(\M\big) \big| \le 
 \sum_{{\bf n} \in \Lambda_n} \big| \mho({\bf n}) \big|   \sum_{\begin{subarray}{r} -\frac{\W^{nd}}{2} \le j  < 0 \end{subarray}}\bigg|  \sum_{\pi\in \Gamma_{j \to j}^{n}}  \M\{\pi\}\bigg| \ .
\end{equation*} 
Note that by  formula (\ref{path_0}), the last sum is equal to $\M_{jj}^n$ and, according to formula (\ref{mho}), we have\begin{equation*}
 \sum_{{\bf n} \in \Lambda_n} \big| \mho({\bf n}) \big|  =  \sum_{\ell=1}^{n-1} \frac{1}{\ell+1} \sum_{k_1+\cdots+k_{\ell+1}=n}  \frac{n!}{k_1! \cdots k_{\ell+1} !}  \le n! 2^n \ ,    
 \end{equation*}
since $2^n$ is the number of (integer) compositions of $n\in \N$. Thus, we get
\begin{equation}  \label{bound_2}
\big| \varpi^{n}_{0}\big(\M\big) \big| \le 
  n! 2^n  \sum_{\begin{subarray}{r} -\frac{\W^{nd}}{2} \le j  < 0 \end{subarray}}\big|  \M_{jj}^{n}\big| \ .
\end{equation} 
Since $\M=F(\L)$ with  $F(x)= \sum_{k\le d} c_k x^k$, we have for any $ n \ge 0$
%$$\sum_{\pi \in \Gamma_{k\to k}^n } \M\{\pi\}  =  \M^{n} \{kk\} $$
$$\M^{n}_{jj}=  \sum_{k_1,\dots k_n \le d} \prod_{i=0}^{n-1} c_{k_i} \L^{k_i}_{jj} \  .$$   
If we let $C= \sup\{ |\L_{ij}| : i, j\in \Z \} $ and $\mathrm{A}= \max\big\{ |c_k| : k\le d\big\}$, by applying the triangle inequality,  we obtain
\begin{equation}\label{bound_3}
\left| \M_{j j}^{n}  \right| \le (\mathrm{A} d)^n   C^{nd} \ . 
\end{equation}
Combining the upper-bounds \eqref{bound_2} -- \eqref{bound_3}, we conclude that 
$ \big| \varpi^{n}_{0}\big(\M\big) \big|  \le n! (2\W^d \mathrm{A}dC^d )^n$. \qed\\

\subsection{Gaussian fluctuations: proof of theorem~\ref{thm:CLT}.} \label{sect:fluctuation}

In this section, we will show that  theorem~\ref{thm:CLT}  follows directly from Soshnikov's Main Combinatorial Lemma. 
If $\J^{(N_k)}  \overset{\ell}{\rightarrow}  \L(s)$, then for any $F\in \R\langle x\rangle$,  theorem~\ref{thm:weak_cvg} implies that 
\begin{equation*}
 \Xi_{N_k}(F)-\E{\Xi_{N_k}(F)} \Rightarrow \mathrm{X}(F)  \hspace{1cm}\text{as }k\to\infty .
\end{equation*}
Hence, it remains to prove that the cumulants $ \varpi^{n}_{0}\big(F(\L(s))\big)$ vanish for all $n \ge 3$ and  
 that the variance is given by formula (\ref{variance_s}).
%By formula (\ref{path}), we know that 
%\begin{equation*}
 %\varpi^{n}_{0}(\M) = \sum_{{\bf n}\in \Lambda_n} \mho({\bf n})
 %\sum_{ k= 1}^{\infty} \sum_{\pi\in \Gamma^n_{-k\to-k}}  \M\{\pi\} \1_{\max\{\pi(n_1),\cdots,\pi(n_{\ell-1})\} \ge 0}  \ . 
%\end{equation*}   
According to lemma~\ref{thm:weight} below, if $\L=\L(s)$ is a Laurent matrix, then $\M=F(\L)$ is also a Laurent Matrix with symbol $F(s)$.
In particular,   the weighted graph $\g(\M)$ is translation invariant and, for any $j \in \Z$, 
 there is a bijection $\Gamma_{0\to0}^{n} \mapsto \Gamma_{j\to j }^{n}$ given by $\pi = \tilde\pi +j$ so that  $\M\{\pi\}=\M\{\tilde\pi\}$ and
$$ \max\{\pi(n_1),\cdots,\pi(n_{\ell-1})\} = \max\{\tilde\pi(n_1),\cdots,\tilde\pi(n_{\ell-1})\} +j  . $$
By formula (\ref{path'}), it implies that 
\begin{align} \notag
 \varpi^{n}_{0}(\M) 
&=   \sum_{{\bf n}\in \Lambda_n} \mho({\bf n}) \sum_{\tilde\pi\in \Gamma_{0\to 0}^{n}}  \M\{\tilde\pi\}  \sum_{ j= 1}^{\infty}  \1_{\max\{\tilde\pi(n_1),\cdots,\tilde\pi(n_{\ell-1})\} \ge j} \\
 & \label{cumulant_4}
 =  \sum_{{\bf n}\in \Lambda_n} \mho({\bf n})
 \sum_{\tilde\pi\in \Gamma_{0 \to 0}^{n}}  \M\{\tilde\pi\}\max\{0,\tilde\pi(n_1),\cdots,\tilde\pi(n_{\ell-1})\}  .
\end{align}  
 
We  let $F(s(z))=\sum_{k\in\Z} \beta_k z^k $.  By definition, $(i,j)$ is an edge of the graph $\g(\M)$,  if and only if $\beta_{j-i} \neq 0$. Hence, given any $\omega \in \Z^n$ such that  $\omega_1+\cdots+\omega_n=0$, if we denote 
\begin{equation}\label{path_5}
\pi_\omega(k)= \omega_1+\cdots +\omega_k \hspace{1cm} \forall\ k=0,\dots,n  ,
\end{equation}
then $\pi_\omega \in \Gamma_{0\to 0}^n$ if and only if $\M\{\pi_\omega\}=\beta_{\omega_1}\cdots \beta_{\omega_n} \ne 0$.
Making the change of variables $\tilde\pi=\pi_\omega$ in formula (\ref{cumulant_4}), we obtain
\begin{equation}  \label{cumulant_5}
 \varpi^{n}_{0}(\M) 
 =  \sum_{{\bf n}\in \Lambda_n} \mho({\bf n}) 
\sum_{\begin{subarray}{c}\omega_1+\cdots+\omega_n = 0\\ \omega_i \in\Z \end{subarray}} 
\beta_{\omega_1}\cdots \beta_{\omega_n} \max\{0,\pi_\omega(n_1),\cdots, \pi_\omega(n_{\ell-1})\}   .
\end{equation}  
It remains to observe that by (\ref{path_5}) and (\ref{G}), for any $\omega \in \Z^n$, we have
\begin{equation*} 
 \sum_{{\bf n} \in \Lambda_n} \mho({\bf n}) \max\{0,\pi_\omega(n_1),\cdots, \pi_\omega(n_{\ell-1})\}  
 = \G_n(\omega_1,\dots, \omega_n) .
 \end{equation*}
Hence, by formula (\ref{cumulant_5}), this implies that
\begin{equation}  \label{cumulant_6}
 \varpi^{n}_{0}(\M) 
 = \sum_{\begin{subarray}{c}\omega_1+\cdots+\omega_n = 0\\ \omega_i \in\Z \end{subarray}} 
\beta_{\omega_1}\cdots \beta_{\omega_n}\ \G_n(\omega_1,\dots, \omega_n)    .
\end{equation}  
Since $\M= \L\big(F(s) \big)$ and  $ \beta_\omega= \widehat{F(s)}_\omega$, this is nothing but formula (\ref{L_path}) given in the introduction.  By symmetry, we can rewrite formula (\ref{cumulant_6}),
\begin{equation*}
 \varpi^{n}_{0}(\M) 
 = \sum_{\begin{subarray}{c}\omega_1+\cdots+\omega_n = 0\\ \omega_i \in\Z \end{subarray}} 
\beta_{\omega_1}\cdots \beta_{\omega_n}\ \frac{1}{n!} \sum_{\sigma \in \Sy(n)} \G_n(\omega_{\sigma(1)},\dots, \omega_{\sigma(n)})   .
\end{equation*}  

By the MCL, lemma~\ref{thm:G}, we conclude that $ \varpi^{n}_{0}(\M)=0$ for all $n\ge 3$ and 
\begin{equation*} 
 \varpi^{2}_{0}(\M) = \frac{1}{2} \sum_{\omega \in \Z} \beta_{\omega}\beta_{-\omega} |\omega|  .
\end{equation*}

This shows that, if $\L=\L(s)$ is a Laurent matrix, then the limit $\mathrm{X}(F)$  of the random variables $\Xi_{N_k}(F)$ as $k\to\infty$ is Gaussian and the variance is given by      
\begin{equation}  \label{variance_1}
 \var\big[ \mathrm{X}(F)\big]= \sum_{\omega>0} \omega \beta_{\omega}\beta_{-\omega}  ,
\end{equation}
where 
$$  \beta_\omega= \widehat{F(s)}_\omega =  \frac{1}{2\pi i} \oint F\big(s(z)\big) z^{-\omega} \frac{dz}{z}  .  $$
This completes the proof of the central limit theorem~\ref{thm:CLT}.

 \begin{lemma} \label{thm:weight}
Let $\L$ be a Laurent matrix with symbol $s(z)$. For any $n \in\N$, we have
\begin{equation}  \label{weight_0}
s(z)^n =  \sum_{j\in\Z}  z^{j}  \L\big\{\Gamma_{0\to j}^n\big\}  .
\end{equation}
Moreover, if $F\in\R\langle x\rangle$, then $F(\L)$ is  a Laurent matrix with symbol $F(s)$. 
\end{lemma}

 \proof Recall that we denote $s(z)=\sum_{\omega\in\Z} \widehat{s}_\omega z^\omega$ so that for any $n\in\N$, we have
  $$ s(z)^n =  \sum_{ \omega \in \Z^n  } \prod_{k=1}^n s_{\omega_k}   z^{\omega_k} . $$
 Then, using the notation (\ref{path_5}),  $\prod_{k=1}^n s_{\omega_k} = \L\{\pi_\omega\}$, and we obtain
\begin{equation} \label{weight_1}
 s(z)^n =  \sum_{ \omega \in \Z^n  } \L\{\pi_\omega\}  z^{ \pi_\omega(n)}  . 
 \end{equation}
Since the sum (\ref{weight_1}) is over all path of length $n$ starting at $0$ on the graph $\g(\L)$, this proves 
formula (\ref{weight_0}). Then, since the weighted graph $\g(\L)$ is translation invariant, by (\ref{path_0}), we have for all $i,j \in\Z$,
$$ \L^{n}_{ij} = \L\big\{\Gamma^n_{i\to j}\big\} =  \L\big\{\Gamma^n_{0\to j-i}\big\} .$$

According to formula  (\ref{weight_0}), this implies that for all $i,j \in\Z$,
$$ \L^{n}_{ij} = \widehat{s(z)^n}_{j-i}  .$$
By definition, we conclude that $\L^n$ is a Laurent matrix with symbol $s(z)^n$. The same property hold for any polynomial  by taking linear combinations. \qed\\

\section{Equilibrium measures and Law of Large Numbers} \label{sect:LLN}

We will now explain how to compute the equilibrium measure and its moments  from the asymptotics of the recurrence coefficients associated to a BOE. Surprisingly, this requires  stronger assumption than that of theorem~2.4. In particular,  we recover a result of Ledoux, \cite{Ledoux}, that the equilibrium measures of the classical OPEs can be expressed as a mixture of the arcsine law with an independent uniform distribution. 
We also prove a law of large numbers for the BOEs considered in 
sections~\ref{sect:product_Ginibre} and~\ref{sect:BME}. 

\begin{theorem} \label{thm:LLN}
 Suppose that the functions $(P_k^N)_{k=0}^\infty$ satisfy a recurrence relation,  that the recurrence matrix $\J^{(N)}$ satisfies the condition \eqref{band} and 
that for all $t\in [0,1]$  and for all $i,j\in \mathbb{Z}$,
\begin{equation} \label{equilibrium_3}
\lim_{N\to \infty}    \J_{\lfloor Nt  \rfloor+i, \lfloor Nt  \rfloor+j}^{(N)}  = \M(t)_{i,j} . 
\end{equation}
Moreover, if the entries of the  matrices $ \M(t)$  are all  bounded by $\A>0$ for all $t\in[0,1]$.  Then, for any polynomial $F\in \R\langle X \rangle$,
\begin{equation*}  %\label{LLN_0}
\lim_{N\to\infty} \frac{ \Xi_N(F)}{N} = \int_0^1 F\big(\M(t)\big)_{0,0} dt ,  
\end{equation*}
almost surely. 
\end{theorem}

\proof By formula~\eqref{Laplace_5},
\begin{equation*} 
\E{\Xi_{N}(F)} = \Cu^1_N[f]  = \tr[\p_N F(\J)\p_N] ,
\end{equation*}
so that 
\begin{equation*}
\frac{\E{\Xi_{N}(F)}}{N} = \int_0^1   F(\J)_{\lfloor Nt  \rfloor, \lfloor Nt  \rfloor} dt . 
\end{equation*}
 The condition \eqref{equilibrium_3} implies that for any polynomial $F$ and for  $t\in [0,1]$,
 $$ \lim_{N\to\infty} F(\J)_{\lfloor Nt  \rfloor, \lfloor Nt  \rfloor} =  F\big(\M(t)\big)_{0,0} . $$
By assumption, the functions $t\mapsto  F(\J)_{\lfloor Nt  \rfloor, \lfloor Nt  \rfloor}$ are measurable and  uniformly bounded on $[0,1]$. So, by the dominated convergence theorem, we obtain
\begin{equation} \label{LLN_1}
\lim_{N\to\infty} \frac{\E{\Xi_{N}(F)}}{N} = \int_0^1  F\big(\M(t)\big)_{0,0}  dt .
\end{equation}
 
  The condition \eqref{equilibrium_3} also implies that $\M(1)$ is the right-limit of the recurrence matrix $\J^{(N)}$
and, by theorem~\ref{thm:weak_cvg}, 
$$ \lim_{N\to\infty} \Var \Xi_{N}(F)  = \Var \mathrm{X}(F) <\infty .$$
Hence, by Chebyshev's inequality and the Borel-Cantelli lemma, the sequence   of random variables 
$\frac{\Xi_N(F)-\E{\Xi_{N}(F)}}{N} $ converges almost surely to 0 and the result follows from  formula \eqref{LLN_1}. \qed\\

\begin{corollary} \label{thm:mu}
 Under the assumption of  theorem~\ref{thm:LLN}, there exists a (unique) probability measure $\mu_*$ with compact support on $\R$ so that for any $F\in \R\langle X\rangle$,  
\begin{equation} \label{LLN_0}
\lim_{N\to\infty} \frac{ \Xi_N(F)}{N} = \int_\R  F(x) d\mu_*(x)
\end{equation}
almost surely. The measure $\mu_*$ is called the equilibrium measure associated to the biorthogonal ensemble $\Xi_N$. Moreover, if $\M(t)$  are Laurent matrices with symbols $s_t$ for all $t\in(0,1)$, then for all $n\ge 1$
\begin{equation} \label{moment_1} 
\mathfrak{M}_*^n := \int_\R x^n d\mu_*(x)  = \frac{1}{2\pi i}  \oint  \bigg( \int_0^1 s_t(z)^n dt  \bigg)\frac{dz}{z} .  
\end{equation}
\end{corollary}

\proof

Using the notation of section~\ref{sect:path},  for any $n\ge 1$,
\begin{equation*} 
 \M(t)^n_{0,0} = \sum_{\pi\in\Gamma_{0\to0}^n} \M(t)\{\pi\}
 \end{equation*}
and an elementary estimate shows that  $\big| \M(t)^n_{0,0} \big| \le (2\W\A)^n $ for all $t\in [0,1]$. Thus, the sequence $\big( \M(t)^n_{0,0} \big)_{n=1}^\infty $ satisfies the so-called {\it Carleman's condition}
and there exists a (unique) probability measure $\mu_t$ with compact support inside $[-2\W\A, 2\W\A]$  so that for all $n\ge 1$, 
%by theorem~3.3.12 in \cite{Durrett}, 
\begin{equation} \label{LLN_3} \M(t)^n_{0,0}  = \int_\R x^n d\mu_t(x) .
\end{equation}
%If it is absolutely continuous, $d\mu_t(x) = \rho_t(x) dx$, 
Then, by theorem~\ref{thm:LLN} and Fubini's theorem, we see that for all $n\ge 1$, 
$$
\lim_{N\to\infty} \frac{ \Xi_N(X^n)}{N} = \int_\R x^n d\mu_* 
$$
where
\begin{equation} \label{equilibrium_4}
\mu_*(x) = \int_0^1 \mu_t(x) dt . 
\end{equation}
Moreover, if $\M(t)$ is a Laurent matrices with symbols $s_t$,
then by formulae \eqref{weight_0}  and \eqref{LLN_3},  we obtain
\begin{equation} \label{LLN_4}
 \int_\R x^n d\mu_t(x)  = \frac{1}{2\pi i} \oint s_t(z)^n \frac{dz}{z} 
\end{equation}
for all $n\ge 1$. Hence, if we integrate this formula over $t\in (0,1)$, we get the expression \eqref{moment_1}. \qed\\  

\begin{remark} To extend the law large numbers from polynomials to all continuous test functions, it suffices to have a good tail bound for the one-point function $u_N(x) = N^{-1} K_N(x,x)$. Namely, if there exists $L,C>0$ so that  for all $n\ge 1$
$$ \int_{\R\backslash[-L,L]} x^{2n} u_N(x) dx  \le C ,  $$
then for any function $f\in C(\R)$ with at most polynomial growth at $\infty$,
 $$ \lim_{N\to\infty} \frac{ \Xi_N(f)}{N} = \int_\R  f(x) d\mu_*(x)$$ 
 in probability. 
\end{remark}

Note that if the conditions~\eqref{equilibrium_1} are satisfied,  then $\M(t)$ is a tridiagonal  Laurent matrix with symbol $s_t(z)=a(t) z^{-1} +b(t) + a(t) z$ and its spectral measure at $0$ is given by
\begin{equation} \label{arcsine}
 d\mu_t(x) =   \frac{dt}{\pi \sqrt{(\alpha_+(t) - x)(x-\alpha_-(t))}}  \1_{\alpha_-(t)<x<\alpha_+(t)}  \hspace{.5cm}\text{where }\
 \alpha_\pm(t)= b(t)\pm 2a(t)  \ . 
 \end{equation}
Hence, using \eqref{equilibrium_4}, we recover the formula \eqref{equilibrium_0} for the equilibrium measure obtained by Kuijlaars and Van Assche, \cite{KV_99}.\\

When the symbol $s_t$ is more complicated,  it is difficult to get an explicit formula for the spectral measure $\mu_t$ and the representation \eqref{equilibrium_4} is not so helpful. However, we can still use formula \eqref{moment_1} and residue calculus  to compute the moments of the equilibrium measure. 
For instance, for the GUE (with the normalization $V(x)=2x^2)$, the recurrence coefficients of the normalized Hermite polynomials are 
$$a_n^N =\frac{1}{2} \sqrt{\frac{n+1}{N}} 
\hspace{.5cm}\text{and}\hspace{.5cm}
b_n^N =0 ,  $$
so that $s_t(z) = \sqrt{t}\frac{z+1/z}{2}$. Thus, we see that $\mathfrak{M}_*^{2n+1} = 0$ and 
$$
\mathfrak{M}_*^{2n} = \frac{1}{2n+1} \frac{1}{2\pi i}  \oint  \bigg( \frac{z+1/z}{2} \bigg)^{2n}\frac{dz}{z}
=  4^{-n} \frac{1}{n+1}  {2n \choose n} $$
for all $n\ge 0$. We recognize the Catalan numbers which proves that the equilibrium measure of the GUE is the semicircle law: $\rho_*^{\text{GUE}}(x) = \frac{2}{\pi} \sqrt{1-x^2} \1_{|x|<1}$. 
We can apply the same method to compute the equilibrium measure for the square singular values of products of rectangular complex Ginibre matrices.

\begin{theorem} \label{thm:Fuss_Catalan}
Using the notations of section~\ref{sect:product_Ginibre}, the point process \eqref{M_Ginibre} satisfies a law of large numbers  of the type \eqref{LLN_0} and the equilibrium measure depends continuously on the parameters $\gamma_1,\dots, \gamma_m$. Moreover, in the case $\gamma_1=\cdots =\gamma_m=1$, its moments are given by
\begin{equation*}
 \int_\R x^n d\mu_{m,*}(x)  = \mathfrak{C}_{m}^n := \frac{1}{n(m+1)+1} {n(m+1)+1 \choose n} . 
 \end{equation*}
The moments $(\mathfrak{C}_{m}^n)_{n=0}^\infty$ are called the Fuss-Catalan numbers  and the equilibrium measure $\mu_{m,*}$ is called the Fuss-Catalan distribution with variance $m$, \cite{FW_15}. 
\end{theorem}

\proof
Proposition~\ref{thm:KZ_3} shows that for any $t>0$,   the sequence $\big(\J_{\lfloor Nt  \rfloor+i, \lfloor Nt  \rfloor+j}\big)_{i, j}$ converges to a  Laurent matrix with symbol 
\begin{equation}\label{FC}
s_t(z) = \sum_{k=-1}^m z^{-k}  \sum_{\begin{subarray}{c} S \subseteq \{0,\dots, m\} \\ |S| = k+1 \end{subarray}}   
 \prod_{l \in S} \gamma_l \prod_{i \notin S} \big(1 - (1-t)\gamma_i \big) = z \prod_{k=0}^m
\big(  \gamma_k/z+ 1 - (1-t)\gamma_k\big) .
\end{equation}
as $N\to\infty$ and $N_l/N \to 1/ \gamma_l$ for all $l=1,\dots, m$.  
Hence, the assumptions of theorem~\ref{thm:LLN} are satisfied and we obtain a law of large numbers for polynomial linear statistics of the rescaled square singular values process. In particular, in the regime where~$\gamma_1=\cdots =\gamma_m=1$,  we have
$$
s_t(z) = z^{-m} \big(1+ zt\big)^{m+1} . 
$$

Then, for any $z\in\C\backslash\{0\}$, 
\begin{align*} \notag
\int_0^1 s_t(z)^n dt  &= z^{-mn} \int_0^1  \big(1+ zt\big)^{n(m+1)} dt \\
&=  \frac{(1+z)^{n(m+1)+1}-1}{(n(m+1)+1) z^{mn+1} }  .
\end{align*}
According to formula \eqref{moment_1}, this implies that for all $n\ge 1$, 
\begin{equation*}
\mathfrak{M}_*^n = \frac{1}{n(m+1)+1} {n(m+1)+1 \choose n} . 
\end{equation*}
 In particular, if $m=1$, we obtain the usual Catalan numbers. \qed\\

%%% Work

\begin{remark} By formula~\eqref{Laguerre_6}, it also follows from this proof that, once properly rescaled, the Laguerre Muttalib-Borodin ensemble with parameter $\theta=1/m$ also satisfies a  law of large numbers and its equilibrium measure is also  the Fuss-Catalan distribution with variance $m$. On the other hand, 
for the Jacobi Muttalib-Borodin ensemble, proposition~\ref{thm:Jacobi_asymptotics} shows that the recurrence matrix has a limit (without scaling) which is a Laurent matrix with symbol
$$ s_t(z) = z^{-1} \left( \frac{z+m}{m+1}\right)^{m+1}  $$ 
for any $t>0$. According to formula \eqref{moment_1}, this implies that the moments of the equilibrium measure are given by for all $n\ge 1$, 
\begin{equation} \label{moment_3} 
\mathfrak{M}_*^n =
\frac{1}{2\pi i}  \oint z^{-n} \left( \frac{z+m}{m+1}\right)^{n(m+1)}    \frac{dz}{z} 
= \beta_m^{-n}  { n(m+1) \choose n } 
\end{equation}
where $\beta_m = m(1+ 1/m)^{m+1}$. In particular, this gives an elementary proof of proposition~3.10 in~\cite{FW_15}.
\end{remark}

It turns out that the support of the Fuss-Catalan distribution $\mu_{m,*}$ is the interval $[0,\beta_m]$ where $\beta_m = m(1+ 1/m)^{m+1}$. The next proposition provides an estimate on the tail of the largest square singular value.
This bound is far from optimal but it follows rather directly from an  estimate on the growth of the recurrence coefficients. In particular, combining theorem~\ref{thm:Fuss_Catalan} and  proposition~\ref{thm:tail_bound}, this implies that
$$
\frac{\mu_{\max}}{ N^m} \to \beta_m
$$
almost surely as $N\to \infty$. 

\begin{proposition} \label{thm:tail_bound}
Let $\mu_{\max}$ be the largest square singular value of the product of $m$ independent $N\times N$ Ginibre matrices (i.e. $\eta_1= \cdots =\eta_m=0$ and $M=N^m$).  
Then there exists a constant  $C>0$ which depends only on the parameter $m$ such that for any $s>0$,
$$
\P{\mu_{\max} \ge N^m \beta_m e^s } \le C N^{3/4} e^{- \sqrt{N} s }  . 
$$
\end{proposition}

\proof 
By Markov's inequality, for any $s>0$ and $q \in \N$, 
\begin{align}
\P{\mu_{\max} \ge N^m \beta_m e^{s} } 
& \notag \le \E{\left(\frac{\mu_{\max}}{N^m}\right)^q } \beta_m^{-q}e^{-q s } \\
& \label{Markov}\le \E{\Xi_N(x^q)}
\beta_m^{-q}e^{-q s}
\end{align}
where $\Xi_N$ is the rescaled square singular-value process \eqref{M_Ginibre}.
By proposition~\ref{thm:KZ_1}, the recurrence coefficients associated to the BOE $\Xi_N$ satisfy
\begin{equation*} 
 \alpha_{n,-k}  = \1_{-1 \le k \le m}  \frac{1}{(k+1)!}  \sum_{i=0}^{k+1} (-1)^{i} {k+1 \choose i}   (n-i+1)^{m+1}  .
\end{equation*}
In particular, $ \alpha_{n,1}= (n+1)^{m+1}$, $ \alpha_{n,m}=1$,  and using a variant of lemma~\ref{Binomial_identity}, we have that for any $0 \le k <m$, 
\begin{equation*} 
 \alpha_{n,-k}  =  (n+1)^{m-k} {m+1 \choose k+1} 
 - (n+1)^{m-k-1 }  {m+1 \choose k+2}\frac{k+2}{2}   + \underset{n\to\infty}{O}(n^{m-k-2 }) .
\end{equation*}
By formula \eqref{KZ_2} with $\kappa_n(N)=N^n$, this implies that there exist a constant $C>0$ such that for all $n \ge 0$
\begin{equation} \label{upper_bound_1}
0\le \J_{n,n-k}  \le  \1_{-1 \le k \le m}  \left(\frac{C \vee(n+1)}{N}\right)^{m-k} {m+1 \choose k+1} . 
\end{equation}
Let $\L_t$ be the Laurent matrix with symbol $ s_t(z)= \sum_{k=-1}^m {m+1 \choose k+1} z^k t^{m-k} = z^{-1} \big(z+ t\big)^{m+1}$. The estimate \eqref{upper_bound_1} implies that for any $q \in \N $ and $0\le j <N$, we have
\begin{equation} \label{upper_bound_2}
\J(\Gamma_{j\to j}^q) \le \L_{1+q/N}(\Gamma_{0\to 0}^q). 
\end{equation}
Moreover, by formula \eqref{weight_0}, it is easy to check that for any $q\in N$,
 $$
 \L_t(\Gamma_{0\to 0}^q) =  \frac{1}{2\pi i} \oint s_t(z)^n \frac{dz}{z} 
=   { q(m+1) \choose q}  t^{qm} .
$$
Thus, by Stirling's formula, if $q$ is sufficiently large and $t= 1+q/N$, we get
 $$
 \L_t(\Gamma_{0\to 0}^n) \le q^{-1/2}  \beta_m^q  e^{m q^2/N} . 
 $$
 Combined with \eqref{upper_bound_2}, this estimate implies that
 $$
 \E{\Xi_N(x^q)} = \sum_{j=0}^{N-1}\J(\Gamma_{j\to j}^q) \le  q^{-1/2}   N  \beta_m^q  e^{m q^2/N}  . 
 $$
 Hence, by \eqref{Markov},  we conclude that
 \begin{equation*}
\P{\mu_{\max} \ge  N^m \beta_m e^{s} } 
\le  q^{-1/2}  N e^{mq^2/N-q s }
\end{equation*}
 and it remains to choose $q = \lceil \sqrt{N} \rceil$. \qed\\

Finally, If $m=1$, the distribution of the square singular values of a rectangular Ginibre matrice is known as the Wishart or Laguerre (unitary) ensemble. Moreover, its equilibrium measure is the celebrated Marchenko-Pastur distribution, \cite{MP67},  and it depends only on the parameter 
$$\gamma := \lim_{N\to\infty} N_1/N \in [0,1] .$$

We conclude  this section by recovering the Marchenko-Pastur law using the  formulae \eqref{equilibrium_4} and \eqref{FC}. In this case, $\gamma_0=1$ and $\gamma_1=\gamma$ so that the symbol of the matrix $\M(t)$ is given by   
\begin{equation} \label{MP_2}
 s_{t}(z) = z t \big(1-(1-t) \gamma\big) + 1-\gamma +2t\gamma + \gamma/z  
 \end{equation}
for any $t\in (0,1]$. Since $ s_{t}(z)$ is a Laurent polynomial of degree 1, the spectral measure of the matrix $\M(t)$ is an arcsine law, \eqref{arcsine},  with parameters 
\begin{equation} \label{MP_1}
a(t) = \sqrt{\gamma t (1-\gamma +t \gamma) }
\hspace{.5cm}\text{and}\hspace{.5cm}
 b(t) = 1-\gamma +2t\gamma .  
 \end{equation}
  
An elementary computation shows that if $ |x-b(t)| < 2 a(t)$, then 
  $$
 \frac{d}{dt} \sqrt{b(t)+2a(t)-x} \sqrt{ x-b(t)+2a(t)} =  \frac{(x-b(t))b'(t)+ 4a(t)a'(t)   }{\sqrt{b(t)+2a(t)-x} \sqrt{ x-b(t)+2a(t)} } . 
  $$
Then, by \eqref{MP_1}, we can check that  
  $$4a(t)a'(t) = 2\gamma \big( 2\gamma t + 1-\gamma \big) = b(t)b'(t)$$
  so that,
    $$
 \frac{d}{dt} \frac{\sqrt{\alpha_+(t)-x} \sqrt{ x-\alpha_-(t)}}{2\gamma x} =  \frac{1  }{\sqrt{\alpha_+(t)-x} \sqrt{ x-\alpha_-(t)} } 
  $$
  where  $\alpha_{\pm}(t) = b(t)\pm 2a(t)$. 
Since the functions $\alpha_+$ and $\alpha_-$ are increasing, respectively decreasing, on the interval $[0,1]$,   by formula \eqref{equilibrium_4}, we obtain
\begin{align*}
  \frac{ d\mu_{\gamma}^{\text{MP}} }{dx}
  &=   \int_0^1\frac{\1_{\alpha_-(t)<x<\alpha_+(t)}  }{\pi\sqrt{(\alpha_+(t)-x)( x-\alpha_-(t))} } dt  \\
&=      \begin{cases}\displaystyle \frac{\sqrt{(\alpha_+(1)-x)( x-\alpha_-(t))}}{2\pi\gamma x}  &\text{if}\ \alpha_-(1)<x<\alpha_+(1) \\
\ 0& \text{else}
\end{cases}
\end{align*}
By definition, $\alpha_{\pm}(1) =1+\gamma \pm 2\sqrt{\gamma} = \left(1\pm \sqrt{\gamma}\right)^2$
 and we recover the well-known formula for the density of the Marchenko-Pastur law:
\begin{equation*}
  \frac{ d\mu_{\gamma}^{\text{MP}} }{dx}=  \frac{\sqrt{(\beta_+ - x)(x-\beta_-)}}{2\pi \gamma x}  \1_{\beta_-<x<\beta_+}  \hspace{.4cm}\text{where }\
 \beta_\pm= \left(1\pm \sqrt{\gamma}\right)^2   . 
 \end{equation*}\\

%Using formulae \eqref{moment_1} and \eqref{MP_2}, it is also possible to compute the moments of the  Marchenko-Pastur law:
%$$
%\int_\R x^n d\mu_\gamma^{\text{MP}}(x) = \sum_{k=0}^{n-1} \frac{\gamma^k}{n} {n \choose k} {n \choose k+1} . 
%$$
%\begin{remark}

{\large\bf Acknowledgment:}
I thank Maurice Duits, Adrien Hardy and Kurt Johansson for many valuable discussions on their related works. In particular, I am grateful Kurt Johansson for reading carefully the draft of this article and to Maurice Duits for suggesting the application to products of independent complex Ginibre matrices which resulted in  section~\ref{sect:product_Ginibre}.\\

\appendix

\section{Appendix: Proof of the Main Combinatorial Lemma} \label{A:MCL}

In this section, we give a proof of the MCL, lemma~\ref{thm:G}.  It is different from the argument of~\cite{Soshnikov_00a} and includes additional details. Moreover, a factor $1/2$ seems to be missing in the formulation of the MCL in~\cite{Soshnikov_00a}.
Based on the seminal results of  Spitzer, \cite{Spitzer_56}, on the geometry of simple random walks, we give an alternative proof of  a formula due to Rudnick and Sarnak, cf.~lemma~\ref{thm:RS}, which plays a key role in Soshnikov's proofs of both the MCL and Spohn's lemma, \cite{Spohn_87}. In fact, this shows that all these  lemmas follow from the DHK formulae, (\ref{Kac}), and  that the combinatorial machinery behind Kac and Soshnikov proofs of the Strong Szeg\H{o} theorem is the same. This fact was certainly known to experts. However it  seems to be missing from the literature. \\

In section~\ref{sect:Szego}, we have seen that, for the CUE, the fluctuations of linear statistics are described by the    
the Strong Szeg\H{o} limit theorem.  In~\cite{Kac_54}, Kac gave a combinatorial proof of the Strong Szeg\H{o} theorem and a continuum analog for Fredholm determinant,  also know as the {\it Ahiezer-Kac formula}.  His  approach is also the basis of further generalizations to pseudo-differential operators on manifolds; see \cite{Gioev_01} or \cite[section~6.5]{Simon_04a} for further references. 
Kac showed that, for sufficiently smooth functions, the Strong Szeg\H{o} theorem boils down to the following identity
 \begin{equation} \label{Kac}
 \sum_{\sigma\in \Sy(n)}\m(0, x_{\sigma(1)},\dots, x_{\sigma(n)}) =  \sum_{\sigma\in \Sy(n)} x_{\sigma(1)}\sum_{k=1}^n\theta(x_{\sigma(1)}+\dots+ x_{\sigma(k)})  \ ,
\end{equation}  
where $\theta$ denotes the Heaviside step function and for any $x\in\R^n$,
\begin{equation}\label{M}
\m(x_1,\dots, x_n)= \max\big\{ x_1,x_1+x_2,\dots, x_1+\cdots+x_n\big\} \ .
\end{equation}

In~\cite{Kac_54},  Kac gave an elementary proof of formula (\ref{Kac}) which is due to Dyson and Hunt. Therefore it is common to call  (\ref{Kac}) the Dyson-Hunt-Kac (DHK) formula, see~\cite[Thm.~6.5.3]{Simon_04a}. He also provided an interpretation of (\ref{Kac}) in terms of the expected value of the maximum of a $n$-step simple random walks, cf.~formula (\ref{Spitzer'}) below.  The general mechanism to keep track of the distribution of the  maximum of a random walk with independent identically distributed increments was understood in~\cite{Spitzer_56}. Namely,  let $X_1,X_2,\dots$ be i.i.d.~real-valued random variables,  let $(S_k)_{k\ge 0}$ be their partial sums, and for any $u\in\R$, let  $u^+=\max\{0,u\}$. Then,  for any $t \in\R$ and $|\lambda|<1$,  
\begin{equation} \label{Spitzer}
\sum_{n=0}^\infty \E{e^{i t  \max\{S_0,\dots, S_n\} }}  \lambda^n =\exp\bigg(\sum_{n=1}^\infty \E{e^{i t S_n^+}}  \frac{\lambda^n}{n} \bigg) \ .
\end{equation} 
The proof of formula (\ref{Spitzer}) is based on a beautiful bijection  due to  Bohnenblust, see \cite{Spitzer_56},  valid for any $x\in\R^n$, between the set $\{\m(0,x_{\sigma(1)},\dots, x_{\sigma(n)})\}_{\sigma \in \Sy(n)}$ and $\{T(x ; \sigma)\}_{\sigma \in \Sy(n)}$ where 
$$T(x ; \sigma)= \sum_{\text{cycles } \tau \text{ of }\sigma } \bigg(\sum_{j\in \tau} x_j \bigg)^+ \ .$$
Formula (\ref{Spitzer}) has several applications in the theory of  random walks. For instance, a formula for the joint distribution of $ \big(\max\{S_0,\dots, S_n\},S_n\big)$ and a nice proof of the strong law of large numbers are given in~\cite{Spitzer_56}; see also~\cite{Steele_02} for a modern reference.
It also leads to formula (\ref{Kac}).
%It is not difficult to show that Spitzer's bijection implies formula (\ref{RS}). In fact,
 If we differentiate formula (\ref{Spitzer}) with respect to the parameter $t$ and evaluate at $t=0$, we obtain
\begin{equation*}
\sum_{n=0}^\infty \E{\m(0,X_1,\dots, X_n) }  \lambda^n =  \frac{1}{1-\lambda} \sum_{n=1}^\infty \E{S_n^+} \frac{ \lambda^n}{n} \ ,
\end{equation*}  
since $\exp\big(\sum_{n=0}^\infty \frac{\lambda^n}{n} \big) = \frac{1}{1-\lambda}$ and $\max\{S_0,\dots, S_n\} =\m(0,X_1,\dots, X_n)$.  If we identify the coefficients of $\lambda^n$ on both sides, we conclude that
\begin{equation} \label{Spitzer'}
 \E{\m(0,X_1,\dots, X_n) } =\sum_{k=1}^n  \frac{\E{(X_1+\cdots+ X_k)^+}}{k} \ .
\end{equation}

Now, let $x\in \R^n$, $p_1,\dots, p_n >0$ be distinct numbers such that $\sum_{j=1}^n p_j=1$ and suppose that $\P{X_1=x_j}=p_j$ for all $j=1, \dots, n$. If we identify the coefficient of $p_1\cdots p_n$ in formula (\ref{Spitzer'}), we obtain 

\begin{equation} \label{DHK}
 \sum_{\sigma\in \Sy(n)}\m(0, x_{\sigma(1)},\dots, x_{\sigma(n)}) =\sum_{k=1}^n  \frac{1}{k}  \sum_{\sigma\in \Sy(n)}\big(x_{\sigma(1)}+\dots+ x_{\sigma(k)}\big)^+  \ .
\end{equation}  

Note that, if $\theta$ denotes the Heaviside step function, then for any $k=1,\dots,n$,  
  \begin{equation*}
 \sum_{\sigma\in \Sy(n)}\big(x_{\sigma(1)}+\dots+ x_{\sigma(k)}\big)^+  = k  \sum_{\sigma\in \Sy(n)}  x_{\sigma(1)}\ \theta(x_{\sigma(1)}+\dots+ x_{\sigma(k)}) \ ,    
\end{equation*} 
and formulae (\ref{Kac}) and (\ref{DHK}) are equivalent.   \\

The connection with random matrix theory was realized in~\cite{RS_96} where the Bohnenblust-Spitzer method was used to derive the following identity.
  \begin{lemma} \label{thm:RS}
 If $x_1+\cdots+x_n= 0$, then 
\begin{equation} \label{RS}
\sum_{\sigma \in \Sy(n)}  \m(x_{\sigma(1)},\dots, x_{\sigma(n)})
 = \frac{n}{4}  \sum_{ F \subset [n] }  (|F|-1)! (n-|F|-1)!\  \bigg| \sum_{k\in F} x_k \bigg|^+ \ ,
\end{equation}
where the sum is over all subsets of  $[n]= \{1,\dots,n\}$ and, by convention $(-1)! =0$. 
 \end{lemma}

%Below, we will show that formulae (\ref{Kac}) and  (\ref{RS}) are equivalent  on the hyperplane $x_1+\cdots+x_n= 0$.

\proof 
First observe that for any $x\in \R^n$ such that $x_1+\cdots+x_n=0$, we have
\begin{equation} \label{M_symmetry} 
 \sum_{\sigma \in \Sy(n)}  \m(-x_{\sigma(1)},\dots, -x_{\sigma(n)}) 
 =\sum_{\sigma \in \Sy(n)}  \m(x_{\sigma(1)},\dots, x_{\sigma(n)})  \ .
\end{equation}
This follows from symmetry and the simple observation that, since $-x_1-\cdots-x_k= x_{k+1}+\cdots+x_n$, by formula (\ref{M}), one has
\begin{equation*}
\m(-x_1,\dots, -x_n)= \max\big\{ x_2+\cdots+x_n, x_3+\cdots+x_n,\dots, x_n,0\big\} \ .
\end{equation*}
Formulae (\ref{M_symmetry}) and (\ref{DHK}) imply that, if $x_1+\cdots+x_n=0$, then
\begin{equation} \label{DHK_1}
 \sum_{\sigma\in \Sy(n)}\m( x_{\sigma(1)},\dots, x_{\sigma(n)}) =\sum_{k=1}^{n-1}  \frac{1}{2k}  \sum_{\sigma\in \Sy(n)}\big|x_{\sigma(1)}+\dots+ x_{\sigma(k)}\big|  \ .
\end{equation}  
Moreover, for any $x\in\R^n$, we have
 \begin{align*}
 \sum_{k=1}^n  \frac{1}{k}  \sum_{\sigma\in \Sy(n)}\big|x_{\sigma(1)}+\dots+ x_{\sigma(k)}\big|  
&= \sum_{k=1}^n \frac{1}{k} \sum_{\begin{subarray}{c} F \subset[n] : |F|=k \end{subarray}}
 \#\bigg\{ \sigma \in \Sy(n) : \{\sigma(1),\dots,\sigma(k)\}=F\bigg\}     \bigg| \sum_{k\in F} x_k \bigg| \\
& =  \sum_{\begin{subarray}{c} F \subset[n]  \end{subarray}}  (|F|-1)! |F^c|!\  \bigg|\sum_{k\in F} x_k \bigg| \ ,
\end{align*}
where $F^c$ denotes the complement of $F$ and $|F|$ is the cardinal of $F$. By convention, the last sum equals 0 if $F=\emptyset$. Finally, if $x_1+\cdots+x_n= 0$,  $|\sum_{k\in F} x_k| =|\sum_{k\in F^c} x_k|$ and we obtain 
 \begin{align} \notag
 \sum_{k=1}^n  \frac{1}{k}  \sum_{\sigma\in \Sy(n)}\big|x_{\sigma(1)}+\dots+ x_{\sigma(k)}\big|  
&=  \frac{1}{2} \sum_{\begin{subarray}{c} F \subset[n] \end{subarray}} \big\{  (|F|-1)! |F^c|!+  (|F^c|-1)! |F|!\big\}     \bigg| \sum_{k\in F} x_k \bigg| \\
& \label{RS'}
=  \frac{n}{2}  \sum_{\begin{subarray}{c} F \subset[n]  \end{subarray}}  (|F|-1)! (|F^c|-1)!\  \bigg|\sum_{k\in F} x_k \bigg|\ .
\end{align}
Hence, combining formulae (\ref{DHK_1}) and (\ref{RS'}), we have proved formula (\ref{RS}). \qed \\

In~\cite{RS_96},  Rudnick and Sarnak investigated the statistical properties of the non-trivial zeroes of primitive $L$-function viewed as a pseudo-random point process. They used lemma~\ref{thm:RS} to prove that in the correct scaling, the correlation functions of the spacings between the zeros agree with the ones of the sine process.
On account of this analogy between the statistical behavior of the zeroes of primitive $L$-functions and the eigenvalues of random matrices from the corresponding compact groups, it is not so surprising that the proof of the MCL in~\cite{Soshnikov_00a} is based on  lemma~\ref{thm:RS}.  \\

{\it Proof of lemma~\ref{thm:G}.} We fix $x\in\R^n$.  We want to prove that, if $x_1+\cdots +x_n=0$, then
 \begin{equation} \label{G''}
 \sum_{\sigma\in\Sy_n} \G_n(x_{\sigma(1)},\dots, x_{\sigma(n)}) 
 = \begin{cases} |x_1| &\text{if}\ n=2 \\ 0 &\text{else} \end{cases} \ ,
 \end{equation} 
 where, according to the definitions (\ref{G}) and (\ref{mho}), we have
\begin{equation} \label{G'}
 \G_n(x_1,\dots, x_n) = \sum_{\ell=1}^n \frac{(-1)^\ell}{\ell} \sum_{k_1+\cdots + k_\ell =n} {n \choose {\bf k}} 
  \max\left\{0,\sum_{i=1}^{k_1} x_{i},\sum_{i=1}^{k_1+k_2} x_{i},\cdots, \sum_{i=1}^{k_1+\cdots+k_{\ell}} x_{i}  \right\} \ .
 \end{equation}

For any finite set $A$, we let $\Pi(A)$ be the set of all partitions of $A$. It means that $\pi=\{\pi_1,\dots, \pi_\ell\}\in \Pi(A)$ if and only if $\emptyset \neq \pi_j \subseteq A$ for all $j \le \ell$ and 
$ \pi_1\uplus \cdots \uplus \pi_\ell= A $. In the following, we will also denote $\ell=|\pi|$.  
We refer to \cite[Example~3.10.4]{Stanley_12}  for an account on the partition lattice $\Pi(A)$. In particular, its M\"{o}bius function is given by 
\begin{equation} \label{Mobius}
 \mu(\pi)= \mu(\hat{0}, \pi)= (-1)^{|\pi|-1} (|\pi|-1)!  
 \end{equation} 
 where $\hat{0}=\{A\}$ is the trivial partition and 
\begin{equation} \label{Mobius_0}
 \sum_{\pi\in \Pi(A)} \mu(\pi) = \1_{|A|=1} \ . 
 \end{equation}
Indeed, if the set $A$ has a single element, then there is only 1 partition and the sum (\ref{Mobius_0}) is equal to $\mu(\hat{0})=1$. Otherwise, the sum vanishes by definition of the M\"{o}bius function. In fact, it is a simple exercise to check that formulae (\ref{Mobius_0}) and  (\ref{Mobius_1}) are the same. 
For any $n\in\N$, we let $\Pi[n]$ be the partition lattice of $[n]=\{1,\dots, n\}$ and for any $\pi=\{\pi_1,\dots, \pi_\ell\}\in \Pi[n]$, we define
 $$ \Upsilon_\pi(x_1,\dots, x_n) = \sum_{\tau\in \Sy(\ell)} \max\bigg\{ \sum_{i\in \pi_{\tau(1)}} x_i +\cdots+\sum_{i\in \pi_{\tau(j)} } x_{i}\  \big|  j=1,\dots \ell   \bigg\}  \ .  $$
Lemma~\ref{thm:RS} implies that if $x_1+\cdots +x_n=0$,
\begin{equation} \label{Soshnikov_2}
 \Upsilon_\pi(x_1,\dots, x_n) = \frac{|\pi|}{4}  \sum_{\begin{subarray}{c} S \subset \{1,\dots,|\pi|\} \end{subarray}}  (|S|-1)! (|S^c|-1)!\  \bigg| \sum_{ i \in \underset{j\in S}{\bigcup}\pi_j} x_i \bigg| \ .
 \end{equation}

On the other hand, for any $x\in\R^n$, given a composition ${\bf k}= (k_1,\dots, k_\ell)$ of $n$, we have
$$ 
\sum_{\sigma\in\Sy(n)}  \max\left\{\sum_{i=1}^{k_1} x_{\sigma(i)},\cdots, \sum_{i=1}^{k_1+\cdots+ k_{\ell}} x_{\sigma(i)}  \right\} 
=  k_1!\cdots k_\ell! \sum_{\begin{subarray}{c}  F_1 \uplus \cdots  \uplus F_\ell = \{1,\dots, n\}   \\ |F_j|=k_j \end{subarray}}  \max\left\{\sum_{i\in F_1} x_{i}\ ,\cdots, \sum_{i\in F_1\cup \cdots\cup F_{\ell}} x_{i}  \right\}  \ .
$$
This formula comes from the fact that for any disjoint subsets $F_1, \dots, F_\ell \subset \{1,\dots,n\}$ such that $|F_j|=k_j $,
$$ \#\left\{ \sigma\in\Sy(n)\ \bigg| \begin{array}{c} \sigma(1),\ \dots,\ \sigma(k_1) \in F_1 \\  \sigma(k_1+1),\ \dots,\ \sigma(k_1+k_2) \in F_2 \\ \vdots \\  \sigma(k_{1}+\cdots+k_{\ell-1}+1),\dots, \sigma(k_1+\cdots+k_\ell) \in F_\ell \end{array}  \right\} = k_1!\cdots k_{\ell}! \ . $$

This implies that  for any $\ell = 1,\dots, n$,
\begin{align*} \sum_{k_1+\cdots + k_\ell =n} {n \choose {\bf k}} 
\sum_{\sigma\in\Sy(n)}  \max\left\{\sum_{i=1}^{k_1} x_{\sigma(i)},\cdots, \sum_{i=1}^{k_1+\cdots+ k_{\ell}} x_{\sigma(i)}  \right\} 
&= n! \sum_{\begin{subarray}{c}  F_1 \uplus \cdots  \uplus F_\ell = \{1,\dots, n\}  \end{subarray}}  \max\left\{\sum_{i\in F_1} x_{i}\ ,\cdots, \sum_{i\in F_1\cup \cdots\cup F_{\ell}} x_{i}  \right\}  \\
&= n! \sum_{\begin{subarray}{c}  F \in \Pi[n] \\ |F|=\ell \end{subarray}} \Upsilon_F(x_1,\cdots, x_n) \ .
\end{align*}

According to formula (\ref{G'}), we have proved that 
\begin{align*}
 \sum_{\sigma\in\Sy(n)} \G_n(x_{\sigma(1)},\dots, x_{\sigma(n)}) = n! \sum_{F \in \Pi[n]} \frac{(-1)^{|F|}}{|F|}  \Upsilon_F(x_1,\cdots, x_n) \ .    
\end{align*}

If $x_1+\cdots +x_n=0$, by formula (\ref{Soshnikov_2}), this implies that
\begin{equation} \label{Soshnikov_1}
 \sum_{\sigma\in\Sy(n)} \G_n(x_{\sigma(1)},\dots, x_{\sigma(n)}) = \frac{n!}{4}  \sum_{F \in \Pi[n]} (-1)^{|F|}    \sum_{\begin{subarray}{c} S \subset \{1,\dots,|F|\} \end{subarray}}  (|S|-1)! (|S^c|-1)!\  \bigg| \sum_{ i \in \underset{j\in S}{\bigcup} F_j} x_i \bigg| \ .
\end{equation}
In the previous sum, we make the change of variables $(F,S)\in \Pi[n] \times 2^{|F|}$ to $(A,\pi,\tilde{\pi}) \in  2^{[n]} \times \Pi(A)\times \Pi(A^c)$ given by
$$ A= \bigcup_{j\in S} F_j \ ,   \hspace{.7cm} 
\pi = \{ F_j  \ |\ j\in S   \}  \ ,  \hspace{.7cm} 
\tilde{\pi} = \{ F_j  \ |\ j\notin S   \} \ .
$$
By (\ref{Mobius}), we have  $(-1)^{|F|} (|S|-1)! (|S^c|-1)! = \mu(\pi)\mu(\tilde\pi)$ and we obtain
\begin{align*}\sum_{F \in \Pi[n]} (-1)^{|F|}    \sum_{\begin{subarray}{c} S \subset \{1,\dots,|F|\} \end{subarray}}  (|S|-1)! (|S^c|-1)!\  \bigg| \sum_{ i \in \underset{j\in S}{\bigcup} F_j} x_i \bigg| 
&=  \sum_{A \subset \{1,\dots,n\}}  \bigg| \sum_{ i \in A} x_i \bigg| 
\sum_{ \pi \in \Pi(A)} \mu(\pi)  \sum_{ \tilde\pi \in \Pi(A^c)} \mu(\tilde\pi)  \ .
 \end{align*}
By formulae (\ref{Soshnikov_1}) and (\ref{Mobius_0}), we conclude that
\begin{align*}
 \sum_{\sigma\in\Sy(n)} \G_n(x_{\sigma(1)},\dots, x_{\sigma(n)}) &=
 \frac{n!}{4} \sum_{A \subset \{1,\dots,n\}}  \bigg| \sum_{ i \in A} x_i \bigg|  \1_{|A|=1} \1_{|A^c|=1} \\
 &= \frac{1}{2} \delta_{n,2} \big( |x_1|+|x_2| \big) \ .
\end{align*}
Since $|x_1|=|x_2|$ when $x_1+x_2=0$, this completes the proof of  formula (\ref{G''}). \qed\\

\section{The covariance structure of OPEs} \label{A:Chebyshev}

In this section, we deal with an OPE in the one-cut regime when the Jacobi matrix $\J(\mu_N)$ has a right-limit $\L$ which is the adjacency matrix of the graph (rooted at $0$)
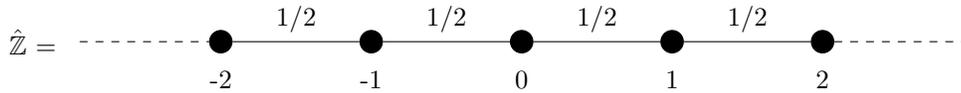
\begin{figure}[h]
\begin{center}
\begin{tikzpicture}
 \node at (-.5,0) {$\hat{\Z}$\ =} ;
 \foreach \x in {1,...,5}{\node(\x) at (2*\x,0) [draw,circle,inner sep=3pt ,fill]  {} ; }
 \node(0) at (0,0) {} ;
 \node(6) at (12,0) {} ;
 \foreach \x in {-2,...,2}{ \node at (2*\x+6,-.5) {\x} ; }
\path
   (0)      edge[dashed]                 node {}    (1)
  (1)      edge                 node [above]  {$1/2$}     (2)
  (2)      edge                 node [above]  {$1/2$}     (3)
  (3)      edge                 node [above]  {$1/2$}     (4)
  (4)      edge                 node [above]  {$1/2$}     (5)
  (5)      edge[dashed]                 node   {}     (6); 
  \end{tikzpicture}.
\end{center}
\end{figure}

As discussed in section~\ref{sect:OPE}, it corresponds to the case where the equilibrium measure for the point process $\Xi_N$  satisfies $\supp(\mu_*)= [-1,1]$ and, by theorem~\ref{thm:OPE},  the fluctuations are described by the  Gaussian field $\mathfrak{G}$, (\ref{noise}).  The aim is to prove that the OPs with respect to the spectral measure  of the  graph $\hat{\Z}$  are the basis which diagonalizes the covariance matrix of the process $\mathfrak{G}(\mathrm{X})=(\mathfrak{G}(x),\mathfrak{G}(x^2), \dots)^T$ where  $\mathrm{X}=(x^1,x^2,x^3, \dots)^T$. 
The results of section~\ref{sect:fluctuation} implies that for any polynomial $F\in\R\langle x\rangle$, the linear statistic $\Xi_N(F)$, once centered, converges in distribution to a mean-zero Gaussian random variable $\mathfrak{G}(F)$ with variance
\begin{equation*}
 \var\big[ \mathfrak{G}(F)\big]= \sum_{k>0} k \beta_{k}\beta_{-k}  \ ,
\end{equation*}
 where $  \beta_k$ are the coefficients of the Laurent polynomial $F(\frac{z+1/z}{2})$. % Let $n,m\in \N$ and  $F(x)= x^n + t x^m$, keeping track of the coefficient of $t$, we obtain
 By polarization, we have
\begin{equation}    \label{variance_2}
\E{\mathfrak{G}(x^n) \mathfrak{G}(x^m)} =  \sum_{k>0} k \alpha^n_{k}\alpha^m_{-k} \ ,
\end{equation}
 where $\alpha^n_k$ denotes  the coefficients of the Laurent polynomial $(\frac{z+1/z}{2})^n$. Using the Binomial formula, we obtain
$$\alpha^n_k = 2^{-n} {n \choose  \frac{n-k}{2}}  \ ,  $$
and by convention ${n \choose j}=0$ if $j \notin \{0,1,\dots, n\}$. These coefficients may also be obtained from formula (\ref{weight_0}) by counting the number of paths of length $n$ on the graph $\hat{\Z}$ which starts at $0$ and  end up at  the vertex $k$.  Observe that for any $k>0$, we have $ \alpha^n_{k}=  \alpha^n_{-k}$, so that if we let $\mathbf{A} = \big( \alpha^n_k \big)_{k,n \in \N}$ and ${\bf \Delta}= \big( \sqrt{n} \delta_{k,n} \big)_{k,n \in \N}$,  by formula (\ref{variance_2}), we obtain 
\begin{equation}\label{variance_3} \E{\mathfrak{G}(\mathrm{X}) \mathfrak{G}(\mathrm{X})^T} = \mathbf{A}{\bf \Delta}^2 \mathbf{A}^T \ . 
\end{equation}
The matrix $\mathbf{A}$ is lower triangular with positive entries on its main diagonal. Thus, all its principle minors are invertible and, for any $\mathrm{c} \in \R^{\N}$, we can define 
\begin{equation} \label{eigenbasis}
\mathrm{Y} = \mathbf{A}^{-1}(\mathrm{X}- \mathrm{c})  \ .
\end{equation}
In fact, we will choose $\mathrm{c}_n = 2^{-n}\alpha^n_{0}= {n \choose n/2}$ but it is not relevant for now. Then, for any $n\in\N$,  $\mathrm{Y}_n$ is polynomial of degree $n$, so that $\mathrm{Y}$ is a basis of $\R\langle x\rangle$. Moreover,  by linearity, $\mathfrak{G}(\mathrm{Y})= \mathbf{A}^{-1}\mathfrak{G}(\mathrm{X})$ since $\mathfrak{G}(\mathrm{c})=0$ because $\mathrm{c}$ is constant. According to formula (\ref{variance_3}), we get
$$ \E{\mathfrak{G}(\mathrm{Y}) \mathfrak{G}(\mathrm{Y})^T} =  {\bf \Delta}^2 \ .
$$
Hence, in the basis $\mathrm{Y}$, the covariance matrix of the process $\mathfrak{G}$ is diagonal.
It remains to check that $\mathrm{Y}(x)$ is the orthogonal basis for the spectral measure $\hat{\nu}$ of the graph $\hat{\Z}$.  By definition,  for any $n, m \ge 0$,
$$ \int x^n d\hat{\nu} = 2^{-n} \#\{\text{paths of length }n\text{ on }\hat{\Z} \text{ starting at 0 and ending at 0}\} = \alpha_0^n \ ,  $$
and 
\begin{equation} \label{nu_moments}
\int (x^n- \alpha^n_0)(x^m-\alpha^m_0) d\hat{\nu}(x) 
= \alpha^{n+m}_0 - \alpha^m_0\alpha^m_0 
\end{equation}
Then, note that 
$$  \sum_{k \in \Z}  \alpha^n_{k}\alpha^m_{-k} = 2^{-n-m} \sum_{k>0} {n \choose  \frac{n-k}{2}}  {m \choose  \frac{m+k}{2}}  = 2^{-n-m}  {n+m \choose  \frac{n+m}{2}} 
= \alpha^{n+m}_0   $$
and by formula \eqref{nu_moments}, since  $ \alpha^m_{k}=  \alpha^m_{-k}$, we obtain
\begin{align*} \int (x^n- \alpha^n_0)(x^m-\alpha^m_0) d\hat{\nu}(x) &= \sum_{k\in \Z, k\neq 0}  \alpha^n_{k}\alpha^m_{-k}  \\
&=2 \big(\mathbf{A} \mathbf{A}^T\big)_{nm}
\end{align*}
By (\ref{eigenbasis}), this implies that
\begin{align*} \int \mathrm{Y} \mathrm{Y}^T  d\hat{\nu} 
&= \mathbf{A}^{-1}\left(  \int (\mathrm{X}- \mathrm{c}) (\mathrm{X}- \mathrm{c})^T d\hat{\nu}  \right) \mathbf{A}^{-T}  \\
&=  2 \mathbf{I}
 \end{align*}
where $\mathbf{I}$ is the identity matrix.   The last identity shows that  $\mathrm{Y}$ is an eigenbasis for $\R\langle x\rangle$ equipped with the inner product inherited from $L^2(\hat{\nu})$. It is a well-know fact that the spectral measure of 
the graph $\hat{\Z}$ is the arcsine distribution $d\hat{\nu}= \frac{1}{\pi}(1-x^2)^{-1/2} \1_{|x|<1} dx$, so that  $\mathrm{Y}_k$ are proportional to the Chebyshev polynomials of the first kind. Actually, if $\mathrm{Y}_0\equiv1$, formula (\ref{eigenbasis}) boils down to the well-known identity, 
\begin{align*}
x^n  =  \sum_{k \ge 0} \alpha^n_{k} \mathrm{Y}_k(x) + \mathrm{c}_n  = 2^{1-n}\sum_{k=1}^n {n \choose \frac{n- k}{2} } T_{k}(x)   
 +2^{-n} {n \choose n/2} \ . \end{align*}

%\bibliography{References}
%\bibliographystyle{acm}

\end{document}